\documentclass[10pt]{article}

\usepackage[T1]{fontenc}
\usepackage[utf8]{inputenc}
\usepackage{lmodern}
\usepackage{graphicx}
\usepackage{subcaption}
\usepackage{amsmath,amssymb,amsfonts,amsthm}
\usepackage{dsfont}
\usepackage{color}
\usepackage{hyperref}
\usepackage{tabularx}
\usepackage{etoolbox}
\usepackage{latexsym}
\usepackage{booktabs}
\usepackage{chngcntr}
\usepackage{cancel}

\counterwithin*{equation}{section}

\usepackage{algorithm}
\usepackage{algorithmic}

\usepackage{tikz}

\usepackage[backend=biber, style=numeric, maxcitenames=99, maxnames=99]{biblatex}
\newtheorem{Theorem}{Theorem}[part]
\newtheorem{Definition}{Definition}[part]
\newtheorem{Proposition}{Proposition}[part]

\newtheorem{Lemma}{Lemma}[part]
\newtheorem{Corollary}{Corollary}[part]
\newtheorem{Remark}{Remark}[part]

\newcommand{\Sum}[2]{\sum\limits_{#1}^{#2}}

\newcommand{\Int}[2]{\int_{#1}^{#2}}
\newcommand{\Frac}[2]{\frac{#1}{#2}}

\newcommand{\Ind}[1]{\mathds{1}_{\left\{ #1 \right\}}}
\newcommand{\Expect}[1]{\mathbb{E}\left[ #1 \right]}
\newcommand{\Prob}[1]{\mathbb{P}\left( #1 \right)}

\newcommand{\ubar}[1]{\bar{#1}}
\newcommand{\lbar}[1]{\underline{#1}}

\newcommand{\N}{\mathbb{N}}
\newcommand{\R}{\mathbb{R}}

\newcommand{\Z}{\mathbb{Z}}
\newcommand{\E}{\mathbb{E}}

\newcommand{\F}{\mathbb{F}}

\newcommand{\T}{\mathbb{T}}

\newcommand{\Ac}{\mathcal{A}}

\newcommand{\Cc}{\mathcal{C}}

\newcommand{\Fc}{\mathcal{F}}

\newcommand{\Lc}{\mathcal{L}}

\newcommand{\Oc}{\mathcal{O}}

\newcommand{\Sc}{\mathcal{S}}

\def \ni{\noindent}

\newcommand{\Dx}[1]{\Frac{\partial #1}{\partial x}}

\newcommand{\Dy}[1]{\Frac{\partial #1}{\partial y}}

\usepackage{geometry}
\title{Optimal dividend and capital injection under self-exciting claims}

\author{
	Paulin Aubert\thanks{Laboratoire de Mathématiques et Modélisation d'Évry, Université Évry Paris-Saclay, Exiom Partners, France, \sf paulin.aubert@univ-evry.fr.} $\:$ 
	Etienne Chevalier\thanks{Laboratoire de Mathématiques et Modélisation d'Évry, Université Évry Paris-Saclay, UMR 8071 CNRS, France, \sf etienne.chevalier@univ-evry.fr.} $\:$ 
	Vathana Ly Vath\thanks{Laboratoire de Mathématiques et Modélisation d'Évry, Université Paris-Saclay, ENSIIE, UMR 8071 CNRS, France, \sf vathana.lyvath@ensiie.fr.}
}

\begin{document}

\maketitle

\begin{abstract} 
	In this paper, we study an optimal dividend and capital-injection problem in a Cramér-Lundberg model where claim arrivals follow a Hawkes process, capturing clustering effects often observed in insurance portfolios.
	We establish key analytical properties of the value function and characterise the optimal capital-injection strategy through an explicit threshold. We also show that the value function is the unique viscosity solution of the associated HJB variational inequality.
	For numerical purposes, we first compute a benchmark solution via a monotone finite-difference scheme with Howard’s policy iteration. We then develop a reinforcement learning approach based on policy-gradient and actor-critic methods. The learned strategies closely match the PDE benchmark and remain stable across initial conditions.
	The results highlight the relevance of policy-gradient techniques for dividend optimisation under self-exciting claim dynamics and point toward scalable methods for higher-dimensional extensions.
\end{abstract}

\vspace{7mm}

\par \bigskip


\section{Introduction} 
\label{SEC:introduction}

	The allocation of an insurer’s surplus between solvency and shareholder remuneration is a central question in actuarial science, traditionally addressed through ruin probabilities and optimal dividend policies. 
	In this context, the surplus process provides a framework for quantifying the trade-off between long-term financial stability and the distribution of profits. Since the seminal contribution of de~Finetti~\cite{DeFinetti1957}, a vast literature has emerged at the intersection of probability theory, stochastic control, and insurance mathematics. \\
	
	Classical studies build upon the Cramér--Lundberg model introduced by Lundberg~\cite{Lundberg1903} and Cramér~\cite{Cramer1930}, and further developed by Gerber~\cite{Gerber1969, Gerber1970}. Over the past decades, the dividend optimisation problem has been analysed using both regular and singular control techniques in models driven by compound Poisson processes or Brownian motion. See for instance Jeanblanc and Shiryaev~\cite{Jeanblanc_Shiryaev_1995}, Asmussen and Taksar~\cite{ASMUSSEN19971}, and Gerber and Shiu~\cite{GerberShiu2004, Gerber2006}.
	Numerous extensions have since been proposed to incorporate investment risk, reinsurance, capital injections, and taxation, as documented in the works of Paulsen and Gjessing~\cite{PAULSEN1997215}, Hojgaard and Taksar~\cite{HOJGAARD199841}, Azcue and Muler~\cite{Azcue2005}, Kulenko and Schmidli~\cite{KULENKO2008270}, Lokka and Zervos~\cite{LOKKA2008954}, and Albrecher and Thonhauser~\cite{ALBRECHER2008134}. 
	Comprehensive reviews of these developments can be found in Albrecher and Thonhauser~\cite{Albrecher2009} and in the monograph by Schmidli~\cite{Schmidli2008}. 
	A persistent assumption in the classical literature is that claim arrivals are independent and identically distributed, typically modelled by a Poisson process. Yet real insurance portfolios---particularly those exposed to catastrophic, environmental, cyber, or systemic risks---often display pronounced clustering, generating temporal dependence in claim occurrences. 
	This has motivated the use of more general point processes, including Cox and shot-noise dynamics~\cite{albrecher2023optimal}, and Hawkes processes, whose relevance for insurance risk modelling was demonstrated early on by Stabile and Torrisi~\cite{StabileTorrisi2010} in the context of ruin probabilities.
	More broadly, Hawkes processes have found a growing range of applications in insurance and risk management. Bensoussan and Chevalier-Roignant~\cite{BensoussanCR2024} provide an overview of stochastic control problems for diffusions with self-exciting jumps, highlighting the analytical challenges specific to this class of models. In the context of cyber risk, Hillairet, Réveillac, and Rosenbaum~\cite{HillairetRR2023} develop an expansion formula for Hawkes processes with applications to cyber-insurance derivatives, while Ren and Zhang~\cite{RenZhang2025} propose a Hawkes-based model for cyber risk insurance incorporating loss covariates. On the control side, Cao, Landriault, and Li~\cite{CaoLandriaultLi2020} study optimal reinsurance-investment strategies under a dynamic contagion claim model, and Brachetta, Callegaro, Ceci, and Sgarra~\cite{Brachetta2024} address optimal reinsurance via BSDEs in a partially observable setting with jump clusters. These works illustrate the relevance of self-exciting dynamics for insurance modelling and control but do not address the joint optimisation of dividends and capital injections, which is the focus of the present paper.
	Dividend optimisation has been analysed in some of these non-Poisson settings. However, the combined optimisation of dividends and capital injections in a Cramér-Lundberg model driven by Hawkes claim arrivals has not been addressed in the existing literature.
	The present work develops a dividend optimisation framework for a Cramér–Lundberg model with Hawkes claim arrivals, allowing for both dividend distributions and capital injections.
	This extends classical results obtained under compound Poisson dynamics, including those of Kulenko and Schmidli~\cite{KULENKO2008270}. 
	From a technical standpoint, the Hawkes dynamics introduce difficulties absent from these earlier works: the jump operator couples both state variables, preventing the concavity arguments that underpin the one-dimensional analyses of~\cite{KULENKO2008270} and~\cite{ALBRECHER2008134}. Moreover, unlike~\cite{KULENKO2008270} where capital injections systematically prevent ruin, our model allows for bankruptcy, giving rise to an additional free boundary.
	From an analytical perspective, we establish fundamental properties of the value function, including bounds, monotonicity, and local Lipschitz continuity, and characterise the optimal capital-injection strategy through an explicit threshold. 
	We then show that the value function is the unique viscosity solution of the associated Hamilton--Jacobi--Bellman variational inequality.
	
	Because Hawkes dynamics considerably increase the analytical complexity of the model, numerical methods are required to approximate the value function and the associated optimal policy. 
	As a classical benchmark, we first compute a reference solution using a monotone finite-difference approximation of the HJB variational inequality combined with Howard’s policy iteration algorithm. This PDE-based approach serves to validate the structure of the optimal strategy in our setting.
	The main numerical contribution of the paper lies in the development of a reinforcement learning methodology tailored to this class of singular stochastic control problems. 
	A growing body of work aims to connect stochastic control theory with reinforcement learning by developing policy-gradient and actor–critic formulations, as illustrated by the contributions of Wang et al.~\cite{Wang2020}, Jia and Zhou~\cite{JiaZhou22_AC, JiaZhou22_TD, JiaZhou23_QL}, as well as the recent advances of Hamdouche et al.~\cite{Hamdouche23} and Pham and Warin~\cite{PhamHuyên2023Acla}.
	Building on these developments, we examine whether parameterised stochastic policies can learn near-optimal dividend and capital-injection strategies in our setting.
	Our methodology is related to the framework of Hamdouche et al.~\cite{Hamdouche23}, who study policy-gradient approaches for control problems with random exit times. 
	The results obtained reinforce the view that policy-gradient algorithms offer a scalable alternative to PDE-based methods, and can be effectively applied to higher-dimensional or path-dependent extensions of the dividend optimisation problem where classical numerical techniques become impractical. \\
	
	The remainder of the paper is structured as follows. 
	Section~\ref{SEC:model_description} introduces the surplus model with capital injections and Hawkes-driven claims. 
	Section~\ref{SEC:theorical_analysis} establishes key analytical properties of the value function and derives the structure of the optimal injection strategy. 
	Section~\ref{SEC:viscosity_solution} shows that the value function is the unique viscosity solution to the associated HJB variational inequality. 
	Section~\ref{SEC:pde_resolution} presents the finite-difference framework and the corresponding numerical results, which serve as a benchmark and illustrate the economic features of the optimal policy. 
	Finally, Section~\ref{SEC:rl_resolution} develops the reinforcement learning methodology and compares the learned strategies with the PDE benchmark.
		
\section{Modelling insurer's portfolio and clustering effect}
\label{SEC:model_description}
	
	\subsection{Uncontrolled surplus dynamics}
	
		Let $(\Omega, \Fc, \mathbb{P})$ be a probability space on which all stochastic processes and random variables are defined and such that $\mathbb F$ is complete and right-continuous. The insurer's cash reserve is represented by a stochastic process $R = (R_t)_{t \geq 0}$, whose dynamics, in the absence of any control, follow the classical Cramér-Lundberg model:
		\begin{equation*}
			R^x_t = x + c t - \Sum{k=1}{N_t} Y_k ,
		\end{equation*}
		where $x \in \R^+$ is the company's initial capital, $c > 0$ is the constant premium income per unit of time, $N = (N_t)_{t \geq 0}$ is a counting process representing the number of claims occurring up to time $t$ and $(Y_k)_{k \in \N}$ is a positive random variable with density $f$, assumed independent of the counting process. 
		
		Traditionally, $N$ is assumed to be a homogeneous Poisson process with constant intensity $\lambda > 0$, which implies independent, exponentially distributed inter-arrival times. While analytically convenient, this framework is not designed to account for temporal dependence in claim arrivals, which motivates the use of more flexible models such as self-exciting processes.
	
	\subsection{Temporal claim dependence via a Hawkes process}
	
		In practice, claim arrivals often exhibit temporal clustering: events such as natural disasters, cyber incidents, or pandemics tend to generate multiple claims in short time intervals. This behaviour, known as the clustering effect, contradicts the memoryless nature of the Poisson process.
		
		To model the clustering behaviour of claims, we choose to represent the arrival process $N$ as a Hawkes process. Hawkes processes are well known for their ability to model clustering effects. In our framework, the claim arrival intensity $\lambda = (\lambda_t)_{t \geq 0}$ evolves dynamically according to the following equation:
		\begin{equation*}
			\lambda_t = a(b - \lambda_t) dt + \eta dN_t,
		\end{equation*}
		where $a, b, \eta > 0$ are model parameters, and the initial condition is $\lambda_0 = y \in [b, +\infty)$. Between claim arrivals, the intensity $\lambda_t$ reverts toward the long-term level $b$ at rate $a$, while each claim at time $t$ increases $\lambda_t$ by $\eta$. This dynamics captures both the self-exciting nature and the memory effects in claim arrivals. 
		We assume $\lambda_0 \geq b$ without loss of generality. Indeed, under exponential kernels and as soon as a few claims occur, the intensity will almost surely exceed $b$ and remain above it due to the accumulation of excitation. This assumption also simplifies several technical arguments in the analysis that follows.
	
		\begin{Remark}[Net profit condition]
			\label{RMK:net_profit_condition}
			Under the stationarity condition $\eta / a < 1$, the claim intensity admits a finite stationary mean $\lambda_\infty = b/(1 - \eta/a)$. The net profit condition for the uncontrolled surplus process reads
			\begin{equation}
				\label{EQN:net_profit_condition}
				c > \frac{b \, \mathbb{E}[Y]}{1 - \eta/a}.
			\end{equation}
			When~\eqref{EQN:net_profit_condition} holds, the probability of ruin is strictly less than one for all positive initial capitals and decays exponentially, with a Cramér--Lundberg bound (see~\cite{StabileTorrisi2010}, Theorem~4.1). When it fails, or when $\eta/a \geq 1$, ruin is certain.
		\end{Remark}
	
	\subsection{Controlled surplus dynamics}
	
		We assume that the company is owned by a group of shareholders whose objective is to extract value from the surplus through dividend distributions, while preserving solvency via capital injections when needed. These two financial levers modify the surplus dynamics, leading to a controlled process.
		
		Let $\alpha = (Z_t, K_t)_{t \geq 0}$ be a control strategy, where $Z$ is a non-decreasing, right-continuous, $\Fc$-adapted process representing the cumulative dividends paid out to shareholders and $K$ is a non-decreasing, left-continuous, $\Fc$-adapted process representing the cumulative capital injections by shareholders.
		Under strategy $\alpha$ the controlled surplus process is given by:
		\begin{align*}
			X_t &= R^x_t - Z_t + K_t \\
			&= x + c t - \Sum{k=1}{N_t} Y_k - Z_t + K_t .
		\end{align*}
		Dividend payments reduce the reserve, while capital injections increase it. These interventions are subject to economic constraints and are only permitted within an admissible set.
		To ensure both economic relevance and mathematical well-posedness of the model, we restrict our attention to a class of admissible strategies defined as follows:
		\begin{Definition}[Set of admissible strategies] A strategy $\alpha_t = (Z_t, K_t)_{t \geq 0}$ is said to be admissible if:
			\begin{itemize}
				\item $Z$ is càd-làg, $\F$-adapted, non-decreasing and such that $Z_{t} - Z_{t^-} \leq X_{t^-}+R^0_t-R^0_{t^-}$ and $Z_{0^-} = 0$,
				\item $K$ is càg-làd, $\F$-adapted, non-decreasing and such that $K_{0^-} = 0$.
			\end{itemize}
			
			When  $(X_0,\lambda_0)=(x,y)\in\mathbb R\times\lbrack b,+\infty)$, the set of admissible strategies is denoted by $\Ac(x, y)$.
		\end{Definition}
		The condition $Z_{t} - Z_{t^-} \leq X_{t^-}+R^0_t-R^0_{t^-}$ enforces that dividends cannot be paid beyond the available reserve at any time. 
		
	\subsection{Ruin and objective function}
	
		As is standard in risk theory, we assume that the company ceases operations at the time of ruin, i.e., when its reserve becomes negative. The ruin time under strategy $\alpha$ is defined as:
		\begin{equation*}
			T^\alpha = \inf \{ t \geq 0, X_{t^+} < 0 \}.
		\end{equation*}
		The shareholders' objective is to maximize the expected discounted net gains until ruin. The gain includes the total discounted dividends and subtracts a penalty proportional to the capital injected. Formally, for $(x,y)\in\mathbb R\times\lbrack b,+\infty)$, the reward associated with a strategy $\alpha \in \Ac(x, y)$ is given by:
		\begin{equation*}
			J_\alpha(x, y) = \Expect{\Int{0}{T^\alpha} e^{-\rho s} dZ_s - \delta \Int{0}{T^\alpha} e^{-\rho s} dK_s},
		\end{equation*}
		where $\rho > 0$ is the discount rate, and $\delta > 1$ is the penalty coefficient reflecting the opportunity cost of capital injections.
		The optimisation problem then consists in maximising $J_\alpha(x, y)$ over all admissible strategies:
		\begin{equation}
			\label{EQN:optimal_control_problem}
			v(x, y) = \underset{\alpha \in \Ac(x, y)}{\sup} J_\alpha(x, y)\quad \textrm{on }\mathbb R\times\lbrack b,+\infty).
		\end{equation}
		
		\begin{Remark}
			The condition $\delta > 1$ is crucial to prevent excessive capital injections, which would otherwise be incentivised if $\delta \leq 1$. Similarly, the discount rate $\rho > 0$ ensures finiteness of the value function and rules out infinite accumulation of dividends over time. See \cite{KULENKO2008270} for a detailed discussion.
		\end{Remark}
	
\section{Theoretical analysis}
\label{SEC:theorical_analysis}

	We now examine the analytical and structural properties of the value function associated with the stochastic control problem~\eqref{EQN:optimal_control_problem}. 
	These results provide the mathematical foundations required to characterize the value function as a viscosity solution of the Hamilton--Jacobi--Bellman (HJB) equation, a task carried out in Section~\ref{SEC:viscosity_solution}.
	
	\subsection{Pre-claim intensity}
	
		A recurring element in our analysis is the conditional behaviour of the claim intensity process prior to the first jump of the counting process $N$. In order to simplify computations involving the law of the first claim time, we introduce the deterministic intensity process $\tilde{\lambda}$, defined on the event $\{ t \leq \tau_1 \}$, where $\tau_1$ denotes the first jump time of $N$.
		\ni On $\{ t \leq \tau_1 \}$, the intensity process satisfies the deterministic ordinary differential equation:
		\begin{equation*}
			d {\lambda}_t = a(b - {\lambda}_t)   dt, \quad \lambda_0=y\geq b ,
		\end{equation*}
		which integrates explicitly to:
		\begin{equation}
			\label{EQN:tilde_lambda}
			\tilde{\lambda}_t :=\lambda_t= b - (b - y) e^{- a t} .
		\end{equation}
		This expression appears frequently in computations involving expectations conditional on the absence of claims. In particular, the following expression for the survival probability will be used repeatedly in the analysis. Let $h > 0$. The probability that no claim occurs up to time $h$ is given by:
		\begin{align*}
			\Prob{\tau_1 \geq h} &= e^{-\Int{0}{h} \tilde{\lambda}_s ds} \nonumber \\ 
			&= e^{-bh - \frac{y - b}{a} \left( 1 - e^{- ah} \right)} \nonumber \\ 
			&\underset{h \rightarrow 0}{=} 1 - yh + o(h) .
		\end{align*}
		This approximation is especially useful when analysing the infinitesimal behaviour of the controlled process, as will be required in the rest of this section.
	
	\subsection{Dynamic programming principle}
		
		We begin by establishing the dynamic programming principle (DPP) associated with the control problem~\eqref{EQN:optimal_control_problem}. This fundamental result expresses the value function in terms of sequentially optimal decisions over time intervals, and serves as the cornerstone for deriving the Hamilton--Jacobi--Bellman equation and for analysing the structural properties of the value function.
		In our problem, the dynamic programming principle can be stated as follows:
		
		\begin{Proposition}[Dynamic Programming Principle] 
			\label{PROP:dynamic_programming_princple}
			Let $\theta$ be any $\mathcal F$-stopping time and $(x,y)\in\mathbb R\times\lbrack b,+\infty)$, it follows from the dynamic programming principle that
			\begin{equation}
				\label{EQN:dynamic_programming_principle}
				v(x, y) = \underset{\alpha \in \Ac(x, y)}{\sup} \Expect{\Int{0}{T^\alpha \wedge \theta} e^{-\rho s} dZ_s - \delta \Int{0}{T^\alpha \wedge \theta} e^{-\rho s} dK_s + e^{-\rho (T^\alpha \wedge \theta)}v(X_{T^\alpha \wedge \theta}, \lambda_{T^\alpha \wedge \theta})} .
			\end{equation}
		\end{Proposition}
		We refer to standard texts (e.g., \cite{KULENKO2008270}) for further details and omit the proof, which follows classical arguments.
		
	\subsection{Analytical properties of the value function}
	
		We first derive upper and lower bounds for the value function. The lower bound is immediate, as the controller can always choose to take no action. The upper bound corresponds to an idealised scenario where all available surplus is instantly paid as dividends without receiving further claims.
		
		\begin{Proposition}[Value function boundaries]
			\label{PROP:value_function_bounds}
			For $x \in \R$ and $y \geq b$, we have
			\begin{equation*}
				x^+ \leq v(x, y) \leq x^+ + \frac{c}{\rho} .
			\end{equation*}
		\end{Proposition}
		
		\begin{proof}
			The lower bound follows by considering a strategy $\hat{\alpha}$ which immediately distributes the whole cash reserve, $x^+$ and then does not distribute any dividends and does not inject capital. We get
			\begin{equation*}
				v(x, y) \geq J_{\hat{\alpha}}(x, y) = x^+,\quad\textrm{for }(x,y)\in\mathbb R\times\lbrack b,+\infty) .
			\end{equation*}
			Let $(x,y)\in\mathbb R^+\times\lbrack b,+\infty)$. We know that, for any strategy $\alpha=(Z,K)$, we have $0\leq Z_u \leq x + cu+K_u$, on before $\{u\leq T^\alpha\}$,  so we deduce that:
			\begin{equation*}
				J_\alpha(x, y) \leq x + \Expect{\Int{0}{T^\alpha} e^{-\rho s} c ds+(1-\delta)\Int{0}{T^\alpha} e^{-\rho s} dKs} \leq  x + \Int{0}{+\infty} e^{-\rho s} c ds = x + \frac{c}{\rho} .
			\end{equation*}
			If $(x,y)\in\mathbb R^-\times\lbrack b,+\infty)$, there are only two admissible actions at time 0: letting the firm going to bankruptcy or injecting capital up to 0. Hence we have
			\begin{equation*}
				v(x,y)\leq \max(x+v(0,y); 0)\leq \frac{c}{\rho}.
			\end{equation*}
		\end{proof}
		
		\ni
		We next establish monotonicity properties, reflecting the natural intuition that higher surplus enhances value, whereas higher claim intensity reduces it.
		
		\begin{Proposition}[Monotonicity in $x$]
			\label{PROP:increase_of_v}
			Let $0 \leq x < x'$ and $y \in \lbrack b,+\infty)$. Then:
			\begin{equation*}
				v(x', y) - v(x, y) \geq x' - x .
			\end{equation*}
		\end{Proposition}
		
		\begin{proof} Let $\varepsilon > 0$ and $\alpha_\varepsilon$ be an $\varepsilon$-suboptimal strategy, i.e., $J_{\alpha_\varepsilon}(x, y) \geq v(x, y) - \varepsilon$. Let $0 \leq x < x'$. We consider the strategy consisting in distributing dividends up to $x$ and then apply strategy $\alpha_\varepsilon$. By the dynamic programming principle \eqref{EQN:dynamic_programming_principle} we obtain:
			\begin{align*}
				v(x', y) &\geq x' - x + J_{\alpha_\varepsilon}(x, y) \\
				&\geq x' - x + v(x, y) - \varepsilon .
			\end{align*}
			Letting $\varepsilon \to 0$ yields the desired result.
		\end{proof}
		
		\begin{Proposition}[Monotonicity in $y$]
			\label{PROP:value_function_monotony}
			Let $x \in \R^+$. The function $y\to v(x,y)$ is non-increasing on $\lbrack b,+\infty)$.
		\end{Proposition}
		
		\begin{proof} 
			Let $x\in\mathbb R$ and $b \leq y < y'$. For $\varepsilon > 0$, there exists  $\varepsilon$-suboptimal strategy $\alpha^\varepsilon=(Z^\varepsilon, K^\varepsilon)\in\mathcal A(x,y')$ such that
			$v(x,y')\leq J^{\alpha^\varepsilon}(x,y')+\varepsilon$.
			
			\ni From theorem 3.2 in \cite{DawsonLi12} we know  that the intensity process is such that $\lambda^{y} \leq \lambda^{y'}$ almost surely, which, from Lemma A.2 in \cite{albrecher2023optimal} implies that $X^{\alpha^\varepsilon, y'} \leq X^{\alpha^\varepsilon, y}$ and therefore $T^{\alpha^\varepsilon, y'} \leq T^{\alpha, y}$. From the dynamic programming principle we have:
			\begin{align*}
				v(x, y) &\geq \Expect{\Int{0}{T^{\alpha^\varepsilon, y'}} e^{- \rho s} d(Z^\varepsilon_s-\delta K_s^\varepsilon) +e^{-\rho T^{\alpha^\varepsilon,y'}}v(X^x_{ T^{\alpha^\varepsilon,y'}}, \lambda^y_{ T^{\alpha^\varepsilon,y'}})}, \\
				& \geq \Expect{\Int{0}{T^{\alpha^\varepsilon, y'}} e^{- \rho s} d(Z^\varepsilon_s-\delta K_s^\varepsilon) }\\
				&\geq v(x, y') + \varepsilon .
			\end{align*}
			As this is true for every $\varepsilon > 0$ we deduce, by letting $\varepsilon$ going to 0, that $v$ is non increasing in $y$.
		\end{proof}
		
		\begin{Corollary}
			\label{COR:yasymptotic}
			For $x\in\mathbb R^+$, we have
			$$\lim_{y\to +\infty}v(x,y)=x.$$
		\end{Corollary}
		
		\begin{proof}
			Let $x\in\mathbb R^+$. For $y> e^b$, we set $t^*(y)=\frac{1}{a}\ln\left(\frac{y-b}{\ln(y)-b}\right)$. $t^*(y)$ is then such that 
			$$\lambda^y_t\geq\tilde \lambda^y_{t} \geq \tilde \lambda^y_{t^*(y)},\quad\textrm{for all }t\leq t^*(y).$$
			Let $\varepsilon>0$, it follows from the dynamic programming principle \eqref{PROP:dynamic_programming_princple}, that 
			\begin{eqnarray*}
				v(x,y) & \leq & \varepsilon+\Expect{\Int{0}{T^\alpha \wedge t^*(y)} e^{-\rho s} d(Z_s - \delta K_s) + e^{-\rho (T^\alpha \wedge t^*(y))}v(X_{T^\alpha \wedge t^*(y)}, \lambda_{T^\alpha \wedge t^*(y)})}	
			\end{eqnarray*}
			As $\lim_{y\to+\infty}t^*(y)=+\infty$ and from the monotonicity of $v$ in intensity, for $y$ big enough, we have
			$$\Expect{\Int{0}{T^\alpha \wedge t^*(y)} e^{-\rho s} d(Z_s - \delta K_s)1_{\{T^\alpha\leq t^*(y)\}}}\leq v^{\ln(y)}(x),$$
			where we denote by $v^\zeta$ the value function of our control problem with a constant intensity equal to $\zeta$. 	We recall that $\lim_{\zeta\to+\infty}v^{\zeta}(x)=x.$ On the other hand, from Propositions \eqref{PROP:value_function_bounds}, \eqref{PROP:increase_of_v} and  \eqref{PROP:value_function_monotony}, we have
			\begin{eqnarray*}
				\Expect{\left(\Int{0}{t^*(y)} e^{-\rho s} d(Z_s - \delta K_s) + e^{-\rho  t^*(y)}v(X_{ t^*(y)}, \lambda_{ t^*(y)})\right)1_{\{T^\alpha> t^*(y)\}}} & \leq & e^{-\rho t^*(y)} \left(x+ct^*(y)+\frac{c}{\rho}\right).
			\end{eqnarray*}
			Hence we can conclude, thanks to Proposition  \eqref{PROP:value_function_bounds}, that
			$$x\leq \lim_{y\to+\infty}v(x,y)\leq  \varepsilon+\lim_{y\to+\infty}	v^{\ln(y)}(x)+e^{-\rho t^*(y)}\left(x+ct^*(y)+\frac{c}{\rho}\right)=\varepsilon+x$$
			We obtain the result by letting $\varepsilon$ going to 0.
		\end{proof}
		
		We establish the local Lipschitz continuity of the value function, which is essential for applying comparison principles in the viscosity solution analysis.
		
		\begin{Proposition}[Local Lipschitz continuity]
			\label{LEM:lipscontinuity}
			The value function $v$ is locally Lipschitz continuous on $\mathbb R\times(b, +\infty)$. More precisely:
			
			\begin{itemize}
				\item[i)] For $b\leq y$ and $0 \leq x < x'$ and $y > 0$, we have:
				\begin{equation*}
					x'-x\leq v(x', y) - v(x, y) \leq \delta(x'-x) .
				\end{equation*}
				If $x<x'\leq 0$, we have
				\begin{equation*}
					0\leq v(x',y)-v(x,y)= \textrm{max}(\delta x'+v(0,y),0)-\textrm{max}(\delta x+v(0,y),0)\leq \delta(x'-x) .
				\end{equation*}
				\item[ii)] Let $x \in \R^+$ and $b < y < y'$. For $\varepsilon>0$ such that $b+\varepsilon\leq y$, we have:
				\begin{equation*}
					0\leq v(x, y) - v(x, y') \leq (x+\frac{c}{\rho}) \frac{a}{\varepsilon}(y' - y) +o(y'-y) .
				\end{equation*}
				If $x\leq 0$, we have $v(x,y)-v(x,y')\leq \textrm{max}(\delta x+v(0,y),0)-\textrm{max}(\delta x+v(0,y'),0)\leq v(0,y)-v(0,y') \leq  \frac{ac}{\rho\varepsilon}(y' - y) +o(y'-y).$
			\end{itemize}
		\end{Proposition}
		
		\begin{proof} 
			We start by showing the first point:
			\begin{itemize}
				\item[i)] Let $y \geq b$ and  $0 \leq x < x'$, we consider the strategy, in $\mathcal A(x,y)$ which consits in immediately inject some capital up to the cash level $x'$. It follows from the dynamic programming principle that:
				\begin{equation*}
					v(x,y)\geq v(x',y)-\delta(x'-x) .
				\end{equation*}
				Set $t_0 = (x' - x)/c$ and
				\begin{align*}
					v(x, y) &\geq \Expect{e^{-\rho (\tau_1 \wedge t_0)} v(X_{\tau_1 \wedge t_0}, \lambda_{\tau_1 \wedge t_0})} \\
					&\geq \Expect{e^{-\rho t_0} v(X_{t_0}, \lambda_{t_0}) \Ind{t_0 \leq \tau_1}} \\
					&\geq e^{-\rho t_0} v(x', \lambda_{t_0}) \Prob{t_0 \leq \tau_1} \\
					&= e^{-\rho t_0} e^{-\Int{0}{t_0} \tilde{\lambda}_s ds} v(x', y),\quad\textrm{with }\tilde{\lambda}_s=b - (b - y) e^{- a s}.
				\end{align*}
				Hence, from Proposition \ref{PROP:increase_of_v}, we have:
				\begin{equation*}
					x'-x\leq v(x', y) - v(x, y) \leq \left( e^{\rho t_0} e^{\Int{0}{t_0} \tilde{\lambda}_s ds} - 1 \right) v(x, y)
				\end{equation*}
				As $z\to ze^z-e^z+1$ takes values in $\mathbb R^+$  and that $\Int{0}{t_0} \tilde{\lambda}_s ds = bt_0 + \frac{y-b }{a}\left( 1 - e^{-a t_0}\right) $, we get`
				\begin{eqnarray*}
					v(x', y) - v(x, y) & \leq  & \left\lbrack(\rho+b)t_0+\frac{y-b }{a}\left( 1 - e^{-a t_0} \right)\right\rbrack e^{(\rho+y) t_0} e^{\Int{0}{t_0} \tilde{\lambda}_s ds} v(x, y)\\
					& \leq & (\rho+y)t_0e^{(\rho+y) t_0} v(x, y)\\
					& \leq & \frac{\rho+y}{c}e^{\frac{\rho+y}{c}(x'-x)}(x+\frac{c}{\rho})(x'-x).
				\end{eqnarray*}
				We conclude that $v$ is locally Lipschitz in $x$ and that
				\begin{equation*}
					x'-x\leq v(x', y) - v(x, y) \leq (x+\frac{c}{\rho})\left((x' - x) \frac{\rho + y}{c}+o(x' - x)\right)
				\end{equation*}
				Hence $v$ is Lipschitz in $x$.
										
				\item[ii)] We now consider $x\in\mathbb R^+$, $b < y < y'$ and $\varepsilon>0$ such that $y\geq b+\varepsilon$ . Let $t_0$ be such that $\tilde\lambda^{y'}_{t_0} = y$. 
				From the definition of $\tilde{\lambda}$ (see \eqref{EQN:tilde_lambda}), we have that:
				\begin{equation*}
					t_0 = - \frac{1}{a} \ln \left( \frac{y - b}{y' - b} \right)\leq \frac{y'-y}{a(y-b)}\leq \frac{1}{a\varepsilon}(y'-y).
				\end{equation*}
				Applying the stragtegy $(0,0)\in\mathcal A(x,y')$, the dynamic programming principle implies that:
				\begin{align*}
					v(x, y') &\geq \Expect{ e^{-\rho (\tau_1 \wedge t_0)} v(X_{\tau_1 \wedge t_0}, \lambda_{\tau_1 \wedge t_0})} \\
					&\geq e^{-\rho t_0} v(x+c{t_0}, y) \Prob{t_0 \leq \tau_1} \\
					&\geq e^{-\rho t_0} e^{-\Int{0}{t_0} \tilde{\lambda}_s ds} v(x+c{t_0}, y) .
				\end{align*}
				As $v(x+c{t_0}, y) \geq v(x,y)$, we deduce from Proposition \ref{PROP:value_function_monotony} that
				\begin{equation*}
					0\leq v(x, y) - v(x, y') \leq v(x, y') \left(  e^{\rho t_0} e^{\Int{0}{t_0} \tilde{\lambda}_s ds} - 1 \right)\leq v(x, y')(\rho+y')t_0 e^{(\rho+y') t_0} ,
				\end{equation*}
				$v$ is then locally lipshitz in its second variable and
				\begin{equation*}
					0\leq v(x, y) - v(x, y') \leq (x+\frac{c}{\rho}) \frac{a}{\varepsilon}(y' - y)  e^{\frac{a}{\varepsilon}(\rho+y') (y'-y)} .
				\end{equation*}
			\end{itemize}
		\end{proof}
	
	\subsection{Capital injection strategies}

		From an economic perspective, capital injections represent a costly measure that should only be used to prevent imminent ruin. Injecting capital before it is strictly necessary is therefore suboptimal, as we formally establish below.
		
		\begin{Proposition}[Capital injection policy]
			\label{PROP:inject_capital}
			Injecting capital is only optimal when strictly necessary to prevent ruin. In particular, capital injection at time $t$ can only be optimal if the controlled surplus satisfies $X_t < 0$.
		\end{Proposition}
		
		\begin{proof}
			For $y\geq b$ and $x\in\mathbb R$, two scenarios must be considered: 
			\begin{itemize}
				\item[i)] If $x \geq 0$, we claim that:
				\begin{equation*}
					v(x + \varepsilon, y) - \varepsilon \delta < v(x, y)
				\end{equation*}
				Let $\kappa > 0$ and $\varepsilon > 0$ and $\alpha_\varepsilon$ a $\varepsilon$-suboptimal strategy, i.e. it is such that:
				\begin{equation*}
					v(x + \kappa, y) \leq J_{\alpha_\varepsilon}(x + \kappa, y) + \varepsilon .
				\end{equation*}
				We let $\hat{\alpha}$ be the strategy consiting in applying $\alpha_\varepsilon$ while $s < T^{\alpha_\varepsilon}$ and increasing capital if $s = T^{\alpha_\varepsilon}$. Then we have: 
				\begin{align*}
					v(x, y) &\geq J_{\hat{\alpha}}(x, y) \\
					&= \Expect{\Int{0}{T^{\alpha_\varepsilon}} e^{-\rho s} \left( dZ_s - \delta dK_s \right) - e^{-\rho T^{\alpha_\varepsilon}} \kappa \delta + e^{-\rho T^{\alpha_\varepsilon}} v(X_{T^{\alpha_\varepsilon}} + \kappa, Y_{T^{\alpha_\varepsilon}})} \\
					&= - \delta \kappa \Expect{e^{-\rho T^{\alpha_\varepsilon}}} + \Expect{\Int{0}{T^{\alpha_\varepsilon}} e^{-\rho s} \left( dZ_s - \delta dK_s \right) + e^{-\rho T^{\alpha_\varepsilon}} v(X_{T^{\alpha_\varepsilon}} + \kappa, Y_{T^{\alpha_\varepsilon}})} \\
					&\geq - \delta \kappa \Expect{e^{-\rho T^{\alpha_\varepsilon}}} + v(x + \kappa, y) - \varepsilon .
				\end{align*}
				Finally, letting $\varepsilon \rightarrow 0$ and as $\Expect{e^{-\rho T^{\alpha_\varepsilon}}} < 1$ we conclude that:
				\begin{equation*}
					v(x, y) > v(x + \kappa, y) - \kappa \delta .
				\end{equation*}
				
				\item[ii)] If $x < 0$, capital injection of at least $|x|$ is needed to avoid ruin, incurring a cost of $\delta |x|$. Then: 
				\begin{itemize}
					\item[$\bullet$] Either $v(0, y) > 0$ and we inject at least capital $| x |$ if and only if $v(0, y) + \delta x > 0$ ,
					\item[$\bullet$] Or $v(0, y) = 0$ and we have $v(0, y) + \delta x < 0$ so we let the firm go bankrupt.
				\end{itemize}
			\end{itemize}
		\end{proof}
		
		We now provide an explicit characterisation of the value function for negative surplus values. The following result introduces a threshold that determines whether capital injection is optimal or if letting the firm go bankrupt is preferable.
		
		\begin{Proposition}[Capital injection threshold] 
			\label{PROP:capital_injection_threshold}
			Let $x < 0$ and $y \in\lbrack b,+\infty)$. We define $\kappa^\star(y) = - \frac{v(0, y)}{\delta}$. Then, the value function satisfies:
			\begin{equation*}
				v(x, y) = 
				\begin{cases}
					0 & \text{if } x < \kappa^\star(y), \\
					v(0, y) + \delta x & \text{if } \kappa^\star(y) < x < 0 .
				\end{cases}
			\end{equation*}
		\end{Proposition}
		
		\begin{proof} Let $(x,y)\in(-\infty,0)\times\lbrack b,+\infty)$. By the dynamic programming principle~\eqref{EQN:dynamic_programming_principle} and Proposition~\ref{PROP:inject_capital}, it is never optimal to inject more than $|x|$ units of capital. Let us define:
			\begin{equation*}
				\kappa^\star(y) = \inf\{ z \in \R^-, v(z, y) > 0 \} .
			\end{equation*}
			Then, for $x < 0$, we have that:
			\begin{equation*}
				v(x, y) = \max(v(0, y) + \delta x, 0) .
			\end{equation*}
			Which implies that capital is injected only if $v(0, y) + \delta x > 0$, or equivalently, if $x \geq - v(0, y)/\delta$. The result follows directly.
		\end{proof}
		
		Thus, the capital injection threshold $\kappa^\star(y)$ clearly delineates the boundary between solvency and bankruptcy, allowing us to precisely describe the insurer’s optimal behaviour in situations of financial distress.
	
	\subsection{Hamilton--Jacobi--Bellman equation}
	\label{SEC:viscosity_solution}
		
		In this section, we first state the HJB equation relateted to our control problem and then we show that the value function is the unique locally Lipschitz viscosity solution of the HJB equation. It will allow us to build a benchmark numerical method, based on the discretisation of the variational inequality satisfied by the value function $v$ in the next section.\\
		
		\noindent{}We set $\mathcal D^+:=\lbrack0,+\infty)\times(b,+\infty)$ and $\mathcal D^-:=(-\infty,0)\times(b,+\infty)$. The HJB equation associated with our control problem is given by the following variational inequality:
			\begin{equation}
				\label{EQN:optimal_dividends_HJB_equation}
				\left\{
				\begin{array}{ll}
					\min\left(\varphi \textbf{1}_{\mathcal D^-},\ (\partial_x \varphi - 1)\textbf{1}_{\mathcal D^+}, \ \delta - \partial_x \varphi, \ - \Lc \varphi \textbf{1}_{\mathcal D^+}\right) = 0\quad\textrm{ on }\mathbb R\times (b,+\infty)=\mathcal D^-\cup\mathcal D^+ \\
				\end{array}
				\right.
			\end{equation}
			where, $\Lc$ denotes the infinitesimal generator of the controlled surplus process, defined by:
			\begin{equation*}
				\Lc \varphi(x, y) := - (\rho + y)\varphi + c \partial_x \varphi + a(b - y) \partial_y \varphi + y \Int{0}{+\infty} \varphi(x - z, y + \eta)dF(z),
			\end{equation*}
			and $F$ denotes the cumulative distribution function of the claim sizes.
	
	\subsection{Viscosity solution characterisation}
	
		The HJB equation stated in Equation~\eqref{EQN:optimal_dividends_HJB_equation} reflects the optimal trade-off between three control actions: paying dividends, injecting capital to prevent ruin, or passively allowing the surplus to evolve under the stochastic environment driven by the claim process.
		Due to the complexity introduced by the two-dimensional state space, classical solutions to the HJB equation are not expected to exist. For this reason, we adopt the framework of viscosity solutions. We now define the notion of viscosity solution used throughout the paper. 
		
		\begin{Definition}[Viscosity subsolution] 
			A function $\lbar{u} : \mathbb R \times \lbrack b, + \infty) \rightarrow \R$ is said to be a viscosity subsolution of \eqref{EQN:optimal_dividends_HJB_equation} at point $(x, y)\in\mathbb R\times(b,+\infty)$ if any continuously differentiable function $\varphi : \mathbb R \times \lbrack b, + \infty) \rightarrow \R$ with $\varphi(x, y) = \lbar{u}(x, y)$ such that $\lbar{u} - \varphi$ reaches a local maximum, 0, at $(x, y)$ satisfies:
			\begin{equation*}
				\min\left(\varphi(x,y)\textbf{1}_{(-\infty,0)}(x),\ (\partial_x \varphi (x,y)- 1)\textbf{1}_{\mathbb R^+}(x), \ \delta - \partial_x \varphi(x,y), \ - \Lc \varphi (x,y)\textbf{1}_{\mathbb R^+}(x)\right) \leq 0 .
			\end{equation*}
		\end{Definition}
		
		\begin{Definition}[Viscosity supersolution]
			A function $\ubar{u} : \mathbb R\times \lbrack b, + \infty) \rightarrow \R$ is said to be a viscosity supersolution of \eqref{EQN:optimal_dividends_HJB_equation} at point $(x, y)\in\mathbb R\times(b,+\infty)$ if any continuously differentiable function $\varphi : \mathbb R \times \brack b, + \infty) \rightarrow \R$ with $\varphi(x, y) = \ubar{u}(x, y)$ such that $\ubar{u} - \varphi$ reaches a local minimum, 0,  at $(x, y)$ satisfies:
			\begin{equation*}
				\min\left(\varphi(x,y)\textbf{1}_{(-\infty,0)}(x),\ (\partial_x \varphi (x,y)- 1)\textbf{1}_{\mathbb R^+}(x), \ \delta - \partial_x \varphi(x,y), \ - \Lc \varphi (x,y)\textbf{1}_{\mathbb R^+}(x) \right) \geq 0 .
			\end{equation*}
		\end{Definition}
		
		\begin{Definition}[Viscosity solution]
			\label{DEF:viscosity_solution}
			A function $u : \mathbb R\times \lbrack b, + \infty) \rightarrow \R$ is said to be a viscosity solution of \eqref{EQN:optimal_dividends_HJB_equation} if it is both a viscosity subsolution and a viscosity supersolution.
		\end{Definition}
		
		We now justify the viscosity characterisation of the value function $v$ by proving that it satisfies the HJB equation \eqref{EQN:optimal_dividends_HJB_equation} both as a subsolution and as a supersolution in the viscosity sense. 
		
		\begin{Lemma}[Value function as viscosity supersolution]
			\label{LEMMA:supersolution}
			The value function $v$ is a viscosity supersolution of (\ref{EQN:optimal_dividends_HJB_equation}) at  every point $(x, y) \in \mathbb R\times (b, +\infty)$
		\end{Lemma}
		
		\begin{proof} Let $\varphi \in \Cc^1$ be a test function such that $v - \varphi$ has a local minimum at $(x, y) \in\mathbb R\times (b, +\infty)$ and $v(x, y) = \varphi(x, y)$. We verify the four conditions defining the viscosity supersolution are in force:
			
			\begin{itemize}
			
				\item[i)] \textbf{Bankruptcy constraint:}  As $\varphi(x,y)=v(x,y)\geq 0$ we have $\varphi(x,y) \textbf{1}_{\mathcal D^-}(x,y)\geq 0$.
			
				\item[ii)] \textbf{Dividend constraint:} Assume that $x>0$. For any $\varepsilon >0$ small enough, it follows from the possibility to immediately pay $\varepsilon$ in dividends, that we have:
				\begin{equation*}
					v(x, y) \geq v(x - \varepsilon, y) + \varepsilon \geq \varphi(x-\varepsilon,y)+\varepsilon.
				\end{equation*}
				
				As at point $(x, y)$ we have $v(x, y) = \varphi(x, y)$ and $\varphi \in \Cc^1$ we deduce:
				\begin{equation*}
					\Dx \varphi (x, y) - 1 \geq 0  .
				\end{equation*}
				For $x=0$, we know that $\varphi(-\varepsilon,y)\leq v(-\varepsilon,y))=\max( v(0,y)-\delta\varepsilon, 0).$ Moreover, $\varphi(0,y)=v(0,y)>0$ therefore, for $\varepsilon$ going to 0 we get $\partial_x\varphi(0,x)\geq \delta\geq 1$.

				\item[iii)] \textbf{Capital injection constraint:} Similarly, for any $\varepsilon > 0$, the possibility to inject capital at cost $\delta$ leads to:
				\begin{equation*}
					v(x, y) \geq v(x + \varepsilon, y) - \delta \varepsilon .
				\end{equation*}
				This implies:
				\begin{equation}
					\delta - \Dx \varphi (x, y) \geq 0  .
				\end{equation}
				
				\item[iv)] \textbf{Generator inequality:} The inequality is obvious for $x<0$, so we shall assume that $x\geq 0$. We define the stopping time $\theta_h$ as:
				\begin{equation*}
					\theta_h := \inf\{ u \geq 0: (x + cu, \tilde{\lambda}_u) \notin B(x, y)\} \wedge h .
				\end{equation*}
				
				And recall that $\tau_1$ is the time of arrival of the first claim. Let $d\tilde{\lambda}_t = a(b - \lambda_t)dt$. Then we have:
				\begin{align*}
					v(x, y) &\geq \Expect{\Int{0}{\tau_1 \wedge \theta_h} e^{-\rho s} dZ_s - \delta \Int{0}{\tau_1 \wedge \theta_h} e^{-\rho s} dK_s + e^{-\rho(\tau_1 \wedge \theta_h)} \varphi(X_{\tau_1 \wedge h}, \lambda_{\tau_1 \wedge \theta_h})} \\
					&= \Expect{\Int{0}{\tau_1 \wedge \theta_h} e^{-\rho s} dZ_s + e^{-\rho(\tau_1 \wedge \theta_h)} \varphi(X_{\tau_1 \wedge \theta_h}, \lambda_{\tau_1 \wedge \theta_h})} \\
					&= \Expect{\Int{0}{\tau_1 \wedge \theta_h} e^{-\rho u} dZ_u + e^{-\rho \theta_h} \varphi(x + c\theta_h, \tilde{\lambda}_{\theta_h}) \Ind{\theta_h < \tau_1} + e^{-\rho \tau_1} \varphi(x + c \tau_1 - Y_1, \tilde{\lambda}_{\tau_1} + \eta) \Ind{\tau_1 \leq \theta_h} } \\
					&= \mathbb{E}[\Int{0}{\tau_1 \wedge \theta_h} e^{-\rho u} dZ_u + \left( \varphi(x, y) - \rho \theta_h \varphi(x, y) + c \theta_h\Dx \varphi + a(b - y) \theta_h \Dy \varphi + o(\theta_h) \right) \Ind{\theta_h < \tau_1} ] \\
					&\quad+ \Int{0}{\theta_h} \left( \Int{0}{+\infty} \varphi(x + cs - u, \tilde{\lambda}_s + \eta) dF(u) \right) e^{-\rho s} p_{\tau_1}(s) ds .
				\end{align*}
				Rearranging the terms and dividing by $\theta_h$ we obtain:
				\begin{align*}
					0 &\geq \Expect{\Int{0}{\tau_1 \wedge \theta_h} \frac{1}{\theta_h} e^{-\rho u} dZ_u} + \left( - \rho \varphi(x, y) + c \Dx \varphi + a(b - y) \Dy \varphi + o(\theta_h) \right) \Expect{\Ind{\theta_h < \tau_1}}  \\
					&\quad+ \Int{0}{\theta_h} \left( \Int{0}{+\infty} \frac{1}{\theta_h}(\varphi(x + cs - u, \tilde{\lambda}_s + \eta) - v(x, y)) dF(u) \right) e^{-\rho s} p_{\tau_1}(s) ds \\
					&\quad+ \frac{v(x, y)}{\theta_h} \Int{0}{\theta_h} e^{-\rho s} p_{\tau_1}(s) ds .
				\end{align*}
				But we have that:
				\begin{equation*}
					\frac{v(x, y)}{\theta_h} \Int{0}{\theta_h} e^{-\rho s} p_{\tau_1}(s) ds \underset{h \rightarrow 0}{\longrightarrow} y v(x, y) .
				\end{equation*}
				Because:
				\begin{equation*}
					p_{\tau_1}(s) = (b + (y - b)e^{- a s})e^{-bs - \frac{y - b}{a}(1 - e^{as})} \underset{s \rightarrow 0}{\longrightarrow} y .
				\end{equation*}
				Finally, letting $h \rightarrow 0$ we obtain:
				\begin{equation*}
					0 \geq \Expect{Z_{0^+} - Z_0}  - (\rho + y)\varphi(x, y)  + c \Dx \varphi + a(b - y) \Dy \varphi + y \Int{0}{+\infty} \varphi(x - u, y + \eta) dF(u) .
				\end{equation*}	
				Finally, choosing a strategy $Z$ such that $\Expect{Z_{0^+} - Z_0} = 0$ we obtain:
				\begin{equation*}
					(\rho + y)\varphi(x, y) - c \Dx \varphi - a(b - y) \Dy \varphi - y \Int{0}{+\infty} \varphi(x - u, y + \eta) dF(u) \geq 0 .
				\end{equation*}
			\end{itemize}
		\end{proof}	
		
		\begin{Lemma}[Value function as viscosity subsolution]
			\label{LEMMA:subsolution}
			The value function $v$ is a viscosity subsolution of \eqref{EQN:optimal_dividends_HJB_equation} at every point $(x, y) \in \mathbb R\times (b, +\infty)^2$
		\end{Lemma}	
		
		\begin{proof} This proof is inspired by \cite{ALBRECHER2008134}. Arguing by contradiction that $v$ is not a viscosity subsolution of \eqref{EQN:optimal_dividends_HJB_equation} at point $(x, y) \in \mathbb R\times (b, +\infty)$. By definition this means that one can find $\nu > 0$ and $\varphi \in \Cc^1$ such that $\varphi(x, y) = v(x, y)$ and $\varphi(x', y') \geq v(x', y')$ for $(x', y') \in \mathbb R \times (b, +\infty)$:
			\begin{equation*}
				\min\left(\varphi(x,y)\textbf{1}_{(-\infty,0)}(x),\ (\partial_x \varphi (x,y)- 1)\textbf{1}_{(0,+\infty)}(x), \ \delta - \partial_x \varphi(x,y), \ - \Lc \varphi (x,y)\textbf{1}_{\lbrack0,+\infty)}(x) \right) > \nu .
			\end{equation*}
			\begin{itemize}
			\item[i)] If $x< 0$, this implies that
			\begin{equation*}
				\min\left(\varphi(x,y), \ \delta - \partial_x \varphi(x,y)\right) > \nu .
			\end{equation*}
			As $v(x,y)=\varphi(x,y)>0$, the optimal policy is to inject capital and there exists $\varepsilon>0$ such that $v(x,y)=v(x+\epsilon)-\delta\epsilon$. On the other hand, $\varphi$ is continuously differentiable, so for  $\varepsilon>0$ small enough $\partial_x \varphi(x',y)\leq \delta-\frac{\nu}{2}$ for $x'\in(x, x+\varepsilon)$. Integrating the last ineqality between $x$ and $x+\varepsilon$, we get
			$$\varepsilon(\delta-\frac{\nu}{2})\geq\varphi(x+\varepsilon,y)-\varphi(x,y)\geq v(x+\varepsilon,y)-v(x,y)=\varepsilon\delta.$$
			That leads to a contradiction.\\
			\item[ii)] For $x=0$, we have
			\begin{equation*}
				\min\left(\partial_x\varphi(0,y)-1,\  \delta - \partial_x \varphi(0,y), \ - \Lc \varphi (0,y) \right) > \nu .
			\end{equation*}
			Let $h>0$ and $\varepsilon > 0$. From the dynamic programming principle, there exists $(Z^\varepsilon, K^\varepsilon)\in\mathcal A(0,y)$ such that:
				\begin{align*}
					v(0,y&)\leq\Expect{\Int{0}{\tau_1\wedge h} e^{-\rho s} dZ^\varepsilon_s - \delta \Int{0}{\tau_1\wedge h} e^{-\rho s} dK^\varepsilon_s + e^{- \rho \tau_1 \wedge h} v(X^\varepsilon_{\tau_1 \wedge h}, \lambda^y_{\tau_1 \wedge h})}+\varepsilon h \\
					 &\leq \Expect{\Int{0}{\tau_1\wedge h} e^{-\rho s} d(Z^\varepsilon_s - \delta K^\varepsilon_s )+ e^{- \rho \tau_1 \wedge h} \varphi(X^\varepsilon_{\tau_1 \wedge h}, \lambda^y_{\tau_1 \wedge h})}+\varepsilon h
				\end{align*}
				It follows from Proposition \ref{PROP:inject_capital} that $K^\varepsilon_u=0$ on $\{u\leq \tau_1\}$. Hence, we have $Z^\varepsilon _u\leq cu$ on $0\leq u\leq 0\tau_1$.  It follows that there exists $\hat c\in\lbrack 0,c\rbrack$ such that, on $\{h< \tau_1\}$:
				\begin{align*}
					\Int{0}{ h} e^{-\rho s} d(Z^\varepsilon_s-\delta K^\varepsilon_s)=\Int{0}{ h} e^{-\rho s} dZ^\varepsilon_s=(c-\hat c)h+\textrm{o}(h)\textrm{ and }X^\varepsilon_h=\hat ch+\textrm{o}(h) .
				\end{align*}
				Hence, as $\varphi(0,y)=v(0,y)$, we have:
				\begin{align*}
				\varphi(0,y)  &\leq \Expect{\left((c-\hat c)h+\textrm{o}(h)+e^{-\rho h} \varphi(\hat ch+\textrm{o}(h), \tilde\lambda^y_h)\right) \Ind{\tau_1 > h} } \\
					&+ \Expect{ \left(\Int{0}{ \tau_1} e^{-\rho s} dZ^\varepsilon_s+e^{-\rho \tau_1} \varphi(\tau_1 c - Y_1, \tilde \lambda^y_{\tau_1}+\eta) \right)\Ind{\tau_1 \leq h}}+\varepsilon h,
				\end{align*}	
				where $\tilde\lambda$ is solution of the following ODE: $d\tilde\lambda_s=a(b-\tilde\lambda_s)ds.$ One can easily check that:
				\begin{equation*}
					\tilde\lambda^y_s=	(y-b)e^{-as}+b;\quad \textrm{for }s\geq 0.
				\end{equation*}
				For $s$ going to 0, we get:
				\begin{equation*}
					\tilde\lambda^y_s=y-as(y-b)+\textrm{o}(s).
				\end{equation*}
				Then we have:
				\begin{align*}	
					\varphi(0, y) &\leq e^{- \rho h} \varphi(\hat ch+\textrm{o}(h) , y - a (y-b) h + \textrm{o}(h)) \Prob{\tau_1 > h} + (c-\hat c)h + \textrm{o}(h) \\
					&\quad +\int_0^h \left( \int_0^{z^\star(y)}e^{-\rho s} \varphi(cs-z,y+\eta-as(y-b)+\textrm{o}(s))p(z) dz \right) \tilde\lambda_se^{-\int_0^s\tilde\lambda_u du} ds +\varepsilon h \\
					&= e^{-\rho h} \left( \varphi(0, y) + h \left[ \hat c \Dx{\varphi}(0, y) + a (b - y) \Dy{\varphi}(0, y) \right] + o(h) \right) e^{- \Int{0}{h} \tilde\lambda_s ds}+(c-\hat c)h+\varepsilon h+\textrm{o}(h) \\
					&\quad +h y\int_0^{z^\star(y)} \varphi(-z,y)p(z) dz+\textrm{o}(h)
				\end{align*}	
				For $h$ and then $\varepsilon$ going to 0, we obtain that:
				\begin{align*}
					\varphi(0, y){\leq} \frac{1}{\rho + y} \left[ c+\hat c( \Dx{\varphi}(0, y) -1)+ a(b - y) \Dy{\varphi}(0, y)+y\int_0^{z^\star(y)} \varphi(-z,y)p(z) dz \right] .
				\end{align*}
				As we have $ \Dx{\varphi}(0, y) -1\geq 0$, we get a contradiction between $\mathcal L \varphi(0,y)<-\nu$ and 
				\begin{equation}
					\label{EQN:boundary_condition_subsolution}
					\varphi(0, y){\leq} \frac{1}{\rho + y} \left[ c\Dx{\varphi}(0, y) + a(b - y) \Dy{\varphi}(0, y)+y\int_0^{z^\star(y)} \varphi(-z,y)p(z) dz \right] .
				\end{equation}
			\item[iii)] Assume that $x>0$ and set $B_r(x, y) \subset (0, +\infty)\times(b,+\infty)$ be a closed ball of radius $r > 0$. We define:
			\begin{equation*}
				\tau_B = \inf\{ t > 0 | X_t^\alpha \not\in B_r(x, y) \} .
			\end{equation*}		
			We denote by $\tau^\star = \tau_B \wedge T$. 
	
			\ni
			\underline{Case 1:} On $\{ \tau^\star = \tau_B \}$ two cases are possible:
			\begin{itemize}
				\item There was no jump and:
				
				\begin{equation*}
					\left\{
					\begin{array}{l}
						X_{\tau^{\star^-}}^\alpha = X_{\tau^{\star}}^\alpha = x + r \Rightarrow y - r \leq \lambda_{\tau^{\star^-}} = \lambda_{\tau^{\star}} \leq x + r, \\
						\lambda_{\tau^{\star^-}} = \lambda_{\tau^{\star}} = y - r \Rightarrow x - r \leq X_{\tau^{\star}}^\alpha \leq x + r .
					\end{array}
					\right.
				\end{equation*}		
				\item There has been a jump and:	
				\begin{equation*}
					\left\{
					\begin{array}{l}
						X_{\tau^{\star^-}}^\alpha \geq X_{\tau^{\star}}^\alpha \text{ and } X_{\tau^{\star}}^\alpha \leq x - r \text{ and } \lambda_{\tau^{\star^-}} \leq \lambda_{\tau^{\star}}, \\
						\lambda_{\tau^{\star^-}} \leq \lambda_{\tau^{\star}} \text{ and } \lambda_{\tau^{\star}} \geq y + r \text{ and } X_{\tau^{\star}}^\alpha \leq x + r .
					\end{array}
					\right.
				\end{equation*}		 
			\end{itemize}
			Taken together, these elements give us $X_{\tau^\star}^\alpha \leq x \pm r := x'$ and $\lambda_{\tau^{\star^-}} \leq \lambda_{\tau^\star}$ and as $v$ is increasing in $x$ and decreasing in $y$ we have:
			\begin{equation*}
				v(X_{\tau^\star}^\alpha, \lambda_{\tau^\star}) \leq v(x', \lambda_{\tau^\star}) \leq \varphi(x', \lambda_{\tau^\star}) \leq \varphi(X_{\tau^{\star^-}}^\alpha, \lambda_{\tau^{\star^-}}) .
			\end{equation*}
			
			\ni
			\underline{Case 2:} On $\{ \tau^\star = T \}$ we have $X_{\tau^{\star^-}}^\alpha \geq X_{\tau^\star}^\alpha$ and $\lambda_{\tau^{\star^-}} \leq \lambda_{\tau^\star}$, then we can write:
			\begin{equation*}
				v(X_{\tau^\star}^\alpha, \lambda_{\tau^\star}) \leq \varphi(X_{\tau^{\star^-}}^\alpha, \lambda_{\tau^{\star^-}}) .
			\end{equation*}
			Then, for both cases one can write:
			\begin{equation*}
				e^{- \rho \tau^\star} v(X_{\tau^\star}^\alpha, \lambda_{\tau^\star}) \leq e^{-\rho \tau^{\star^-}} \varphi(X_{\tau^{\star^-}}^\alpha, \lambda_{\tau^{\star^-}}) .
			\end{equation*}
			Recall that:
			\begin{equation*}
				\left\{
				\begin{array}{l}
					X_{\tau^{\star^-}}^\alpha = x + c\tau^{\star^-} - \Sum{k=1}{N_{\tau^{\star^-}}} Y_k - Z_{\tau^{\star^-}} + K_{\tau^{\star^-}}, \\
					d\lambda_{\tau^{\star^-}} = a(b - \lambda_{\tau^{\star^-}})d\tau^{\star^-} + \eta dN_{\tau^{\star^-}} .
				\end{array}
				\right.
			\end{equation*}
			By by applying Itô's formula on $e^{-\rho \tau^{\star^-}} \varphi(X_{\tau^{\star^-}}^\alpha, \lambda_{\tau^{\star^-}})$ we obtain:
			\begin{align*}
				e^{-\rho \tau^{\star^-}} \varphi(X_{\tau^{\star^-}}^\alpha, \lambda_{\tau^{\star^-}}) - \varphi(x, y) &= \Int{0}{\tau^{\star^-}} e^{-\rho s} \Dx \varphi(X_{s^-}^\alpha, \lambda_{s^-}^\alpha) \left[ c ds - dZ_s + dK_s \right] \\
				&\quad+\Int{0}{\tau^{\star^-}} e^{-\rho s} \Dy \varphi(X_{s^-}^\alpha, \lambda_{s^-}^\alpha) \left[ a(b - y) ds \right] - \rho  \Int{0}{\tau^{\star^-}} e^{-\rho s}  \varphi(X_{s^-}^\alpha, \lambda_{s^-}^\alpha) ds \\
				&\quad+ \Sum{\substack{0 \leq s \leq \tau^{\star^-}\\X_{s}^\alpha \neq X_{s^-}^\alpha}}{} \left( \varphi(X_{s}^\alpha, \lambda_{s}^\alpha) - \varphi(X_{s^-}^\alpha, \lambda_{s^-}^\alpha) \right)e^{-\rho s} \\
				&\quad+ \Sum{\substack{0 \leq s \leq \tau^{\star^-}\\X_{s^+}^\alpha \neq X_{s}^\alpha}}{} \left( \varphi(X_{s^+}^\alpha, \lambda_{s^+}^\alpha) - \varphi(X_{s}^\alpha, \lambda_{s}^\alpha) \right)e^{-\rho s} . 
			\end{align*}
			By construction before $\tau^\star$ we are in $B_r(x)$, so this can not be optimal to inject capital, which leads to $K = 0$ before $\tau^\star$. 	
			\begin{align*}
				e^{-\rho \tau^{\star^-}} \varphi(X_{\tau^{\star^-}}^\alpha, \lambda_{\tau^{\star^-}}) - \varphi(x, y) &= \Int{0}{\tau^{\star^-}} e^{-\rho s} \left[ \Dx \varphi(X_{s^-}^\alpha, \lambda_{s^-}^\alpha) c  + \Dy \varphi(X_{s^-}^\alpha, \lambda_{s^-}^\alpha) a(b - y)  \right]  ds  \\
				&\quad - \Int{0}{\tau^{\star^-}} e^{-\rho s} \Dx \varphi(X_{s^-}^\alpha, \lambda_{s^-}^\alpha) dZ_s  -  \rho  \Int{0}{\tau^{\star^-}} e^{-\rho s}  \varphi(X_{s^-}^\alpha, \lambda_{s^-}^\alpha) ds \\
				&\quad + \Sum{\substack{0 \leq s \leq \tau^{\star^-}\\X_{s}^\alpha \neq X_{s^-}^\alpha}}{} \left( \varphi(X_{s}^\alpha, \lambda_{s}^\alpha) - \varphi(X_{s^-}^\alpha, \lambda_{s^-}^\alpha) \right)e^{-\rho s} \\
				&\quad + \Sum{\substack{0 \leq s \leq \tau^{\star^-}\\X_{s^+}^\alpha \neq X_{s}^\alpha}}{} \left( \varphi(X_{s^+}^\alpha, \lambda_{s^+}^\alpha) - \varphi(X_{s}^\alpha, \lambda_{s}^\alpha) \right)e^{-\rho s}  .
			\end{align*}	
			
			\begin{itemize}
				\item $X_{s^+}^\alpha - X_{s}^\alpha \neq 0$ corresponds to the case where dividends has been distributed. So, we have:
				\begin{equation*}
					X_{s^+}^\alpha - X_{s}^\alpha = - (Z_{s^+} - Z_{s}),
				\end{equation*}
				hence:
				\begin{align*}
					\Sum{\substack{0 \leq s \leq \tau^{\star^-}\\X_{s^+}^\alpha \neq X_{s}^\alpha}}{} \left( \varphi(X_{s^+}^\alpha, \lambda_{s^+}) - \varphi(X_{s}^\alpha, \lambda_{s}) \right)e^{-\rho s} &\leq \Sum{\substack{0 \leq s \leq \tau^{\star^-}\\X_{s^+}^\alpha \neq X_{s}^\alpha}}{} \left( \varphi(X_{s^+}^\alpha, \lambda_{s}) - \varphi(X_{s}^\alpha, \lambda_{s}) \right)e^{-\rho s} \\
					&= - \Sum{\substack{0 \leq s \leq \tau^{\star^-}\\Z_{s^+} \neq Z_{s}^\alpha}}{} \left( \Int{0}{Z_{s^+} - Z_s} \Dx \phi (X_s^\alpha - u, \lambda_s) du  \right) .
				\end{align*}
				Using the fact that $\partial_x \varphi - 1 > \nu$ i.e. $\partial_x \varphi > 1$ we obtain:
				\begin{align*}
					- \left( \Int{0}{\tau^{\star^-}} \right. & \left. e^{-\rho s} \Dx \varphi(X_{s^-}^\alpha, \lambda_{s^-}^\alpha) dZ_s + \Sum{\substack{0 \leq s \leq \tau^{\star^-}\\Z_{s^+} \neq Z_{s}^\alpha}}{} e^{-\rho s} \left( \Int{0}{Z_{s^+} - Z_s} \Dx \phi (X_s^\alpha - u, \lambda_s) du  \right) \right) \leq \\
					&- \left( \Int{0}{\tau^{\star^-}} e^{-\rho s} dZ_s + \Sum{\substack{0 \leq s \leq \tau^{\star^-}\\Z_{s^+} \neq Z_{s}^\alpha}}{} e^{-\rho s} \left( Z_{s^+} - Z_s  \right) \right) = - \Int{0}{\tau^\star} e^{-\rho s} dZ_s .
				\end{align*}		
				
				\item $X_{s}^\alpha - X_{s^-}^\alpha$ corresponds to the case where there has been a jump in the cash process (claims has arrived). As pointed out in \cite{KULENKO2008270} and used in \cite{ALBRECHER2008134}, the following process is a martingale:
				\begin{align*}
					\Sum{\substack{0 \leq s \leq \tau^{\star^-}\\X_{s}^\alpha \neq X_{s^-}^\alpha}}{} ( \varphi(X_{s}^\alpha, \lambda_{s}^\alpha) &- \varphi(X_{s^-}^\alpha, \lambda_{s^-}^\alpha) )e^{-\rho s}  \\
					& - y \Int{0}{\tau^{\star^-}} \left( \Int{0}{+\infty} \left(\varphi(X_{s^-}^\alpha - u, \lambda_{s^-} + \eta) - \varphi(X_{s^-}^\alpha, \lambda_{s^-}) \right) dF(u) \right)  e^{-\rho s} ds .
				\end{align*}		
			\end{itemize}		
			
			We obtain:
			\begin{align*}
				e^{-\rho \tau^{\star^-}} & \varphi(X_{\tau^{\star^-}}^\alpha, \lambda_{\tau^{\star^-}}) + \Int{0}{\tau^\star} e^{-\rho s} dZ_s \leq  \varphi(x, y)  \\
				&+\Int{0}{\tau^{\star^-}} e^{-\rho s} \left[ c \Dx \varphi(X_{s^-}^\alpha, \lambda_{s^-}^\alpha) + a(b - y) \Dy \varphi(X_{s^-}^\alpha, \lambda_{s^-}^\alpha) - \rho \varphi(X_{s^-}^\alpha, \lambda_{s^-}^\alpha) \right. \\
				&\left. + y \Int{0}{+\infty} \left(\varphi(X_{s^-}^\alpha - u, \lambda_{s^-} + \eta) - \varphi(X_{s^-}^\alpha, \lambda_{s^-}) \right) dF(u) \right] ds .
			\end{align*}				
			Finally:
			\begin{align*}
				\varphi(x, y) = v(x, y) &\leq \Expect{\Int{0}{\tau^{\star^-}} e^{-\rho s} \left( dZ_s - \delta dK_s \right) + e^{\tau^{\star^-}} v(X_{\tau^{\star^-}}^\alpha, \lambda_{\tau^{\star^-}})} + \varepsilon \\
				&\leq \Expect{\Int{0}{\tau^{\star^-}} e^{-\rho s} dZ_s + e^{\tau^{\star^-}} v(X_{\tau^{\star^-}}^\alpha, \lambda_{\tau^{\star^-}})} + \varepsilon \\
				\leq \varphi(x, y) + &\mathbb{E}\left[ \Int{0}{\tau^{\star^-}} e^{-\rho s} \left[ c \Dx \varphi(X_{s^-}^\alpha, \lambda_{s^-}^\alpha) + a(b - y) \Dy \varphi(X_{s^-}^\alpha, \lambda_{s^-}^\alpha) \right. \right. \\
				&- \left. \left. \rho \varphi(X_{s^-}^\alpha, \lambda_{s^-}^\alpha) + y \Int{0}{+\infty} \left(\varphi(X_{s^-}^\alpha - u, \lambda_{s^-} + \eta) - \varphi(X_{s^-}^\alpha, \lambda_{s^-}) \right) dF(u) \right] ds \right] + \varepsilon .
			\end{align*}	
			This implies the following contradiction:
			\begin{equation*}
				0 \leq - \Expect{\Int{0}{\tau^{\star^-}} e^{- \rho s} \Lc \varphi(X_{s^-}, \lambda_{s^-}) ds} \leq - \nu \Expect{\Int{0}{\tau^{\star^-}} e^{-\rho s} ds } < 0 .
			\end{equation*}
			\end{itemize}
		\end{proof}
		A direct consequence of Lemmas \ref{LEMMA:subsolution} and \ref{LEMMA:supersolution} is the following result.
		\begin{Theorem}[Value function as viscosity solution]
			\label{THM:viscosity_solution}
			The value function $v$ is a viscosity solution of the Hamilton--Jacobi--Bellman equation \eqref{EQN:optimal_dividends_HJB_equation} on the domain $\mathbb R \times (b, +\infty)$.
		\end{Theorem}
		
		\begin{Remark} Notice that the variational inequality at points $(0,y)$ could be considered as a boundary condition because it could be written as
		\begin{equation}
				\label{EQN:boundary_condition}
					(\rho+y) \varphi(0, y) = c \frac{\partial \varphi}{\partial x}(0, y) + a(b - y) \frac{\partial \varphi}{\partial y} (0, y) + y\Int{0}{\varphi(0, y+\eta)/\delta} \left( \varphi(0, y+\eta) -  \delta z \right) dF(z) .
			\end{equation}
		
				\end{Remark}
		
		Following the proof of Proposition 4.2 in\cite{albrecher2023optimal}, we can now give a characterisation of the value function as the smallest viscosity supersolution of equation \eqref{EQN:optimal_dividends_HJB_equation}.
		
		\begin{Theorem}
			\label{THM:comparison_uniqueness}
			$v$ is the smallest viscosity supersolution of \eqref{EQN:optimal_dividends_HJB_equation} that is non-increasing in $y$, locally Lipschitz continuous and satisfies the growth condition established in Proposition \ref{PROP:value_function_bounds}.
		\end{Theorem}
	
\section{Finite-difference estimate}
	\label{SEC:pde_resolution}
	In this section, we present the classical finite-difference scheme used as a numerical benchmark for the solution of the HJB variational inequality~\eqref{EQN:optimal_dividends_HJB_equation}. 
	The method relies on a monotone discretisation of the state dynamics combined with Howard’s policy iteration algorithm to obtain the stationary solution. 
	This framework provides a consistent and interpretable reference against which the reinforcement learning approach introduced later can be compared.
	\subsection{Discrete HJB variational inequality}
		\subsubsection*{Computational grid and domain truncation}
			To approximate the value function numerically, we truncate the state space and construct a uniform grid over the resulting bounded domain. 
			Let $X_{\min} < 0 < X_{\max}$ and $Y_{\max} \in (b,+\infty)$, and define $\mathcal{D} := [X_{\min} ,X_{\max}] \times [b, Y_{\max}]$ as the computational domain for the surplus and intensity variables. 
			The domain $\mathcal{D}$ is discretised using $N_x \in \N$ and $N_y \in \N$ spatial subdivisions along the $x$ and $y$ directions, respectively, leading to the mesh sizes
			\begin{equation*}
				\Delta x := \frac{X_{\max}-X_{\min}}{N_x}, 
				\qquad
				\Delta y := \frac{Y_{\max}-b}{N_y}.
			\end{equation*}
			The corresponding grid points are defined by
			\begin{equation*}
				x_i := X_{\min} + i \Delta x, 
				\qquad
				y_j := b + j \Delta y.
			\end{equation*}
			The full grid is therefore
			\begin{equation*}
				\mathcal{G} := \{(x_i,y_j): 0 \le i \le N_x,  0 \le j \le N_y\},
			\end{equation*}
			with its interior nodes denoted by
			\begin{equation*}
				\mathcal{G}^\circ := \{(x_i,y_j)\in\mathcal{G}:  1 \le i \le N_x-1,  1 \le j \le N_y-1\}.
			\end{equation*}
			At each grid point $(x_i, y_j) \in \mathcal{G}$, the numerical approximation of the value function is denoted by $V_{i,j} \approx v(x_i, y_j)$ and will be used consistently throughout the discrete formulation. 
			The index corresponding to the origin $x=0$ is denoted by $i_0$, so that $x_{i_0}=0$. 
			
		\subsubsection*{Discretisation of differential operators}
			We approximate the differential operators in the HJB variational inequality by means of a monotone finite-difference discretisation. 
			One-sided differences are employed in each direction to preserve the directionality of the underlying drift terms. 
			For $x$ and $y$ coordinates, the discrete first-order operators are defined as
			\begin{equation*}
				D^-_x V_{i,j} := \frac{V_{i,j}-V_{i-1,j}}{\Delta x}, 
				\qquad
				D^+_x V_{i,j} := \frac{V_{i+1,j}-V_{i,j}}{\Delta x},
			\end{equation*}
			and similarly
			\begin{equation*}
				D^-_y V_{i,j} := \frac{V_{i,j}-V_{i,j-1}}{\Delta y}, 
				\qquad
				D^+_y V_{i,j} := \frac{V_{i,j+1}-V_{i,j}}{\Delta y}.
			\end{equation*}
			For a generic convective term $s \partial_\xi V$ with $\xi \in \{x,y\}$, an upwind discretisation is adopted:
			\begin{equation*}
				s \partial_\xi V  \approx 
				\begin{cases}
					s D^-_\xi V, & \text{if } s \ge 0, \\[3pt]
					s D^+_\xi V, & \text{if } s < 0.
				\end{cases}
			\end{equation*}
			In particular, since $c>0$, we approximate $-c \partial_x V$ by $c D^+_x V$. 
			For the intensity dynamics, the drift satisfies $a(b-y)\le 0$ on $[b, Y_{\max}]$, yielding
			\begin{equation*}
				-a(b-y_j) \partial_y V_{i,j} \approx -a(b-y_j) D^-_y V_{i,j}.
			\end{equation*}
		
		\subsubsection*{Approximation of the jump integral}
			The discrete infinitesimal generator $\mathcal{L}_h$ acting on the grid interior $\mathcal{G}^\circ$ is then defined by
			\begin{equation}
				\label{eq:disc_generator}
				-\mathcal{L}_h V_{i,j} := (\rho + y_j) V_{i,j}
				- c D^+_x V_{i,j}
				- a(b-y_j) D^-_y V_{i,j}
				- y_j \mathcal{Q}_h[V]_{i,j},
			\end{equation}
			where $\mathcal{Q}_h$ denotes the discrete approximation of the jump operator
			\begin{equation*}
				\mathcal{Q}[V](x,y) := \int_0^{+\infty} V(x-z, y+\eta) dF(z).
			\end{equation*}
			To approximate this integral, we truncate the support of $f$ at $Z_{\max} = (M+\tfrac{1}{2})\Delta x$ and apply a midpoint quadrature rule on $[0, Z_{\max}]$:
			\begin{equation*}
				\int_{0}^{Z_{\max}} V(x-z, y+\eta) f(z) dz
			 \approx 
				\sum_{m=0}^{M} 
				V\!\big(x-(m+\tfrac{1}{2})\Delta x, y+\eta\big) 
				f\!\big((m+\tfrac{1}{2})\Delta x\big) \Delta x,
			\end{equation*}
			When $x_i - (m+\tfrac{1}{2})\Delta x < 0$, the capital injection condition given in Proposition~\ref{PROP:capital_injection_threshold} is enforced to evaluate $V$, while off-grid values are obtained by bilinear interpolation.
		
		\subsubsection*{Discrete HJB variational inequality}
			Combining the spatial and integral approximations introduced above, the discrete counterpart of the HJB variational inequality takes the form
			\begin{equation}
				\label{eq:disc_VI}
				\min\Big(
				\underbrace{D^-_x V_{i,j}-1}_{\text{dividends}},
				\underbrace{\delta - D^+_x V_{i,j}}_{\text{capital injection}},
				\underbrace{-\mathcal{L}_h V_{i,j}}_{\text{continuation}}
				\Big) = 0,
				\qquad (i,j) \in \mathcal{G}^\circ,
			\end{equation}
			The resulting non-linear system is monotone and consistent with the viscosity framework, providing a robust basis for numerical resolution.
			
	\subsection{Numerical implementation}
	
		\subsubsection*{Local update rules}
			The discrete variational inequality~\eqref{eq:disc_VI} defines, at each grid node, the local optimality condition between the three possible regimes: dividend payment, capital injection, and continuation. 
			In practice, this translates into a set of region-specific update formulas that can be used to iteratively compute the value function over the grid. 
			The expressions below follow directly from the monotone discretisation introduced in the previous subsection.
			
			In the dividend region, the optimal action corresponds to an immediate payout, leading to the first-order condition $D^-_x V_{i,j} = 1$, which yields
			\begin{equation*}
				V_{i,j} = V_{i-1,j} + \Delta x.
			\end{equation*}
			In the continuation region, the process evolves according to the controlled surplus dynamics without intervention, and the value function satisfies the discrete HJB equation obtained from~\eqref{eq:disc_generator}:
			\begin{equation*}
				V_{i,j}
				= \left(\rho + y_j + \frac{c}{\Delta x} + \frac{a(y_j - b)}{\Delta y}\right)^{-1} \left( \frac{c}{\Delta x} V_{i+1,j}
				+ \frac{a(y_j - b)}{\Delta y} V_{i,j-1}
				+ y_j \mathcal{Q}_h[V]_{i,j} \right).
			\end{equation*}
			The capital injection region requires a specific treatment, as its behaviour is entirely characterised by Proposition~\ref{PROP:capital_injection_threshold}. 
			According to this result, the value function is known explicitly for $x < 0$, where injections occur whenever the surplus lies below the optimal boundary. 
			Hence, the relation
			\begin{equation*}
				V_{i,j} = \max\big(0,  V_{i_0,j} + \delta x_i \big),
			\end{equation*}
			is imposed directly as a boundary condition for all grid points with $x_i < 0$, ensuring consistency with the theoretical characterisation of the optimal policy.

		\subsubsection*{Boundary conditions}
			In the negative surplus region $x < 0$, the value function is entirely determined by the theoretical characterisation established in Proposition~\ref{PROP:capital_injection_threshold}, which directly governs the capital injection mechanism. 
			Hence, no numerical update is required in this area, and the boundary relation at $x=0^+$ serves as the effective entry condition for the computational domain.
			According to the HJB equation~\ref{EQN:optimal_dividends_HJB_equation}, for $j \in \{0, \dots, N_y\}$ the value function at $x = 0$ satisfies
			\begin{equation*}
				(\rho + y_j) V_{i_0,j}
				= c D^+_x V_{i_0,j} + a(b - y_j) D^-_y V_{i_0,j} 
				+ y_j I_h\left(V_{i_0, j}^{(\eta)}/\delta\right),  
			\end{equation*}
			where we recall that $i_0$ denotes the index corresponding to $x_{i_0} = 0$.
			Here, the term $V_{i_0,j}^{(\eta)}$ represents the numerical approximation of $v(0,\, y_j+\eta)$ obtained by linear interpolation, while $I_h(\cdot)$ denotes the approximation of the integral term arising from the infinitesimal generator at $x = 0$, using the injection characterisation given in Proposition~\ref{PROP:capital_injection_threshold}.
			In practice, $I_h$ can be evaluated using a midpoint or trapezoidal rule depending on the discretisation of $F$, although for many standard claim size distributions, this integral admits a closed-form expression, allowing for an exact and computationally efficient evaluation.

			At the upper boundary of the intensity domain, the asymptotic behaviour derived in Corollary~\ref{COR:yasymptotic} implies $\lim_{y \to \infty} v(x,y) = x$ for all $x \in \R^+$, which translates numerically into
			\begin{equation*}
				V_{i, N_y} = x_i,
			\end{equation*}
			for all $i \ge i_0$.
			Together, these two boundary conditions fully close the discrete problem and ensure the well-posedness of the numerical resolution of the value function.
			
		\subsubsection*{Howard policy iteration}
			The non-linear discrete system~\eqref{eq:disc_VI} is solved using Howard’s policy iteration algorithm. 
			The method alternates between a policy evaluation step, where the value function is computed for a fixed control configuration, and a policy improvement step, where the control is updated pointwise according to the minimisation operator in~\eqref{eq:disc_VI}. 
			Starting from an initial value function $V^{(0)}$ and an initial policy $\pi^{(0)}$, the iteration proceeds as follows:
			\begin{enumerate}
				\item[\textit{(i)}] \textbf{Policy evaluation:} 
				For a fixed policy $\pi^{(k)}$, the corresponding value function $V^{(k+1)}$ is obtained by solving the discrete HJB system~\eqref{eq:disc_generator} induced by this policy. 
				This consists in applying, at each grid point, the update rule associated with the prescribed regime. 
				The system is solved by fixed-point iteration until convergence, under the boundary conditions described above.
				
				\item[\textit{(ii)}] \textbf{Policy improvement:}
				The policy is then updated pointwise by selecting the locally optimal regime,
				\begin{equation*}
					\pi^{(k+1)}(i,j) = \arg\min_{\pi \in \{\text{dividend, injection, continuation}\}} 
				\mathcal{H}_h^\pi[V^{(k+1)}]_{i,j},
				\end{equation*}
				where $\mathcal{H}_h^\pi$ denotes the local discrete HJB operator associated with regime $\pi$.
			\end{enumerate}
			The algorithm iterates between these two steps until the policy stabilizes, that is, when $\pi^{(k+1)} = \pi^{(k)}$ over the grid, indicating convergence to the stationary optimal control.
			The convergence of the numerical scheme follows from standard arguments for monotone finite-difference approximations and policy iteration methods. 
			Under the usual monotonicity, consistency, and stability assumptions on the discrete operator~$\mathcal{L}_h$, the fixed-point evaluation step preserves the viscosity solution framework of the continuous HJB variational inequality~\cite{Barles1991}. 
			Moreover, the outer Howard iteration, alternating between policy evaluation and improvement, converges to the unique stationary solution of the discrete control problem under these same structural conditions~\cite{Puterman1994, KushnerDupuis2001}.
			Overall, the scheme is guaranteed to converge to the discrete viscosity solution, which consistently approximates the continuous value function as the mesh is refined.

	\subsection{Numerical results and sensitivity analysis}
		\subsubsection{Reference configuration and qualitative analysis}
			\paragraph{Model and grid setup} ~\\
				We begin with a balanced baseline configuration of parameters, chosen to represent a typical regime where claim arrivals, excitation effects, and premium inflows are of comparable magnitude. In particular, claim sizes are assumed to follow an exponential distribution with parameter $\beta$.
				This setup serves as a reference for the numerical results presented below and will later be used to assess the sensitivity of the optimal policy to individual model parameters.
				The corresponding values are reported in Table~\ref{TAB:model_parameters_PDE}, while the discretisation settings are summarised in Table~\ref{tab:grid_params}.
				\begin{table}[H]
					\centering
					\begin{tabular}{ccccccc}
						\toprule
						$a$ & $b$ & $\eta$ & $\rho$ & $c$ & $\delta$ & $\beta$ \\
						\midrule
						$2.0$ & $2.0$ & $0.4$ & $0.1$ & $1.0$ & $1.8$ & $3.0$ \\
						\bottomrule
					\end{tabular}
					\caption{Baseline configuration of model parameters.}
					\label{TAB:model_parameters_PDE}
				\end{table}
				Instead of fixing $N_x$ and $N_y$ directly, we define the grid resolution through the auxiliary parameters $M$ and $n_{\eta}$, which determine the number of discretisation steps relative to $Z_{\max}$ and $\eta$. 
				This construction ensures that $\Delta y$ is an exact multiple of $\eta$ and $\Delta x$ an exact multiple of $Z_{\max}$, thereby avoiding interpolation errors when evaluating the jump and excitation terms. 
				The origin $x=0$ is explicitly enforced to belong to the grid, with minor adjustments of the bounds if necessary.
				\begin{table}[H]
					\centering
					\begin{tabular}{cccccc}
						\toprule
						$X_{\min}$ & $X_{\max}$ & $Y_{\max}$ & $n_{\eta}$ & $M$ & $Z_{\max}$ \\
						\midrule
						$-5.0$ & $4.0$ & $25.0$ & $8$ & $80$ & $5.0$ \\
						\bottomrule
					\end{tabular}
					\caption{Grid parameters used for the numerical discretisation.}
					\label{tab:grid_params}
				\end{table}
				
\begin{figure}[H]
					\centering
					\includegraphics[width=0.55\linewidth]{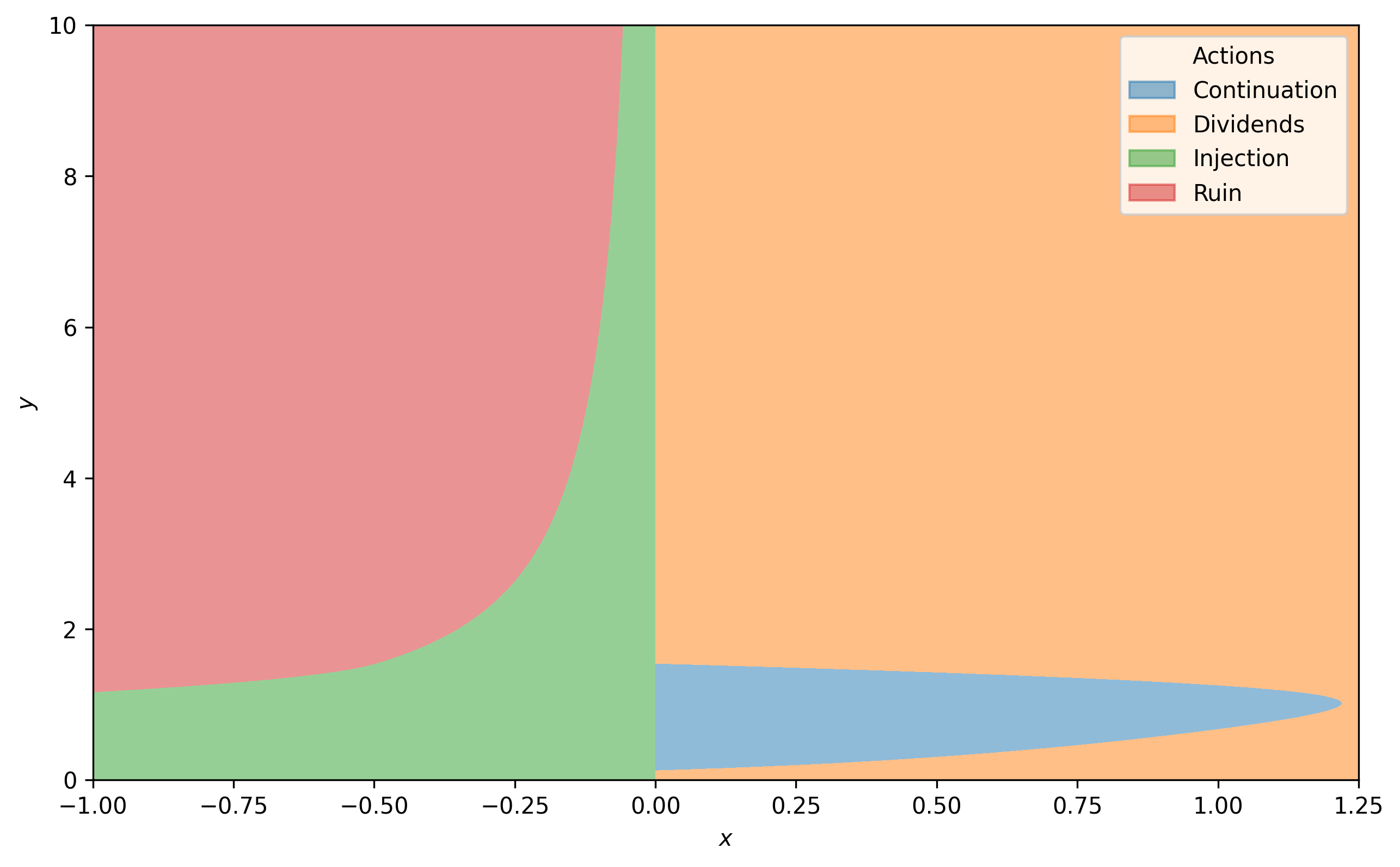}
					\caption{Optimal policy in the compound Poisson model for varying intensity levels.}
					\label{FIG:poisson_policy}
				\end{figure}

			\paragraph*{Poisson benchmark} ~\\
				When $\eta = 0$ and $\lambda_0 = b$, the intensity remains constant at~$b$ and the model reduces to the classical compound Poisson framework studied in~\cite{KULENKO2008270}. In this case, the optimal dividend barrier and injection threshold admit closed-form expressions when claim sizes are exponentially distributed. 
				Figure~\ref{FIG:poisson_policy} displays the optimal policy for a range of constant intensity values. In the Poisson model, the state space reduces to the single variable~$x$ for each fixed~$\lambda$: the dynamics evolve along the horizontal axis only, and the control regions depend solely on the surplus level. The policy is shown in the $(x,y)$ plane to facilitate comparison with the Hawkes case studied below, where the intensity $\lambda_t \geq b$ becomes a dynamic second state variable.
								
\begin{figure}[H]
					\centering
					\begin{subfigure}{0.48\linewidth}
						\centering
						\includegraphics[width=\linewidth]{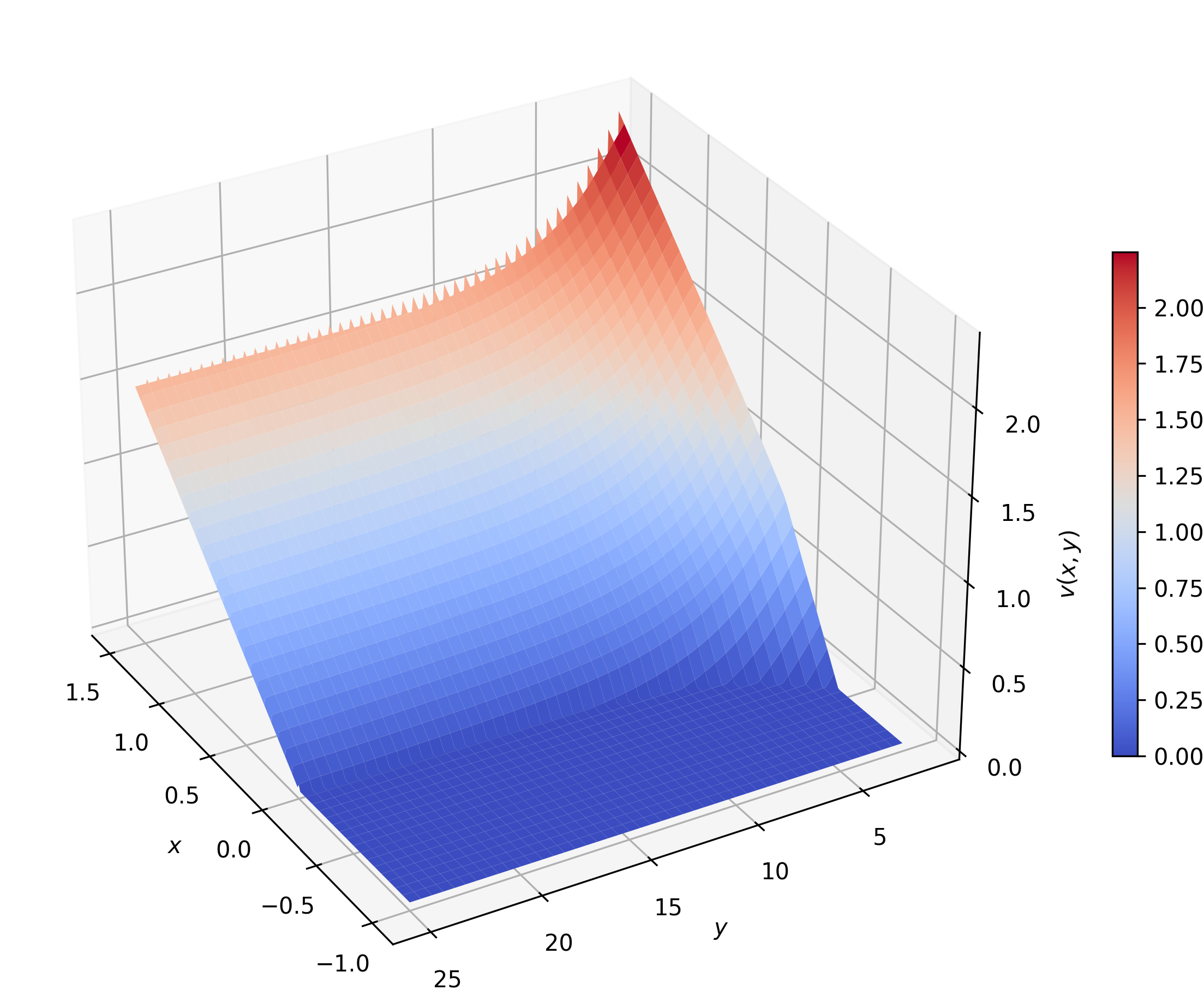}
						\caption{Estimated value function under the baseline parameter configuration.}
						\label{FIG:baseline_value_function_PDE}
					\end{subfigure}
					\hfill
					\begin{subfigure}{0.48\linewidth}
						\centering
						\includegraphics[width=\linewidth]{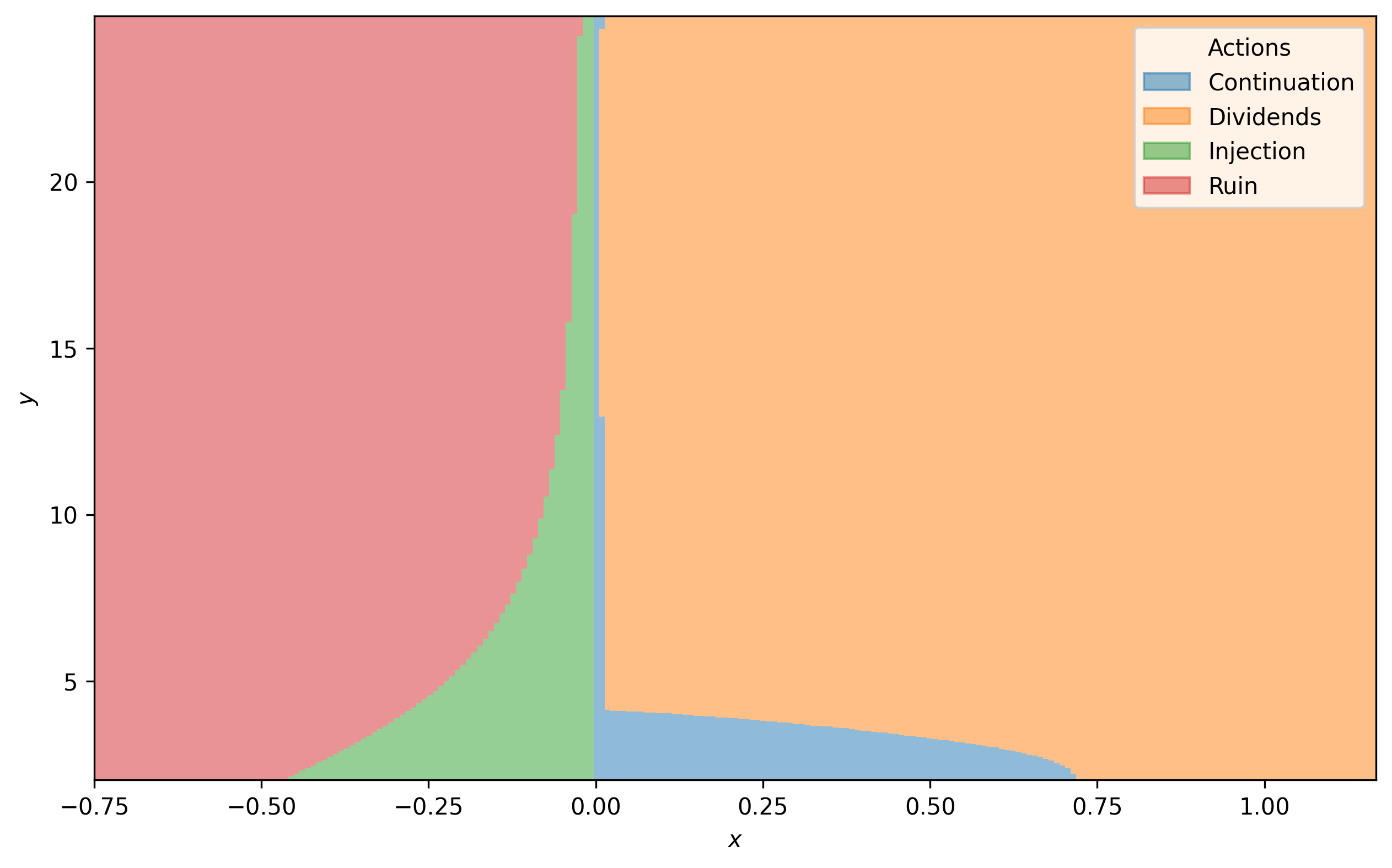}
						\vspace{0.1cm}
						\caption{Optimal control policy under the baseline parameter configuration.}
						\label{FIG:baseline_optimal_policy}
					\end{subfigure}
					\caption{Estimated value function and corresponding optimal control policy under the baseline parameter configuration.}
					\label{FIG:baseline_value_and_policy}
				\end{figure}

			\paragraph*{Value function and associated optimal policy} ~\\
				Figure~\ref{FIG:baseline_value_function_PDE} displays the estimated value function obtained from the finite-difference scheme. 
				The surface exhibits the expected qualitative behaviour: the value increases with the surplus $x$ and decreases with the claim intensity $y$, reaching its highest levels for large surpluses and low intensities. 
				These numerical patterns are fully consistent with the theoretical monotonicity properties established in Propositions~\ref{PROP:increase_of_v} and~\ref{PROP:value_function_monotony}, confirming the accuracy and stability of the discretisation procedure.
				Figure~\ref{FIG:baseline_optimal_policy} shows the optimal control policy under the baseline configuration.
				The solution exhibits a threshold structure, with two distinct regions for $x < 0$ (a ruin region and a capital-injection region) and two regions for $x \ge 0$ (a continuation region and a dividend region), yielding a clear and interpretable partition of the state space.			
				In the positive surplus region, the policy exhibits the expected two-zone structure: a continuation region and a dividend region.
				For sufficiently large surplus levels, it is always optimal to distribute dividends. This behaviour reflects the fact that the insurer holds enough reserves to absorb potential losses, making the immediate distribution of excess capital preferable.
				By contrast, the continuation region corresponds to states in which the activity remains exposed to significant risk. In this zone, it is optimal to retain earnings until the reserve reaches a safer level, at which point dividend payments resume.
				Comparing Figure~\ref{FIG:baseline_optimal_policy} with the Poisson benchmark of Figure~\ref{FIG:poisson_policy} highlights the impact of self-excitation on the optimal policy.
				In the Poisson model, each horizontal slice $y = \lambda$ corresponds to an independent problem with fixed intensity. For intensity levels above the baseline $b = 2$, the dividend barrier is close to zero: it is nearly always optimal to distribute immediately.
				In the Hawkes model, by contrast, the intensity is dynamic: a high current value of~$\lambda_t$ signals an elevated but transient risk, since the intensity reverts toward~$b$ between claims. The optimal policy exploits this by prescribing a substantial continuation region at moderate intensity levels, where the insurer retains reserves in anticipation of the intensity subsiding. This mechanism has no counterpart in the memoryless Poisson setting, where each intensity level is permanent.
				A critical feature emerging from the numerical solution is the existence of an intensity threshold $y$ above which the optimal action is to liquidate the surplus down to $x = 0^+$. 
				In this high-intensity regime, there exists an increased and persistent risk of claim occurrences, leading to a high likelihood of large loss clusters and little chance that the intensity will decline rapidly enough to restore profitability.
				Operating under such conditions is no longer profitable, and the optimal strategy is to distribute all available capital before the firm is driven to ruin.

				In the negative surplus region, the numerical policy reproduces the expected qualitative behaviour, featuring a clear capital-injection region and a ruin region. The boundary separating these two zones coincides with the one derived in Proposition~\ref{PROP:capital_injection_threshold}. Since $\kappa^\star(y) = -v(0,y)/\delta$ and Corollary~\ref{COR:yasymptotic} establishes that $v(0,y) \to 0$ as $y \to \infty$, the convergence of the injection boundary toward $0$ for large intensities is fully consistent with the theoretical predictions.
				Beyond this boundary, capital injection is no longer optimal. When incoming claims push the cash reserves past this threshold, the activity becomes too costly to refinance. Injecting capital up to $x=0$ would not generate future earnings sufficient to offset the cost of the refinancing itself. 
				In such circumstances, further investment is economically dominated, and the optimal decision is to let ruin occur.
				
		\subsubsection{Sensitivity of the optimal policy}
			\label{SEC:sensitivity_policy}
			We now examine how the optimal control policy reacts to changes in the model parameters.
			Each parameter is varied independently around its baseline value while keeping the others fixed. 
			The resulting policy maps illustrate how the intervention thresholds adapt to the underlying economic and risk conditions. 
			Overall, the numerical outcomes remain consistent with theoretical expectations and economic intuition.
			\paragraph*{Impact of Hawkes dynamics parameters} ~\\
				The parameters $(a,b,\eta)$ govern the temporal behaviour of the claim intensity process. 
				An increase in the mean-reversion rate $a$ accelerates the return of $\lambda_t$ to its baseline level $b$, reducing the persistence of high-intensity episodes. 
				This results in wider continuation and injection regions, as the system spends less time in high-risk states.
				In contrast, a higher excitation parameter $\eta$ amplifies clustering effects, making the environment significantly more risky. 
				When the self-excitation of future claim arrivals makes the business unprofitable, the optimal strategy shifts toward full liquidation: distributing all available surplus rather than continuing operations.
				\begin{figure}[H]
					\centering
					\begin{subfigure}[t]{0.48\linewidth}
							\centering
							\includegraphics[width=\linewidth]{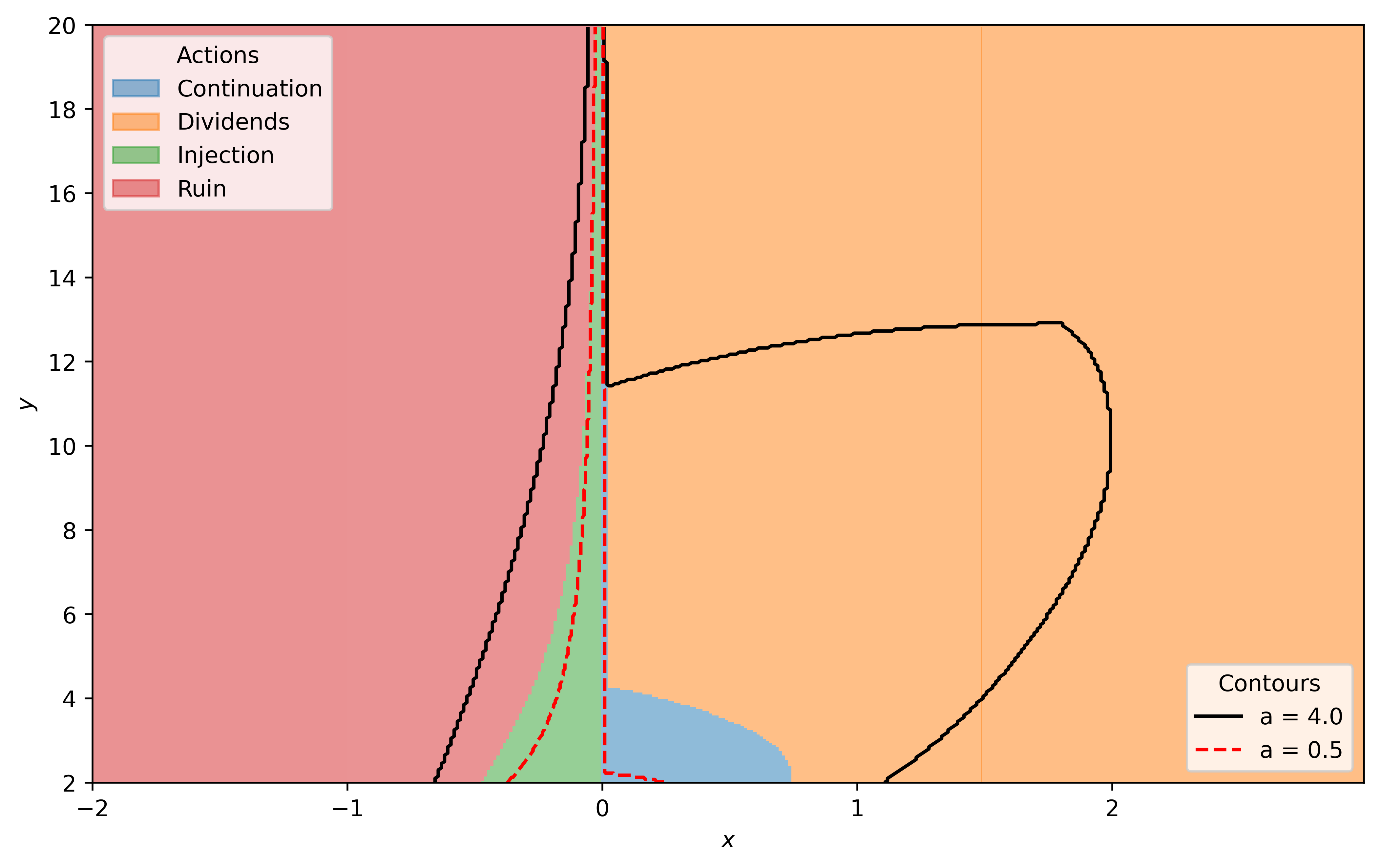}
							\caption{Sensitivity to the mean-reversion rate $a$.}
							\label{FIG:policy_sensi_to_a}
						\end{subfigure}
					\hfill
					\begin{subfigure}[t]{0.48\linewidth}
							\centering
							\includegraphics[width=\linewidth]{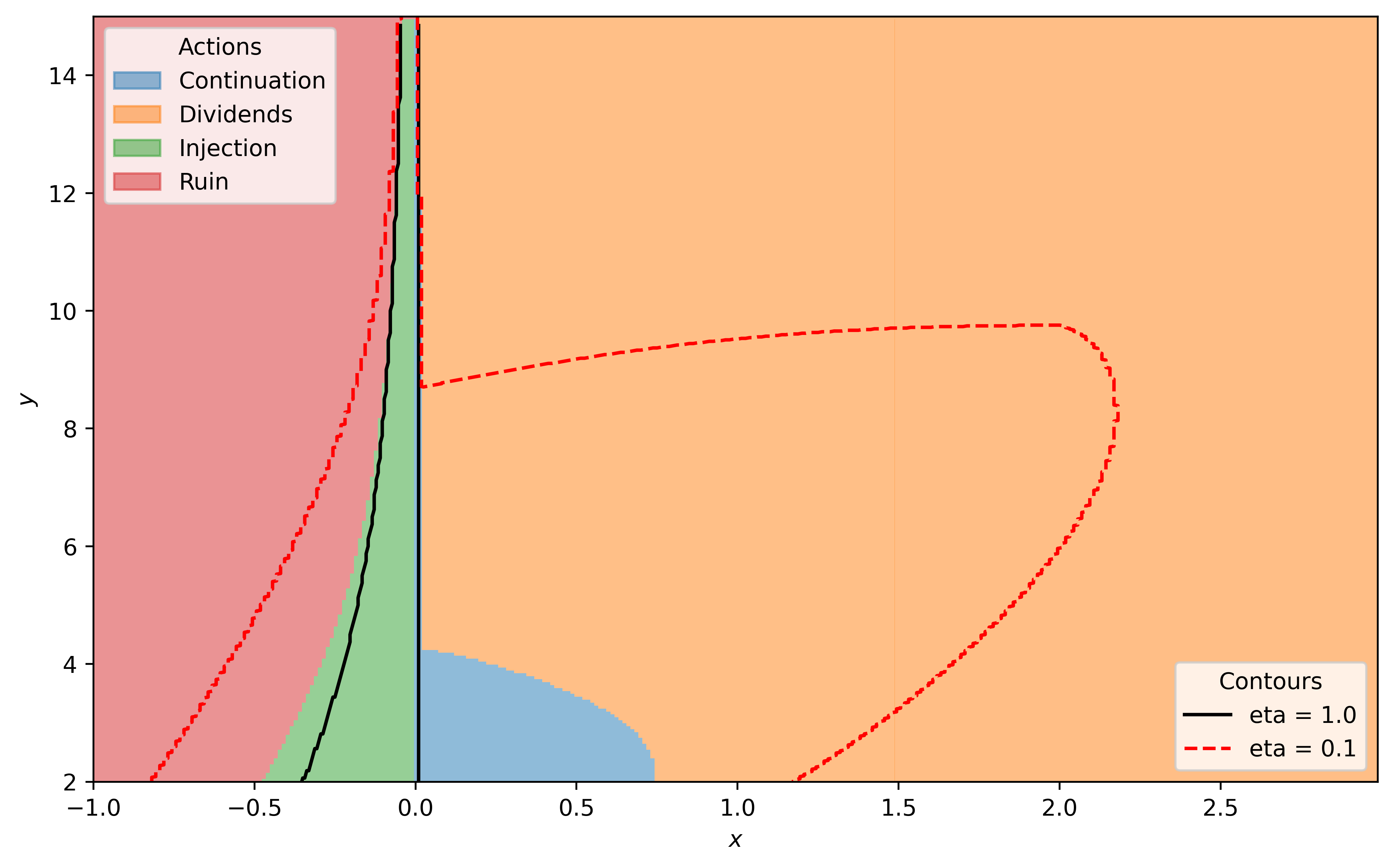}
							\caption{Sensitivity to the excitation parameter $\eta$.}
							\label{FIG:policy_sensi_to_elta}
						\end{subfigure}
					\caption{Sensitivity of optimal policy to Hawkes dynamics parameters.}
					\label{FIG:optimal_policy_sensi_hawkes_parameters}
				\end{figure}
				A closer inspection of the continuation region also reveals a distinctive shape that depends sensitively on the model parameters. 
				For a fixed but sufficiently high intensity level, the optimal policy in the positive surplus region may switch from dividend distribution to continuation and then back to dividend distribution as $x$ increases. 
				This non-monotone pattern appears for specific parameter configurations, such as $a = 4$ and $\eta = 0.1$, but also emerges under other parameter variations in the subsequent sensitivity analyses.
				
				The initial dividend region observed at low surplus levels reflects situations where the intensity has risen too sharply for profitability to be restored. Such states necessarily arise from a sequence of adverse claims originating in the continuation region, which simultaneously depletes the surplus and drives the intensity upward. Under these conditions, continued operation is no longer viable, and the optimal action is to liquidate the available surplus immediately.
				For the same intensity level, a slightly higher surplus would allow the insurer to absorb potential short-term losses while waiting for the intensity to revert, making continuation preferable. 
				As the surplus becomes large, the policy reverts to its usual behaviour: the company holds enough reserves to withstand adverse shocks, and distributing dividends again becomes optimal.
				This layered structure of the continuation region thus captures a subtle interplay between short-term risk exposure and long-term mean reversion in the intensity dynamics, and aligns with the economic interpretation of the Hawkes-driven claim environment.
	
			\paragraph*{Impact of the premium–claim balance} ~\\
				The premium rate $c$ determines the rate of surplus accumulation, directly affecting the insurer’s capacity to sustain operations. 
				Higher values of $c$ expand the continuation region and postpone both injections and dividend payments.  
				In contrast, the claim size parameter $\beta$ affects the expected cost of claims, with larger $\beta$ (smaller expected losses) 
				leading to higher profitability and a broader dividend region.  
				\begin{figure}[H]
					\centering
					\begin{subfigure}[t]{0.48\linewidth}
						\centering
						\includegraphics[width=\linewidth]{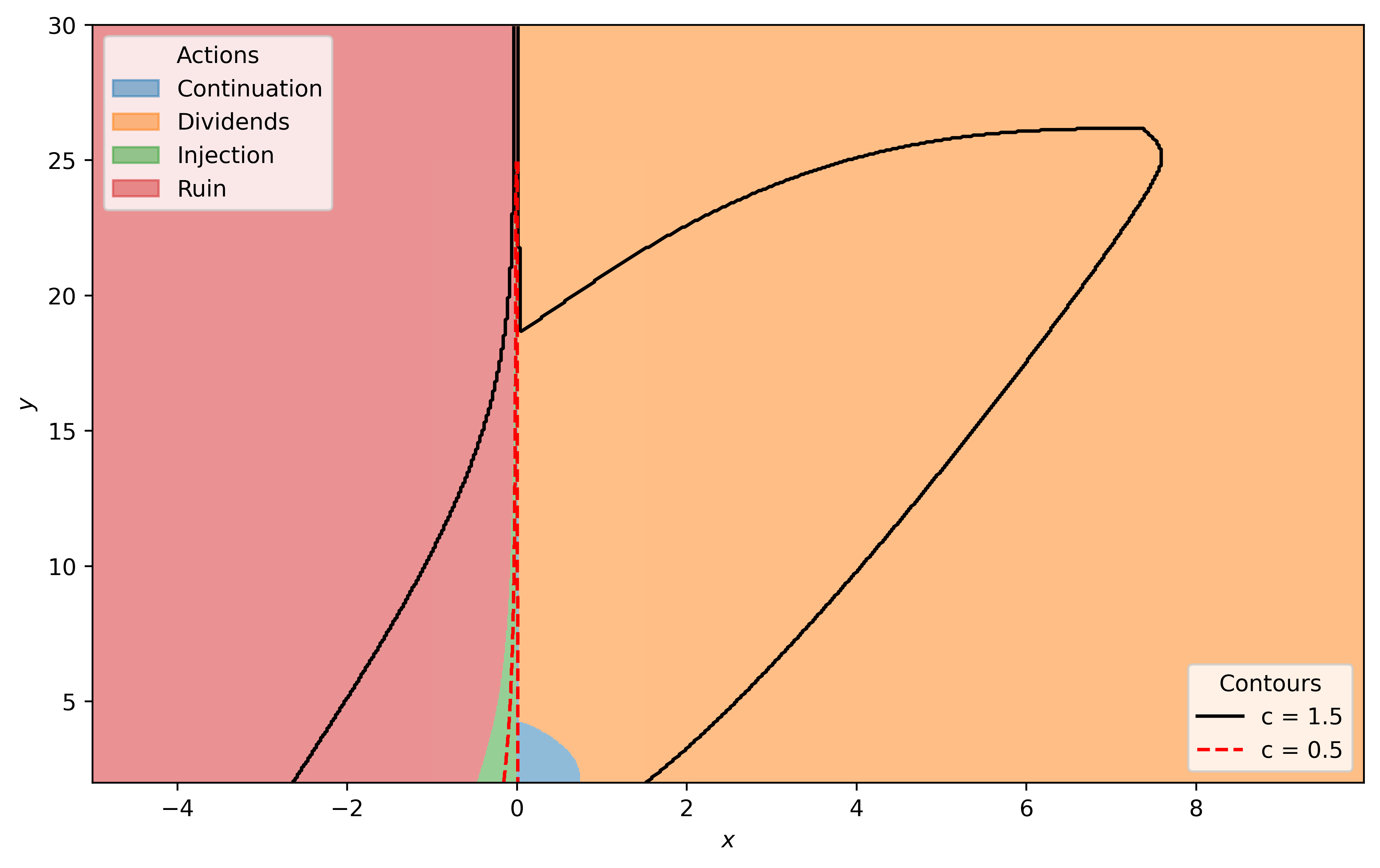}
						\caption{Sensitivity to the premium rate $c$.}
						\label{FIG:policy_sensi_to_c}
					\end{subfigure}
					\hfill
					\begin{subfigure}[t]{0.48\linewidth}
						\centering
						\includegraphics[width=\linewidth]{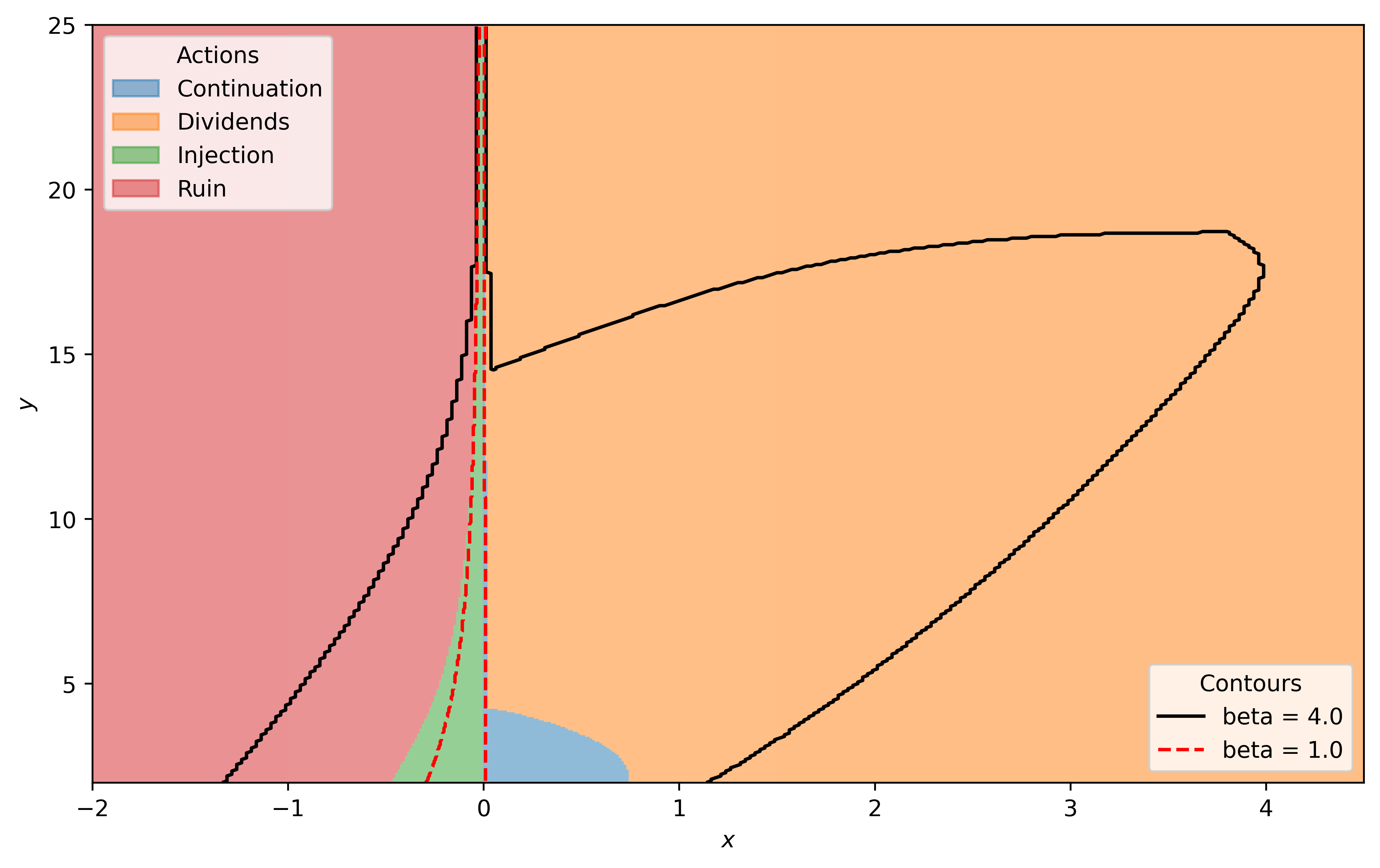}
						\caption{Sensitivity to the claim size parameter $\beta$.}
						\label{FIG:policy_sensi_to_beta}
					\end{subfigure}
					\caption{Sensitivity of the optimal policy to insurance parameters.}
					\label{FIG:optimal_policy_sensi_insurance_policy}
				\end{figure}
				
			\paragraph*{Impact of financing and valuation parameters} ~\\
				The discount rate $\rho$ and the capital injection cost $\delta$ capture financial and valuation effects.  
				A higher discount rate reduces the present value of future profits, leading the insurer to liquidate earlier rather than maintaining operations with limited expected returns.  
				This translates into a contraction of the continuation region and an expansion of the dividend area.  
				Conversely, increasing the injection cost $\delta$ discourages recapitalisation and makes the firm more reluctant to support temporary losses, thereby enlarging the liquidation region and shrinking the domain where capital injections are optimal.
				\begin{figure}[H]
					\centering
					\begin{subfigure}[t]{0.48\linewidth}
						\centering
						\includegraphics[width=\linewidth]{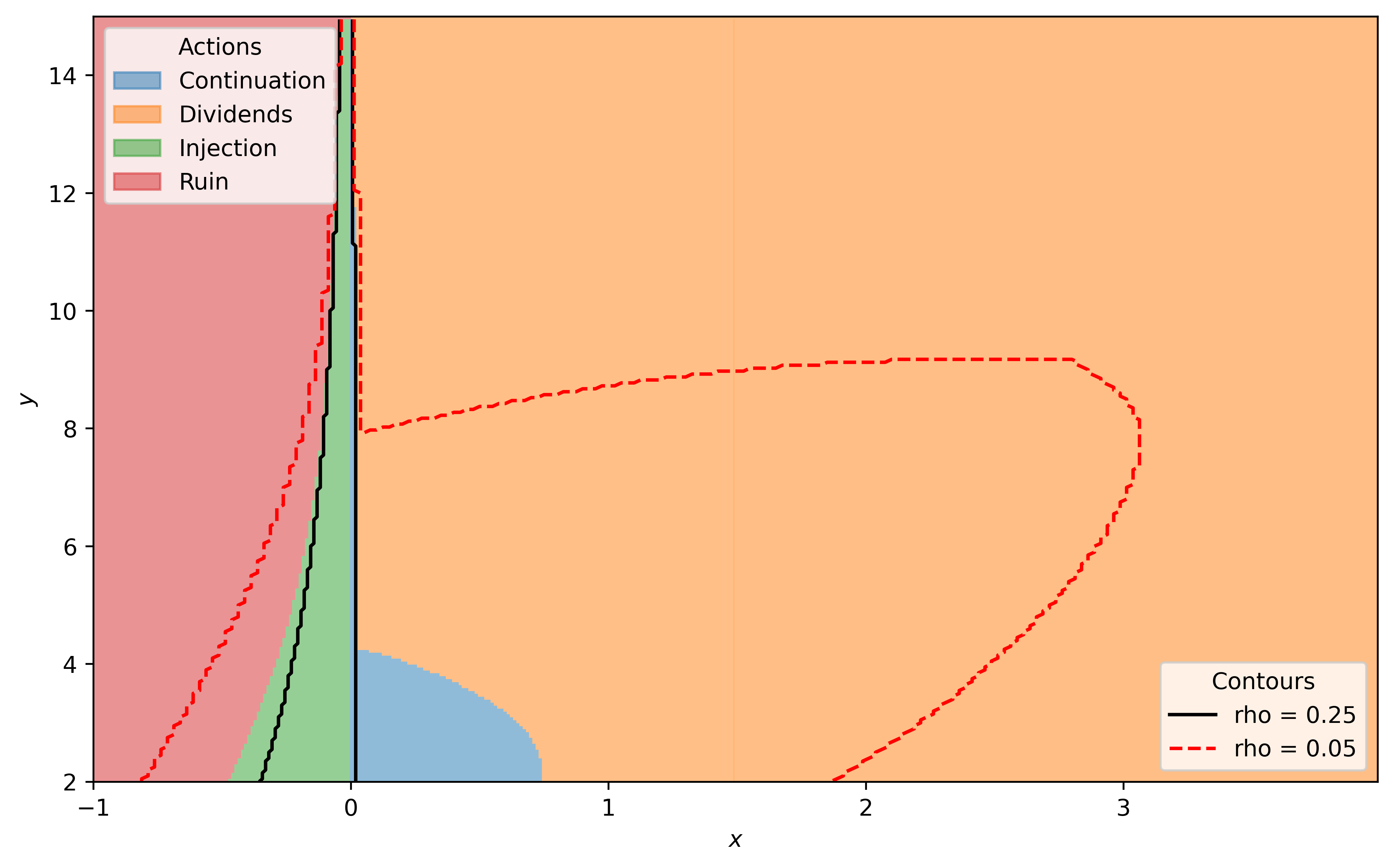}
						\caption{Sensitivity to the discount rate $\rho$.}
						\label{FIG:policy_sensi_to_rho}
					\end{subfigure}
					\hfill
					\begin{subfigure}[t]{0.48\linewidth}
						\centering
						\includegraphics[width=\linewidth]{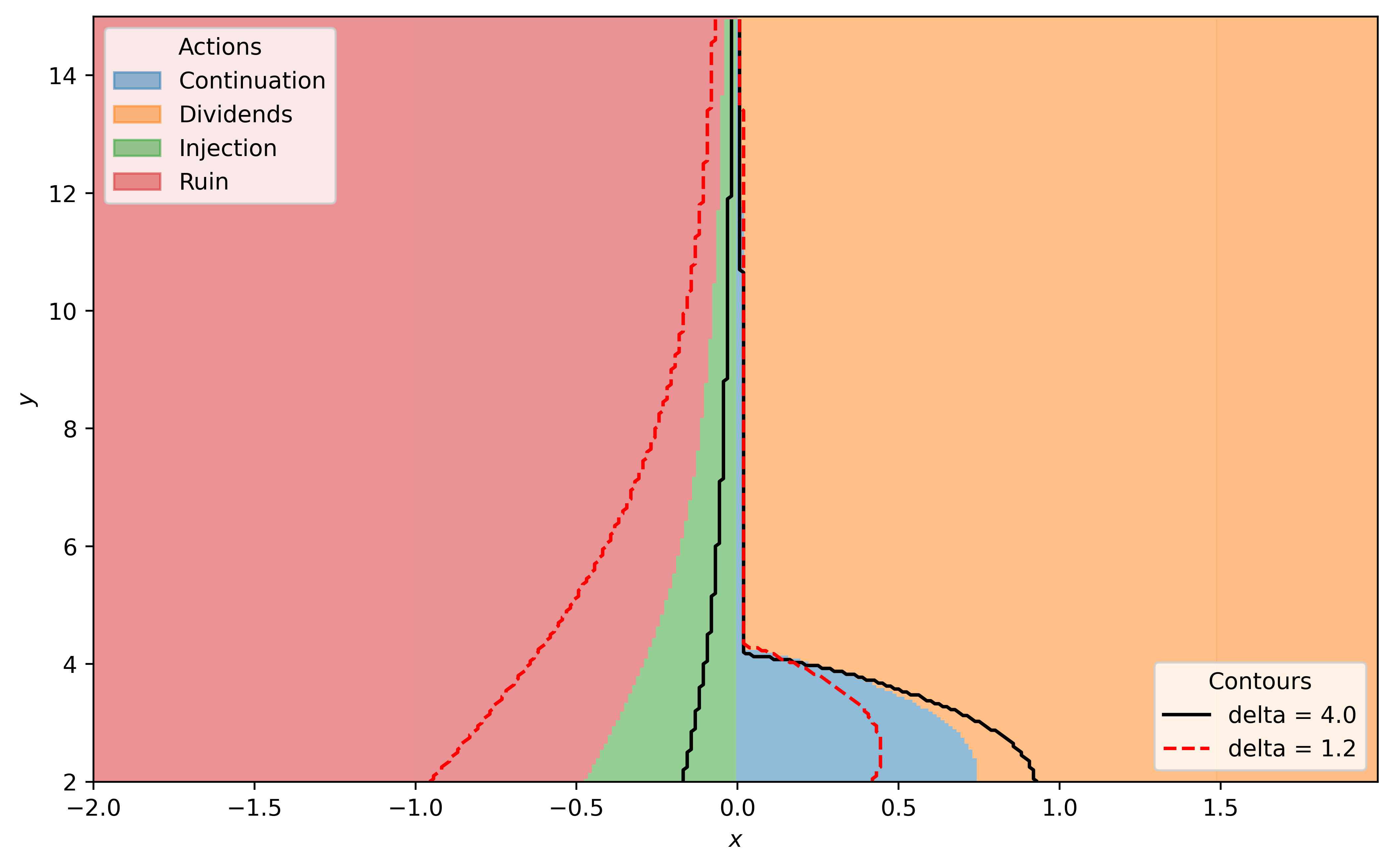}
						\caption{Sensitivity to the injection cost $\delta$.}
						\label{FIG:policy_sensi_to_delta}
					\end{subfigure}
					\caption{Sensitivity of the optimal policy to financing and valuation parameters.}
					\label{FIG:optimal_policy_sensi_conjoncture}
				\end{figure}
				
		\subsubsection{Ruin analysis under the optimal policy}
				
			We now investigate the ruin behaviour of the controlled surplus process and compare it with the uncontrolled case.
			While the previous sections characterise the optimal policy and its sensitivity to model parameters, the present analysis focuses on the distributional properties of the ruin time $T^\alpha$ under the optimal strategy.
					
\begin{figure}[H]
				\centering
				\begin{subfigure}[t]{0.32\linewidth}
					\centering
					\includegraphics[width=\linewidth]{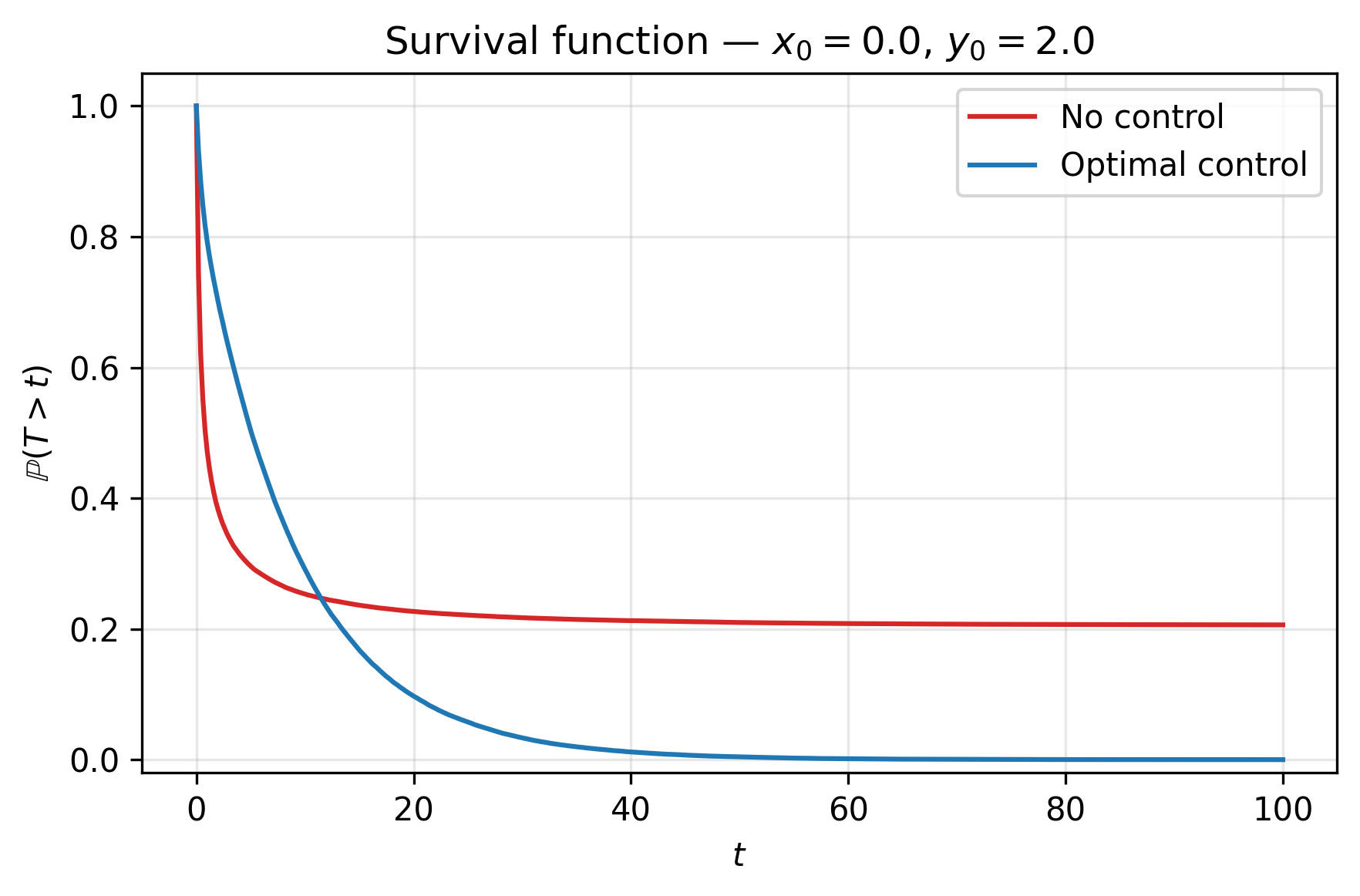}
					\caption{$x_0 = 0$, $y_0 = 2$.}
				\end{subfigure}
				\hfill
				\begin{subfigure}[t]{0.32\linewidth}
					\centering
					\includegraphics[width=\linewidth]{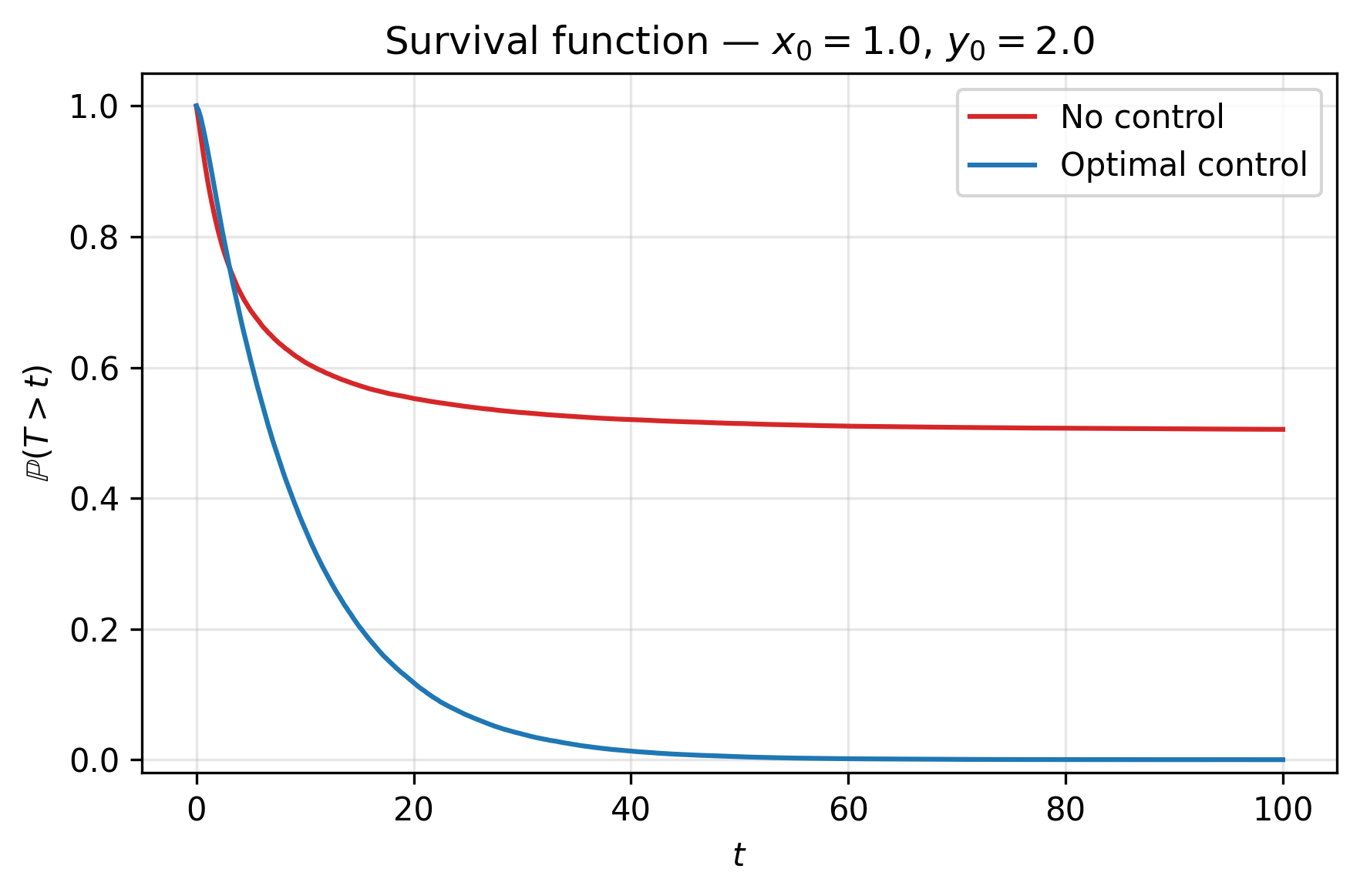}
					\caption{$x_0 = 1$, $y_0 = 2$.}
				\end{subfigure}
				\hfill
				\begin{subfigure}[t]{0.32\linewidth}
					\centering
					\includegraphics[width=\linewidth]{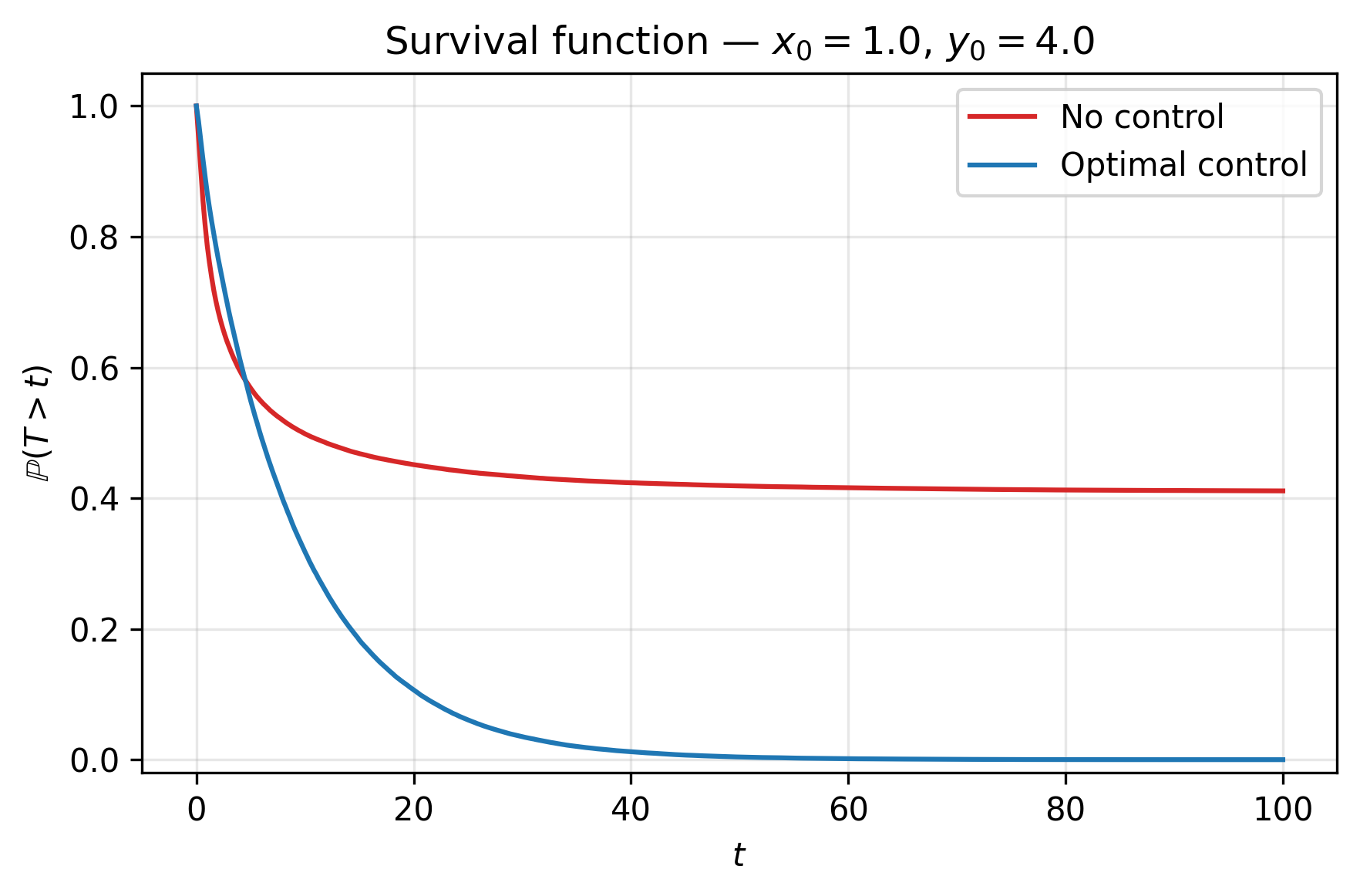}
					\caption{$x_0 = 1$, $y_0 = 4$.}
				\end{subfigure}
				\caption{Survival functions under no control and optimal control for several initial states.}
				\label{FIG:survival_comparison}
			\end{figure}

\begin{figure}[H]
				\centering
				\begin{subfigure}[t]{0.32\linewidth}
					\centering
					\includegraphics[width=\linewidth]{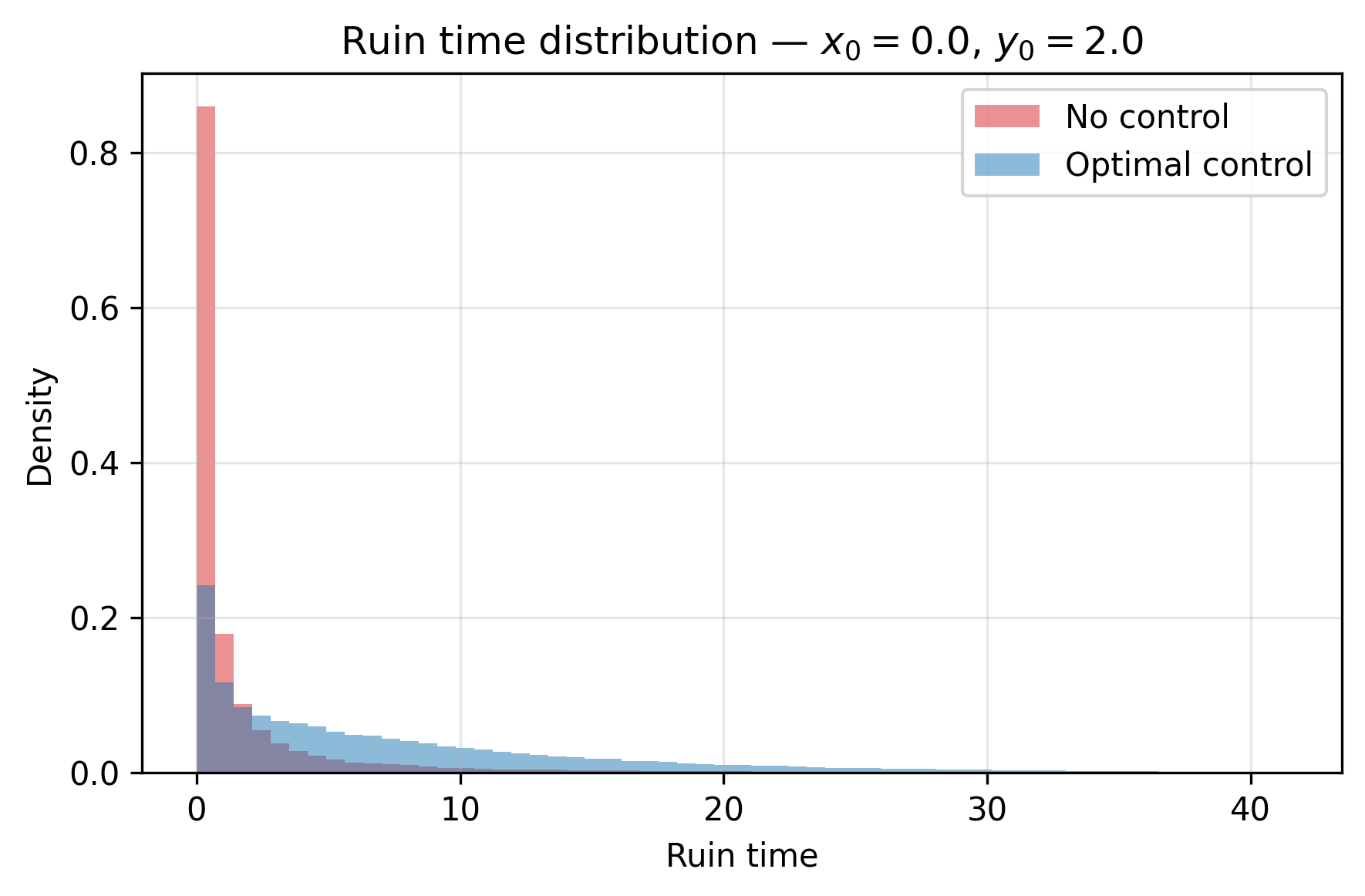}
					\caption{$x_0 = 0$, $y_0 = 2$.}
				\end{subfigure}
				\hfill
				\begin{subfigure}[t]{0.32\linewidth}
					\centering
					\includegraphics[width=\linewidth]{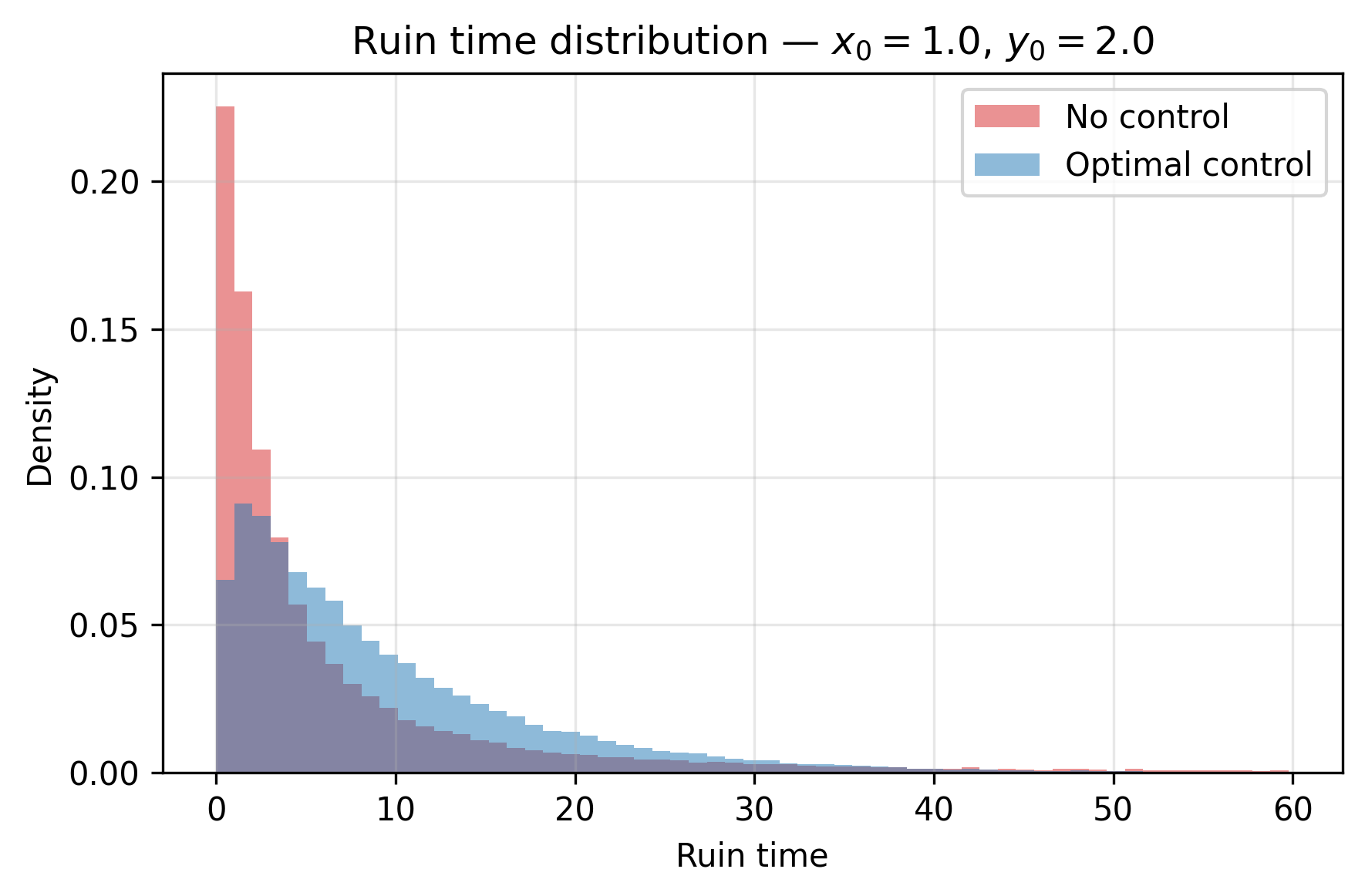}
					\caption{$x_0 = 1$, $y_0 = 2$.}
				\end{subfigure}
				\hfill
				\begin{subfigure}[t]{0.32\linewidth}
					\centering
					\includegraphics[width=\linewidth]{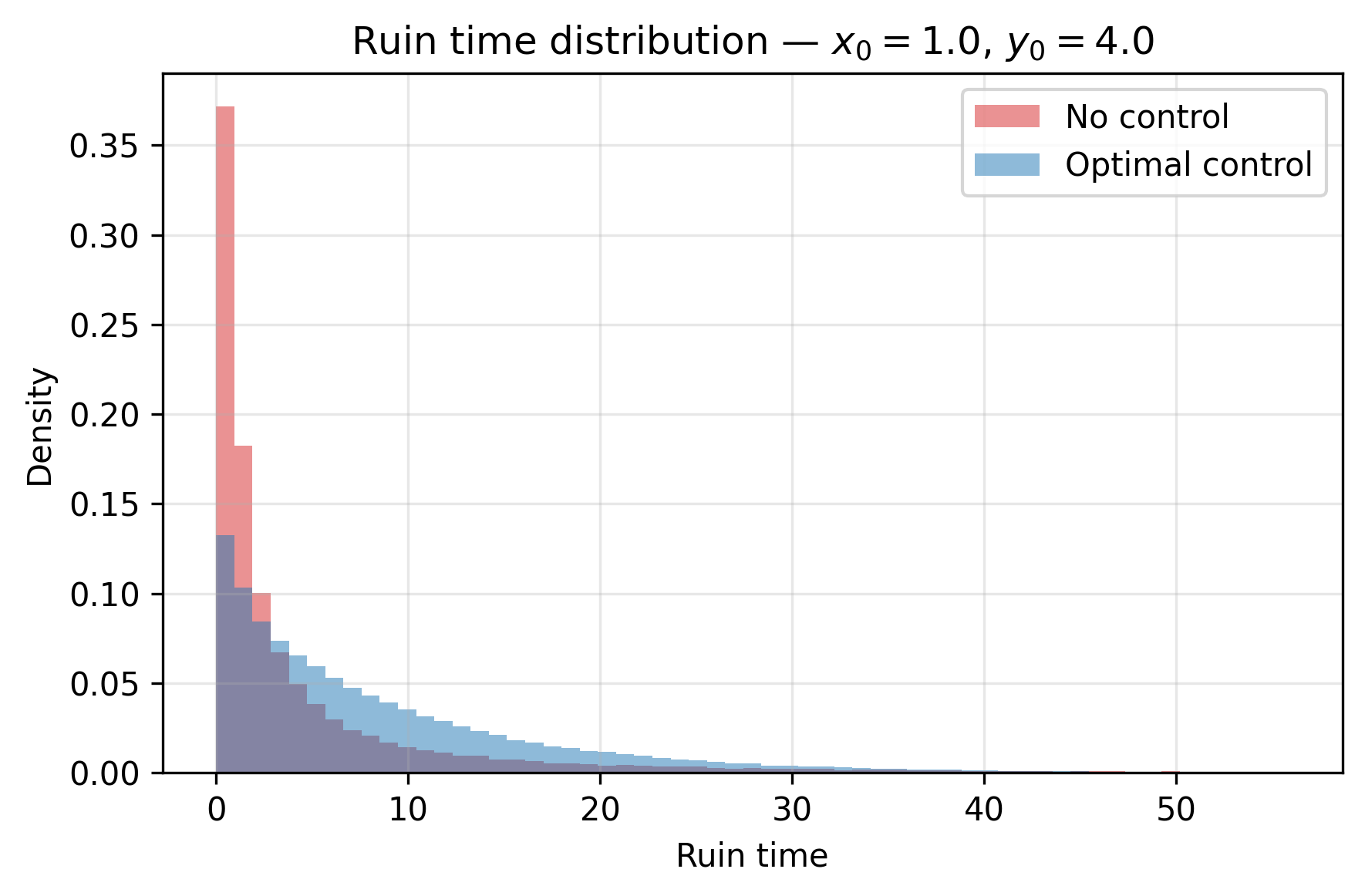}
					\caption{$x_0 = 1$, $y_0 = 4$.}
				\end{subfigure}
				\caption{Empirical distribution of ruin times conditional on ruin, under no control and optimal control.}
				\label{FIG:ruin_histograms}
			\end{figure}

			\paragraph*{Controlled versus uncontrolled ruin} ~\\
			We simulate $100{,}000$ trajectories of the surplus process under the baseline parameter configuration (Table~\ref{TAB:model_parameters_PDE}), both with and without the optimal control, for several initial states $(x_0, y_0)$.
			The simulation is run up to a horizon $T = 100$.
			Figure~\ref{FIG:survival_comparison} displays the corresponding survival functions.
			Under the optimal policy, ruin tends to occur earlier than in the uncontrolled case, confirming that the optimal strategy is to prioritise shareholder value extraction over survival.
			The dividend distribution mechanism actively depletes the surplus, increasing the insurer's exposure to adverse claim sequences.
			Capital injections partially counterbalance this effect by preventing ruin when the surplus falls below zero but remains above the threshold $\kappa^\star$.
			For higher initial intensities, the injection threshold approaches zero (cf.\ Corollary~\ref{COR:yasymptotic}), further reducing the scope for
			recapitalisation and accelerating ruin.
			Table~\ref{TAB:ruin_baseline} reports the associated ruin statistics.
			The ruin time moments and histograms below are computed conditional on ruin occurring before $T$.
			Under the optimal policy, all trajectories eventually reach ruin, 
			whereas a significant fraction survives indefinitely in the uncontrolled case when the net profit condition is satisfied (see Remark~\ref{RMK:net_profit_condition}). 
			This is a direct consequence of the dividend barrier, which prevents the surplus from accumulating without bound. 
			However, conditional on ruin occurring, the optimal policy delays it: capital injections sustain the firm through adverse episodes that would cause immediate ruin in the absence of control.
			\begin{table}[H]
				\centering
				\begin{tabular}{cc|ccc|ccc}
					\toprule
					& & \multicolumn{3}{c|}{No control} & \multicolumn{3}{c}{Optimal control} \\
					$x_0$ & $y_0$ & Ruin \% & $\mathbb{E}[T^\alpha|\text{ruin}]$ & Med.$|\text{ruin}$ & Ruin \% & $\mathbb{E}[T^\alpha|\text{ruin}]$ & Med.$|\text{ruin}$ \\
					\midrule
					$0.0$ & $2.0$ & $0.806$ & $2.7$ & $0.5$ & $1.000$ & $7.7$ & $4.8$ \\
					$1.0$ & $2.0$ & $0.524$ & $8.0$ & $3.2$ & $1.000$ & $9.2$ & $6.5$ \\
					$1.0$ & $4.0$ & $0.610$ & $5.8$ & $1.8$ & $1.000$ & $8.3$ & $5.6$ \\
					\bottomrule
				\end{tabular}
				\caption{Ruin statistics under no control and under the optimal policy for the baseline configuration.}
				\label{TAB:ruin_baseline}
			\end{table}
			Figure~\ref{FIG:ruin_histograms} shows the empirical distribution of ruin times conditional on ruin, for both regimes. In the uncontrolled case, the distribution is concentrated near the origin: trajectories that do reach ruin tend to do so early, driven by an unfavourable claim sequence before the surplus can accumulate. 
			The remaining trajectories, not represented in the histogram, survive indefinitely. Under the optimal policy, the distribution is shifted to the right and more spread out: capital injections sustain the firm through episodes that would otherwise be fatal, but the continuous extraction of dividends ensures that ruin eventually occurs on every trajectory.
			
\begin{figure}[H]
				\centering
				\begin{subfigure}[t]{0.48\linewidth}
					\centering
					\includegraphics[width=\linewidth]{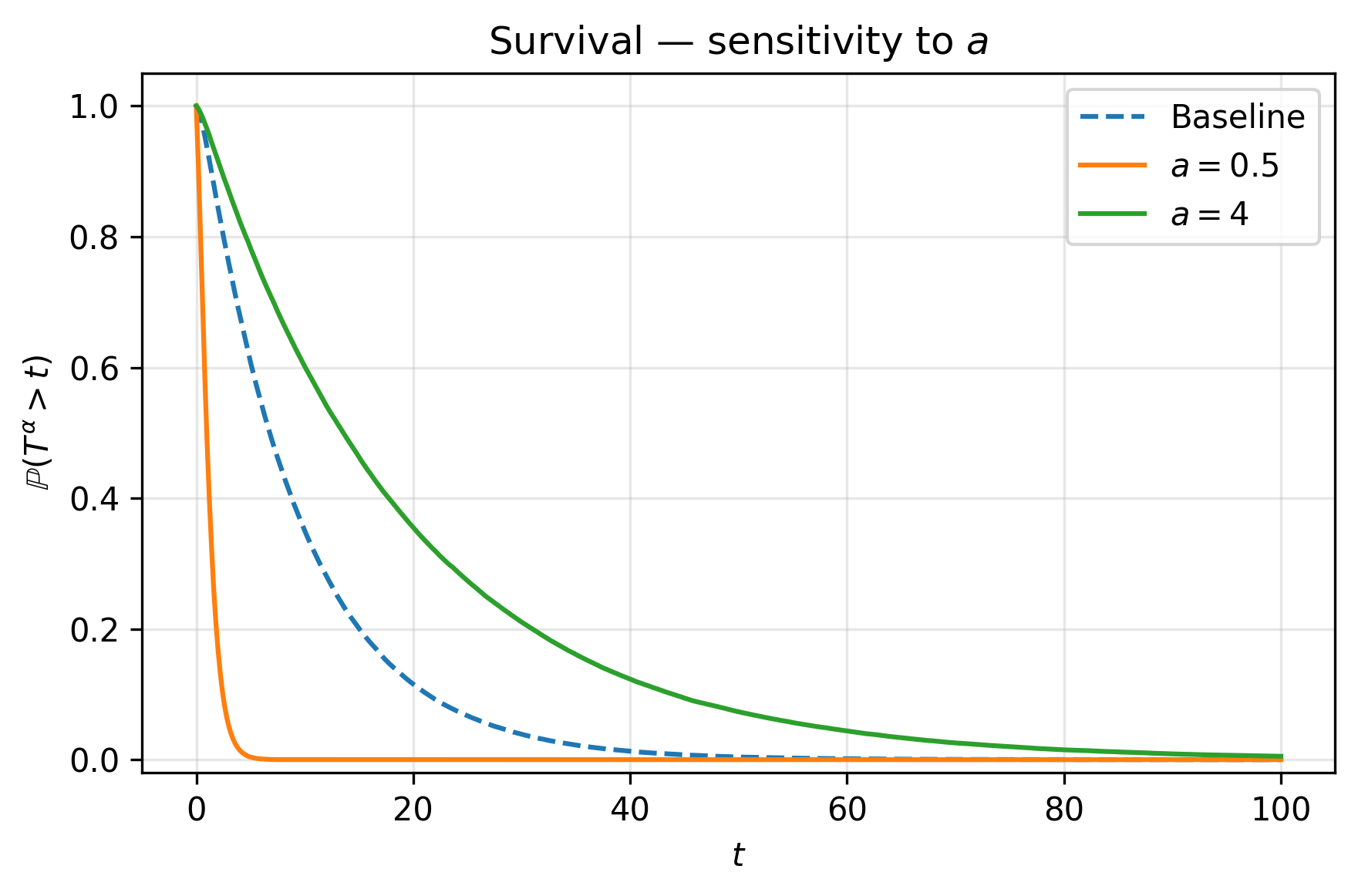}
					\caption{Sensitivity to $a$.}
				\end{subfigure}
				\hfill
				\begin{subfigure}[t]{0.48\linewidth}
					\centering
					\includegraphics[width=\linewidth]{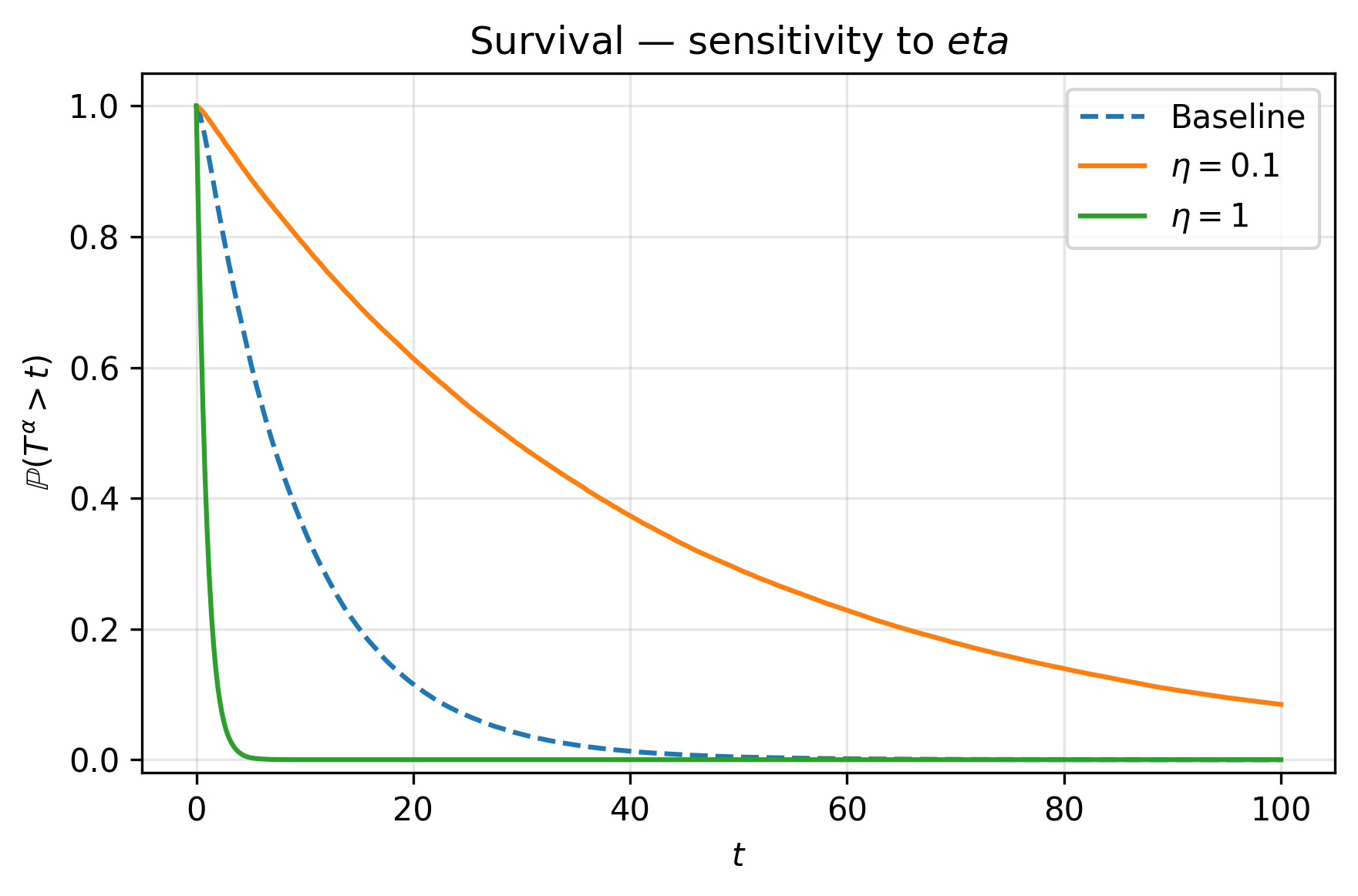}
					\caption{Sensitivity to $\eta$.}
				\end{subfigure}
				\\[0.3cm]
				\begin{subfigure}[t]{0.48\linewidth}
					\centering
					\includegraphics[width=\linewidth]{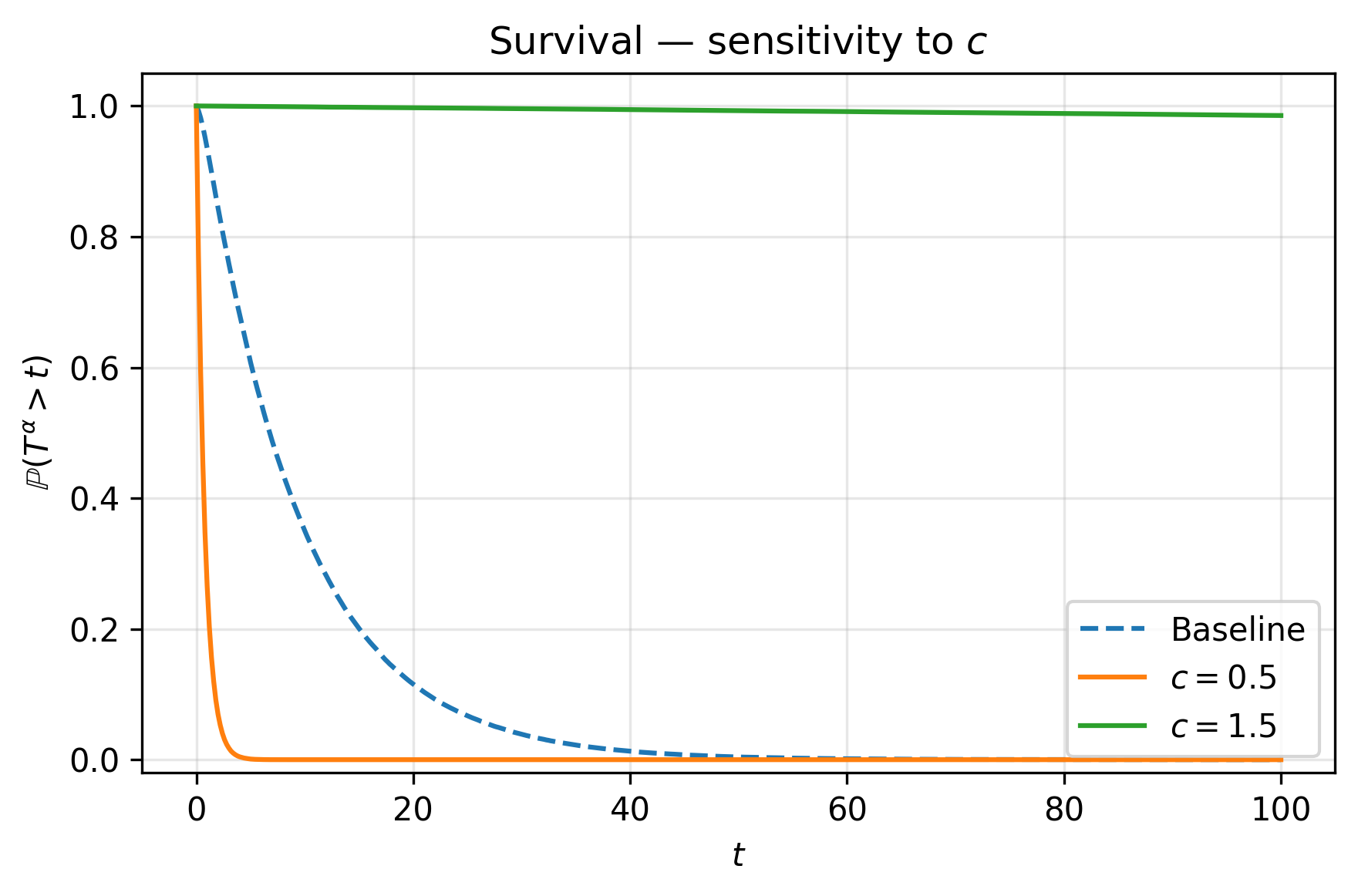}
					\caption{Sensitivity to $c$.}
				\end{subfigure}
				\hfill
				\begin{subfigure}[t]{0.48\linewidth}
					\centering
					\includegraphics[width=\linewidth]{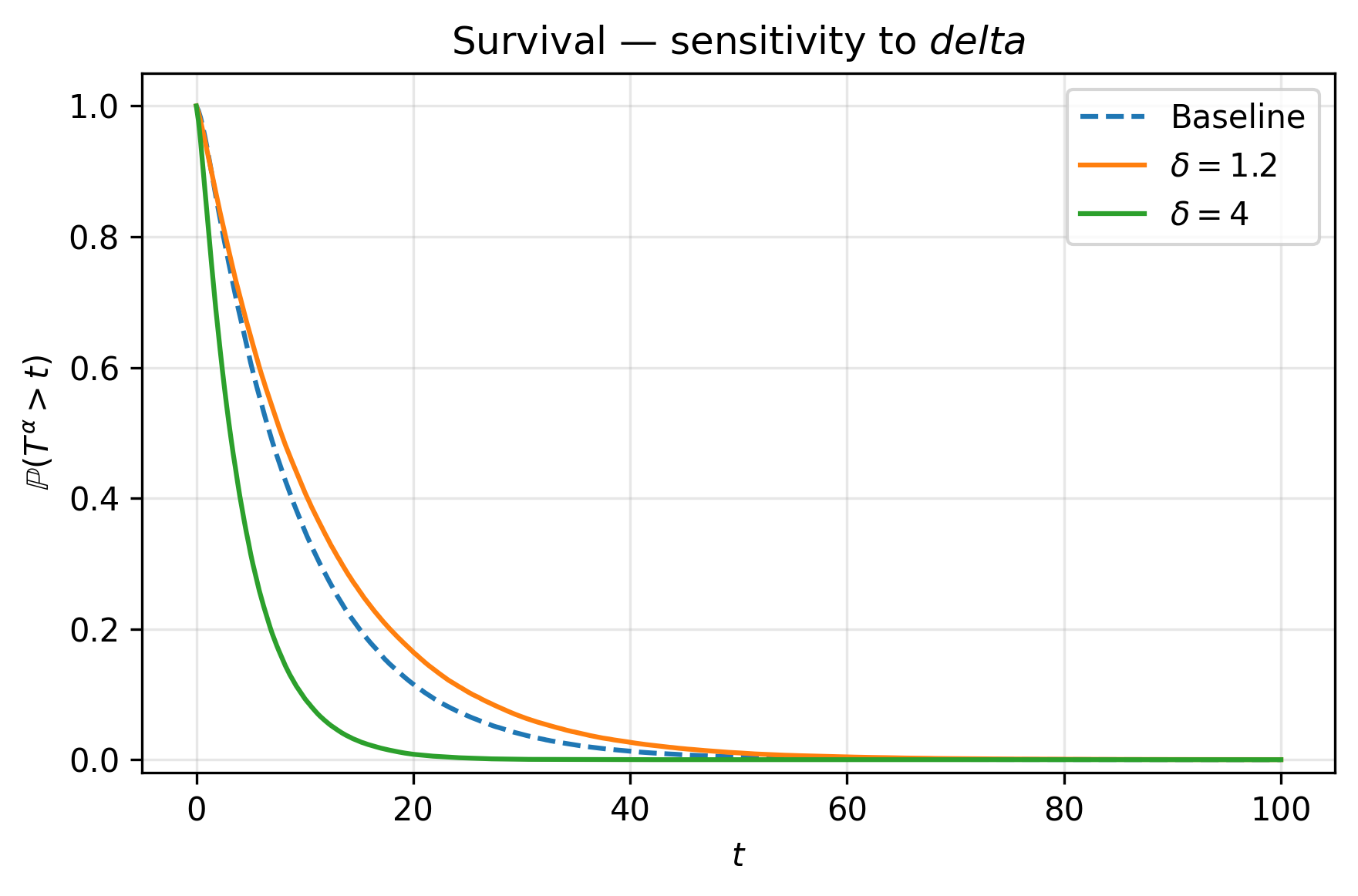}
					\caption{Sensitivity to $\delta$.}
				\end{subfigure}
				\caption{Survival functions under the optimal policy for selected parameter variations.}
				\label{FIG:sensi_survival}
			\end{figure}

			\paragraph*{Sensitivity of ruin metrics to model parameters} ~\\
			We now examine how the ruin behaviour under the optimal policy responds to changes in selected model parameters.
			For each variation, we simulate trajectories starting from $(x_0, y_0) = (1, 2)$ under the corresponding optimal strategy.
			Figure~\ref{FIG:sensi_survival} displays the survival functions, with the baseline shown as a dashed reference.
			The parameters with the strongest impact on ruin are those governing the net profit condition~\eqref{EQN:net_profit_condition}. 
			Increasing $\eta$ or decreasing $c$ deteriorates the long-run profitability and leads to substantially earlier ruin, even under the optimal policy. 
			The mean-reversion rate $a$ has a similar effect: slower reversion prolongs high-intensity episodes and accelerates ruin. In contrast, the injection cost $\delta$ primarily affects the shape of the optimal strategy but has a more moderate effect on ruin times.
			
\section{Reinforcement learning estimate}
	\label{SEC:rl_resolution}
	In this section, we develop a numerical method based on policy optimisation techniques from reinforcement learning to solve the control problem~\eqref{EQN:optimal_control_problem}. 
	While the finite-difference scheme of Section~\ref{SEC:pde_resolution} provides an accurate benchmark, it relies on a grid discretisation of the state space whose computational cost grows rapidly with the dimension. 
	The reinforcement learning approach developed below does not require a grid discretisation of the state space and scales naturally to higher-dimensional or path-dependent extensions of the model, making it a complementary tool for problems where PDE-based methods become impractical.
	
	\subsection{Discrete-time reformulation of the control problem}
		
		\subsubsection*{Formulation as a Markov Decision Process}
			We begin by reformulating the problem within a general MDP framework, following the approach of \cite{Hamdouche23} for control problems with random exit times. 
			Let $(S_t^\alpha)_{t \ge 0}$ denote the controlled state process taking values in a domain $\Sc \subset \mathbb{R}^d$, and let $\alpha$ be an admissible control with values in a subset of $\mathbb{R}^m$. The process evolves until the random exit time
			\begin{equation*}
				T^\alpha := \inf \{ t \geq 0 : S_t^\alpha \notin \Oc \}, 
			\end{equation*}
			where $\Oc \subset \Sc$ is an open domain.
			\begin{Remark}
				In our setting, the state variable is $S_t = (X_t, \lambda_t)$, where $X_t$ denotes the surplus process and $\lambda_t$ the Hawkes intensity. The control is $\alpha_t = (Z_t, K_t)$, consisting of cumulative dividends and capital injections.
			\end{Remark}
			The performance criterion is defined in terms of a running reward function $f$ and a terminal reward function $g$.  
			Given an initial state $S_0 = s$, the expected return under a control $\alpha$ is
			\begin{equation}
				\label{EQN:objective_function_MDP}
				J_\alpha(s)
				= \E\!\left[
				\int_0^{T^\alpha} e^{-\rho t} f(S_t^\alpha, d\alpha_t) 
				+ g(S_{T^\alpha}^\alpha)
				\right].
			\end{equation}
			The associated value function is
			\begin{equation*}
				v(s) = \sup_{\alpha \in \Ac(s)} J_\alpha(s) .
			\end{equation*}
			\begin{Remark}
				In our model, there is no terminal reward, i.e., $g \equiv 0$.  
				The running reward is
				\begin{equation*}
					f(S_t^\alpha, d\alpha_t) = dZ_t - \delta dK_t,
				\end{equation*}
				reflecting dividend payments and penalised capital injections.
			\end{Remark}
			We allow for general controlled state dynamics, potentially involving drift, diffusion, jumps, and control actions. A typical form is
			\begin{equation*}
				dS_t = \mu(S_t)dt + \sigma(S_t)dW_t + \eta(S_t)dN_t + \sum_{i=1}^m d\alpha_t^i,
			\end{equation*}
			where $W_t$ is a Brownian motion and $N_t$ a jump process (e.g., a Hawkes process).  
			In our model, the state $S_t = (X_t,\lambda_t)$ evolves with deterministic drift and jump-driven increments: $\lambda_t$ follows Hawkes dynamics, while $X_t$ is affected by premium inflows, claim jumps, and the control $(Z_t,K_t)$.
			
		\subsubsection*{Discretisation and finite-horizon MDP approximation}
			Let $\T = \{ t_0 = 0 < t_1 < \cdots < t_N \}$ be a uniform time grid with step size $h>0$.  
			The state space is $\Sc \subset \mathbb{R}^d$, and we denote by $s \in \Sc$ the initial state.  
			At each state $s_i$, the set of admissible controls is $\Ac(s_i) \subset \mathbb{R}^m$.  
			
			We consider the discretised controlled process $(S_{t_i})_{i=0}^N$, where the transition from $S_{t_i}=s_i$ to $S_{t_{i+1}}$ under control $a \in \Ac(s_i)$ is specified by a transition kernel
			\begin{equation*}
				p(\cdot \mid t_i, s_i, a),
			\end{equation*}
			that is, $p(\cdot \mid t_i, s_i, a)$ is the law of $S_{t_{i+1}}$ given $(S_{t_i},a)$.
			\begin{Definition}[Randomised policy]
				A randomised policy is a measurable transition kernel
				\begin{equation*}
					\pi : (t_i, s_i) \in \T \times \Sc \longmapsto 
					\pi(\cdot \mid t_i, s_i) \in \mathcal{P}(\Ac(s_i)),
				\end{equation*}
				assigning to each state a probability distribution over admissible actions.  
				We write $\alpha \sim \pi$ for the random control sequence generated under~$\pi$.
			\end{Definition}
			We denote by $\Pi_h$ the set of all admissible discrete-time randomised policies.  
			Under $\pi \in \Pi_h$, the controlled state process is $(S_{t_i}^\pi)_{i=0}^N$.  
			The discrete-time exit time is defined as
			\begin{equation*}
				\tau := \inf \{ t_i \in \T : S_{t_i}^\pi \notin \Oc \}.
			\end{equation*}
			We introduce the corresponding exit index
			\begin{equation*}
				N(\tau) := \inf \{ i \in \{0,\ldots,N\} : S_{t_i}^\pi \notin \Oc \}.
			\end{equation*}
			To obtain a discrete-time counterpart of the objective \eqref{EQN:objective_function_MDP}, 
			let $ A_{t_{i+1}}^\pi := \alpha_{t_{i+1}} - \alpha_{t_i}$ denote the action increment on $[t_i,t_{i+1}]$. 
			The expected cumulative reward under a policy $\pi$ is then
			\begin{equation}
				\label{EQN:value_function_MDP}
				J(\pi)
				= \E_{\alpha \sim \pi}\!\left[
				\sum_{i=0}^{N(\tau)-1} 
				f(S_{t_i}^\pi, A_{t_{i+1}}^\pi)
				+ g(S_\tau^\pi)
				\right],
			\end{equation}
			where $f$ is the instantaneous reward function and $g$ the terminal reward.
	
	\subsection{Policy gradient estimators}
	    \subsubsection*{Gradient representations for policy optimisation}
	    	Policy optimisation methods developed in the reinforcement learning literature provide an alternative way to approximate the solution of stochastic control problems. 
	    	In this framework, the control is modelled through a parametrised family of stochastic policies 
	    	$\{\pi_\theta : \theta \in \mathbb{R}^p\}$, where each policy assigns to a state a probability distribution over actions. 
	    	Such distributions are typically represented by neural networks, whose parameters depend on the current state.  
	    	Sampling actions from these distributions yields unbiased gradient estimators of the expected return, giving rise to policy gradient algorithms.
	    	
	    	We now introduce a formal definition of a parametrised stochastic policy.
	    	\begin{Definition}[Parametrised stochastic policy]
	    		Let $\theta \in \mathbb{R}^p$.  
	    		A stochastic policy $\pi_\theta$ is said to be parametrised if, for each $(t_i,s_i)$, it admits a density with respect to a reference measure $\nu$ on $\Ac(s_i)$:
	    		\[
	    		\pi_\theta(da \mid t_i, s_i) = \rho_\theta(t_i, s_i, a)\nu(da),
	    		\]
	    		where $\rho_\theta : \T \times \Sc \times \Ac \to (0,+\infty)$ is a measurable function.
	    	\end{Definition}
	    	We restrict attention to parametrised policies of the form $\pi_\theta$, and the objective becomes to optimise the parameter $\theta \in \mathbb{R}^p$ so as to maximise the discrete-time functional~\eqref{EQN:value_function_MDP}.
			\begin{Theorem}[Objective function gradient]
				\label{TH:gradient_classic_representation}
				Let $\theta \in \R^p$ and $\pi_\theta$ be a randomised parametrised policy. Then, the gradient of (\ref{EQN:value_function_MDP}) with respect to $\theta$ is given by:
				\begin{equation}
					\label{EQN:policy_gradient_representation}
					\nabla_\theta J(\pi_\theta) = \mathbb{E}_{\alpha \sim \pi_\theta}\left[\left(\Sum{i = 0}{N(\tau) - 1} f(S_{t_i}^{\pi_\theta}, A_{t_{i+1}}^{\pi_\theta}) + g(S_{\tau}^{\pi_\theta})\right)\left(\Sum{i=0}{N(\tau) - 1} \nabla_\theta \log(\rho_\theta(t_i, S_{t_i}^{\pi_\theta}, A_{t_{i+1}}^{\pi_\theta}))\right)\right],
				\end{equation}
				where we recall that $A_{t_{i+1}}^{\pi_\theta} = \alpha_{t_{i+1}}^{\pi_\theta} - \alpha_{t_i}^{\pi_\theta}$
			\end{Theorem}
			
			\begin{proof}				
				Recall that:
				\begin{equation*}
					J(\pi_\theta) = \mathbb{E}_{\alpha \sim \pi_\theta}\left[\Sum{i = 0}{N(\tau) - 1} f(S_{t_i}^{\pi_\theta}, A_{t_{i+1}}^{\pi_\theta}) + g(S_{\tau}^{\pi_\theta})\right]
				\end{equation*}
				The proof relies on  the arguments presented in the work by Hamdouche et al. \cite{Hamdouche23} . In our setting we need to increase the dimension of the dynamics. Hence, we consider the process $Y$ defined as follows:
				\begin{equation*}
					Y_t = \Int{0}{t} e^{-\rho s} f(S_s^\alpha, d\alpha_s) + g(S_t^\alpha),\quad\textrm{for }t\geq 0
				\end{equation*}
				As the process $\tilde{S}_t = (S	_t, Y_t)_{t\geq 0}$ is Markovian and we can apply Theorem (2.1) of Hamdouche et al. to get that			
				\begin{equation*}
					\nabla_\theta J(\pi_\theta) = \mathbb{E}_{\alpha \sim \pi_\theta}\left[\left(\Sum{i = 0}{N(\tau) - 1} f(S_{t_i}^{\pi_\theta}, A_{t_{i+1}}^{\pi_\theta}) + g(S_{\tau}^{\pi_\theta})\right)\left(\Sum{i=0}{N(\tau) - 1} \nabla_\theta \log(\rho_\theta(t_i, S_{t_i}^{\pi_\theta},  A_{t_{i+1}}^{\pi_\theta}))\right)\right]
				\end{equation*}
			\end{proof}
			
			The representation of Theorem~\ref{TH:gradient_classic_representation} expresses the gradient of the performance functional in terms of the cumulative realised reward and the score of the policy.  
			An alternative and often more stable estimator can be obtained by exploiting the dynamic programming structure of the value process.  
			To this end, following \cite{Hamdouche23}, we introduce a dynamic version of the performance functional under the policy~$\pi_\theta$. \\
			For each index $i \in \{0,\dots,N\}$ and state $s \in \Sc$, define
			\begin{equation*}
				v_i^\theta(s)
				:= \E_{\alpha \sim \pi_\theta}\!\left[
				\sum_{j=i}^{N(\tau_i)-1}
				f(S_{t_j}^{\pi_\theta}, A_{t_{j+1}}^{\pi_\theta})
				+ g\!\left(S_{\tau_i}^{\pi_\theta}\right)
				\middle|
				S_{t_i}^{\pi_\theta}=s
				\right],
			\end{equation*}
			where the local exit time is
			\begin{equation*}
				\tau_i
				:= \inf\{ t_j \in \T : t_j \ge t_i, S_{t_j}^{\pi_\theta} \notin \Oc \}
				\wedge t_N.
			\end{equation*}
			Clearly, $v_N^\theta(s)=g(s)$, and $v_i^\theta(s)=g(s)$ for all $i<N$ whenever
			$s \notin \Oc$.  
			Moreover, by the discrete-time dynamic programming principle for $s \in \Oc, i=0,\dots,N-1$,
			\begin{equation*}
				v_i^\theta(s)
				= \E_{\alpha \sim \pi_\theta}\!\left[
				v_{i+1}^\theta\!\left(S_{t_{i+1}}^{\pi_\theta}\right)
				\middle|
				S_{t_i}^{\pi_\theta}=s
				\right] .
			\end{equation*}
			
			\begin{Theorem}[Martingale representation] 
				\label{TH:gradient_martingale_representation}
				We have:
				\begin{equation}
					\nabla_\theta J(\pi_\theta) = \mathbb{E}_{\alpha \sim \pi_\theta}\left[\Sum{i=0}{N(\tau) - 1} v_{i+1}^\theta(S_{t_{i+1}}^{\pi_\theta}) \nabla_\theta \log ( \rho_\theta (t_i, S_{t_i}^{\pi_\theta}, A_{t_{i+1}}^{\pi_\theta} ))\right]
				\end{equation}
			\end{Theorem}
			
			\begin{proof}
				For a trajectory controlled by $\pi_\theta$, define the cumulative reward process on the discrete grid by
				\begin{equation*}
					Y_{t_0} := 0,\qquad
					Y_{t_{i+1}} := Y_{t_i} + f\big(S_{t_i}^{\pi_\theta}, A_{t_{i+1}}^{\pi_\theta}\big),
					\quad i = 0,\dots,N(\tau)-1,
				\end{equation*}
				and set $Y_\tau := Y_{t_{N(\tau)}} + g\big(S_\tau^{\pi_\theta}\big)$.
				Then $J(\pi_\theta) = \E_{\alpha \sim \pi_\theta}[Y_\tau]$, and the augmented process $\tilde S_{t_i} := (S_{t_i}^{\pi_\theta}, Y_{t_i})$ is Markovian. Applying the results in Hamdouche et al.~\cite{Hamdouche23} to $(\tilde S_{t_i})_{i\ge 0}$ yields the gradient representation stated in Theorem~\ref{TH:gradient_martingale_representation}.
			\end{proof}
			
			\begin{Remark}
				An equivalent expression based on the temporal differences of the value function is given by:
				\begin{equation*}
					\nabla J(\pi_\theta) = \mathbb{E}_{\alpha \sim \pi_\theta}\left[\Sum{i=0}{N(\tau) - 1} \left( v_{i+1}^\theta(S_{t_{i+1}}^{\pi_\theta}) - v_i^\theta(S_{t_i}^{\pi_\theta}) \right) \nabla_\theta \log ( \rho_\theta (t_i, S_{t_i}^{\pi_\theta}, A_{t_{i+1}}^{\pi_\theta} ))\right]
				\end{equation*}
				This form is particularly relevant in actor-critic methods where $v_i^\theta$ is replaced by a learned critic.
			\end{Remark}
			
			After time discretisation with step $h>0$, the controlled process $(X_t, \lambda_t)$
			induces a Markov decision process with continuous action space.  
			Rather than discretising the actions, we restrict attention to a parametrised class
			of stochastic policies, typically implemented through neural networks.  
			From a theoretical perspective, the work of Kushner and Dupuis \cite{KushnerDupuis2001}
			shows that, when the full admissible action space is retained, the value functions of
			the discrete-time control problems converge to their continuous-time counterpart as
			$h \to 0$.  
			The use of parametrised stochastic policies introduces a second level of approximation: the optimisation is now restricted to a subset of all admissible randomised controls.
			This additional approximation does not affect the consistency of the time
			discretisation itself, but it may prevent the algorithm from attaining the true
			optimal value if the optimal policy lies outside the chosen parametrised class.
		
		\subsubsection*{Gradient-based learning algorithm}
		
			We now leverage Theorems~\ref{TH:gradient_classic_representation} and~\ref{TH:gradient_martingale_representation} to design learning algorithms aimed at approximating optimal policies.
			Our first method is a direct policy gradient algorithm based on Theorem~\ref{TH:gradient_classic_representation}. This approach corresponds to an extension of the well-known REINFORCE algorithm introduced by Sutton~\cite{Sutton1998}, and has been adapted in recent works. 
			The algorithm proceeds as follows: the policy is initialised and used to generate a collection of sample paths. For each path, the cumulative reward and the log-probabilities of the actions taken are recorded. These are then used to compute a Monte Carlo estimate of the gradient, which serves to update the policy parameters via stochastic gradient ascent.
			
			\begin{algorithm}[H]
				\caption{Policy gradient algorithm}
				\label{ALG:off_line_policy_gradient}
				\begin{algorithmic}
					\STATE Number of episodes $E$, number of Monte Carlo trajectories $K$ and learning rate $\eta$
					\STATE Initialize policy $\pi_\theta$ with its parameters $\theta$ 
					\FOR{epoch $e = 1, ..., E$}
					\FOR{trajectory $k = 1, ..., K$}
					\STATE Apply current policy up to the end of trajectory $N(\tau_k): \{(s_i, a_i)\}_{i=0}^{\tau_k}$
					\STATE Calculate total reward and log probabilities:
					\STATE \quad $G_k = \Sum{i=0}{N(\tau_k) - 1} f(s_{t_i}, a_{t_i}) + g(s_{\tau_k})$
					\STATE \quad $\Lambda_k = \Sum{i=0}{N(\tau_k) - 1} \nabla_\theta \log(\rho_\theta (t_i, s_{t_i}, a_{t_i}))$
					\ENDFOR    
					\STATE Compute total loss and update policy parameters by gradient ascent:
					\STATE \quad $\theta \leftarrow \theta + \eta \frac{1}{K} \Sum{k=1}{K} G_k \Lambda_k$
					\ENDFOR
				\end{algorithmic}
			\end{algorithm}
			
			While the policy gradient algorithm is straightforward to implement and only requires that the policy admit a differentiable density, it does not rely on any value function approximation. This simplicity is one of its main advantages.
			However, a well-known drawback of REINFORCE-type methods is the high variance of the gradient estimators, which can lead to slow and unstable convergence. To address this issue, several variance-reduction techniques have been proposed. A common strategy is to subtract a baseline from the return: typically, an estimate of the value function. It helps reduce variance without introducing bias. This idea motivates the actor-critic methods presented in the next section.
			
		\subsubsection*{Actor-critic algorithm}
			
			The second approach is an actor-critic algorithm, based on the gradient formula provided by Theorem~\ref{TH:gradient_martingale_representation}. This method combines elements of both value-based and policy-based methods: the actor updates the policy, while the critic estimates the value function. This dual update often results in improved sample efficiency and convergence stability.
			We follow the methodology introduced in~\cite{Sutton1998} and adapted in~\cite{Hamdouche23}, using two neural networks: one for the policy $\pi_\theta$ (the actor) and one for the value function $\hat{q}_\omega$ (the critic).
			
			\begin{algorithm}[H]
				\caption{Off-line actor critic policy gradient algorithm}
				\label{ALG:actor_critic_gradient_desc}
				\begin{algorithmic}
					\STATE Number of episodes $E$, number of trajectories to use $K$ and $\eta_\theta$ and $\eta_\omega$ the learning rates
					\STATE Initialize policy and value function $\pi_\theta$ and $\hat{q}_\omega$ with their parameters $\theta$ and $\omega$
					\FOR{epoch $e = 1, ..., E$}
					\FOR{trajectory $k = 1, ..., K$}
					\STATE Apply current policy $\pi_\theta$ up to the end of trajectory $N(\tau_k): \{(s_{t_i}, a_{t_i})\}_{i=0}^{N(\tau_k)}$
					\STATE Calculate total advantage and log probabilities:
					\STATE \quad $\Phi_k = \Sum{i=0}{N(\tau_k) - 1} (\hat{q}_\omega(s_{t_{i+1}}) - \hat{q}_\omega(s_{t_i})) \nabla_\theta \log (\rho_\theta(t_i, s_{t_i}, a_{t_{i+1}}))$
					\STATE \quad $\Psi_k = \Sum{i=0}{N(\tau_k) - 1} (\hat{q}_\omega(s_{t_{i+1}}) - \hat{q}_\omega(s_{t_i})) \nabla_\omega \hat{q}_\omega(s_{t_i})$
					\ENDFOR    
					\STATE Compute total losses and update policy and value function parameters by gradient ascent:
					\STATE \quad $\theta \leftarrow \theta + \eta_\theta \frac{1}{K} \Sum{k=1}{K} \Phi_k$
					\STATE \quad $\omega \leftarrow \omega + \eta_\omega \frac{1}{K} \Sum{k=1}{K} \Psi_k$
					\ENDFOR
				\end{algorithmic}
			\end{algorithm}
		
	\subsection{Reinforcement learning setup}
	 
		In our framework, the observation space consists of two state variables: the current value of the Hawkes intensity process and the insurer's available cash reserves. To simulate the stochastic dynamics of the claim process and its intensity, we rely on Ogata's thinning algorithm \cite{Ogata81}. 
		Notably, the evolution of the intensity process is independent of the agent’s actions, and thus remains unaffected by the control policy. \\
		On the other hand, the agent directly influences the cash reserve through its actions. It is therefore essential to clearly define how the chosen policy impacts the surplus process.
		
		In the theoretical formulation of the problem, the exit time is random and may potentially never be reached. To address this issue in our numerical implementation, we introduce a maximum time horizon $T > 0$ and define the stopping time as:
		\begin{equation*}
			\tau = T \wedge \inf \{t_i \in \T, X_{t_i}^\pi < 0 \} .
		\end{equation*}
		Naturally, the introduction of $T$ modifies the original problem and may introduce a bias if not handled carefully. To mitigate this, we choose $T$ large enough so that, in the absence of any control intervention, ruin occurs before time $T$ with high probability. Formally, we select $T$ such that
		\begin{equation*}
			\mathbb{P} \left(\inf \{ t_i \in \T, X_{t_i}^\pi < 0 \} \geq T\right) \leq \varepsilon ,
		\end{equation*}
		where $\varepsilon > 0$. This ensures that the finite-horizon approximation remains faithful to the structure of the original problem.
			
		\subsubsection*{Naïve setup}
		
			We follow the MDP framework introduced in the previous section, where $\T$ denotes the discretised time grid, and $S_{t_i}$ represents the state at time $t_i \in \T$. In the most basic setup, we define the observation space as the pair $S_{t_i} = (X_{t_i}, \lambda_{t_i})$, and let the agent sample an action $A_{t_i}^\pi$ from a policy $\pi$, constrained to the interval $(-\infty, X_{t_i}]$. A positive action corresponds to a dividend payment, while a negative action corresponds to a capital injection. The cash reserve then evolves according to:
			\begin{equation*}
				X_{t_{i+1}}^\pi = X_{t_i}^\pi + h c - \Sum{k=1}{N_{t_{i+1}} - N_{t_i}} Y_k - A_{t_i}^\pi . 
			\end{equation*}
			The agent's expected reward under policy $\pi$ is then defined by:
			\begin{equation*}
				J(\pi) = \Expect{\Sum{j=1}{N(\tau)} e^{-\rho t_j} \left( A_{t_j}^\pi \Ind{A_{t_j}^\pi \geq 0} + \delta A_{t_j}^\pi \Ind{A_{t_j}^\pi < 0} \right)} .
			\end{equation*}
			While this approach is theoretically valid, it grants the agent considerable freedom, which can significantly slow down learning due to the difficulty of balancing exploration and exploitation. In particular, it becomes challenging for the agent to discover optimal intervention timings. For this reason, we propose a more structured approach that incorporates theoretical insights derived from the analytical study presented in the first part of the paper.
		
		\subsubsection*{Setup based on theoretical knowledge}
		
			This second approach restricts the admissible controls by imposing a two-barrier structure. 
			For capital injections, Proposition~\ref{PROP:capital_injection_threshold} provides an explicit optimal threshold~$\kappa^\star$. 
			For dividend payments, guided by the numerical solution of the HJB variational inequality, we postulate the existence of a state-dependent payout threshold $x^\star(y)$ for $y \in [b,+\infty)$. 
			Such a threshold is economically natural: once the surplus becomes sufficiently large, an optimal strategy must eventually prescribe dividend distributions.
			
			We define the observation space as $S_{t_i} = \lambda_{t_i}$ and use the policy to predict the values of the optimal boundaries $\kappa^\star(y)$ and $x^\star(y)$. The surplus process then evolves according to:
			\begin{equation}
				\label{EQN:RL_cash_dynamics}
				X_{t_{i+1}}^\pi = X_{t_i}^\pi + h c - \Sum{k=1}{N_{t_{i+1}} - N_{t_i}} Y_k - (X_{t_i}^\pi - x^\star(y)) \Ind{X_{t_i}^\pi \geq x^\star(y)} - X_{t_i}^\pi \Ind{\kappa^\star(y) \leq X_{t_i}^\pi < 0 }  .
			\end{equation}
			In this context, the agent's expected reward is:
			\begin{equation}
				\label{EQN:objective_with_barrier}
				J(\pi) = \Expect{\Sum{j=1}{N(\tau)} e^{-\rho t_j} \left( (X_{t_j}^\pi - x^\star(y)) \Ind{X_{t_j}^\pi \geq x^\star(y)} + \delta X_{t_j}^\pi \Ind{\kappa^\star(y) \leq X_{t_j}^\pi < 0 }  \right)} .
			\end{equation}
			This approach reduces the complexity of the learning task by restricting the agent’s output to the prediction of the two optimal boundaries, rather than a full-range action. As a result, it helps accelerate training and improves the stability of the learned policy.
					
	\subsection{Numerical results in the two-dimensional setting}
	
		\subsubsection{Learning boundaries}
			
			We implement both reinforcement learning algorithms together with standard regularisation techniques—such as entropy bonuses—to stabilise training and improve convergence. For comparability with the PDE-based results, we adopt the same model parameters as those reported in Table~\ref{TAB:model_parameters_PDE}.
			
			The learning procedure proceeds as follows.  
			Given a parameter vector $\theta \in \R^p$, we initialise a neural network policy~$\pi_\theta$ that takes as input the current value of the Hawkes intensity and outputs four real numbers corresponding to the parameters used to sample the control.  
			The network architecture consists of two hidden layers of $64$ neurons with ReLU activations.  
			Trajectory generation under the policy is carried out through the following steps:
			\begin{itemize}
				\item[i)] At each time step, the current intensity is observed and passed through $\pi_\theta$, which returns the parameters $(\mu_1,\sigma_1,\mu_2,\sigma_2)$.
				\item[ii)] From these parameters, we construct two normal distributions $\mathcal{N}(\mu_1,\sigma_1)$ and $\mathcal{N}(\mu_2,\sigma_2)$.
				\item[iii)] One sample is drawn from each distribution, and the log-probabilities of the sampled actions are recorded.
				\item[iv)] The corresponding control boundaries are constructed and applied to the surplus process according to Equation~\ref{EQN:RL_cash_dynamics}, after which steps~(i)--(iii) are repeated until the ruin time is reached.
			\end{itemize}
			Each simulated trajectory yields a total reward together with its associated sequence of log-probabilities.  
			Repeating this procedure $M$ times provides a Monte Carlo estimate of the policy gradient, using either Theorem~\ref{TH:gradient_classic_representation} or Theorem~\ref{TH:gradient_martingale_representation}.  
			The policy parameters are then updated via stochastic gradient ascent.
			In the actor–critic setting, the
			procedure remains identical except that a second neural network, with the same architecture as the policy network, is introduced to approximate the value function and serve as a learned baseline for variance reduction.
			
		\subsubsection{Comparison to PDE benchmark}
		
			To assess the performance of the reinforcement learning methods, we train the agents under the benchmark parameter set reported in Table~\ref{TAB:model_parameters_PDE} and compare the learned values to the reference solution obtained from the numerical resolution of the HJB variational inequality. 
			The time discretisation step is set to $h = 1/50$, and the time maximum horizon to $T = 50$, which corresponds to a maximum of $T/h = 2,500$ time steps, an upper limit that is never reached in practice due to earlier ruin (see Figure~\ref{FIG:survival_comparison} and Table~\ref{TAB:ruin_baseline}).
			Each policy update relies on $M = 2048$ Monte Carlo trajectories generated in parallel, with learning rates of order $10^{-3}$ for both the actor and the critic. Training is performed over $200$ epochs for each algorithm. 
			We consider two initial surplus–intensity states, $(x_0,y_0) = (1,2.8)$ and $(x_0,y_0) = (0,2.8)$, representing respectively a comfortably capitalised position and a near-boundary initial surplus.

\begin{figure}[H]
				\centering
				\begin{subfigure}[t]{0.49\textwidth}
					\centering
					\includegraphics[scale=0.38]{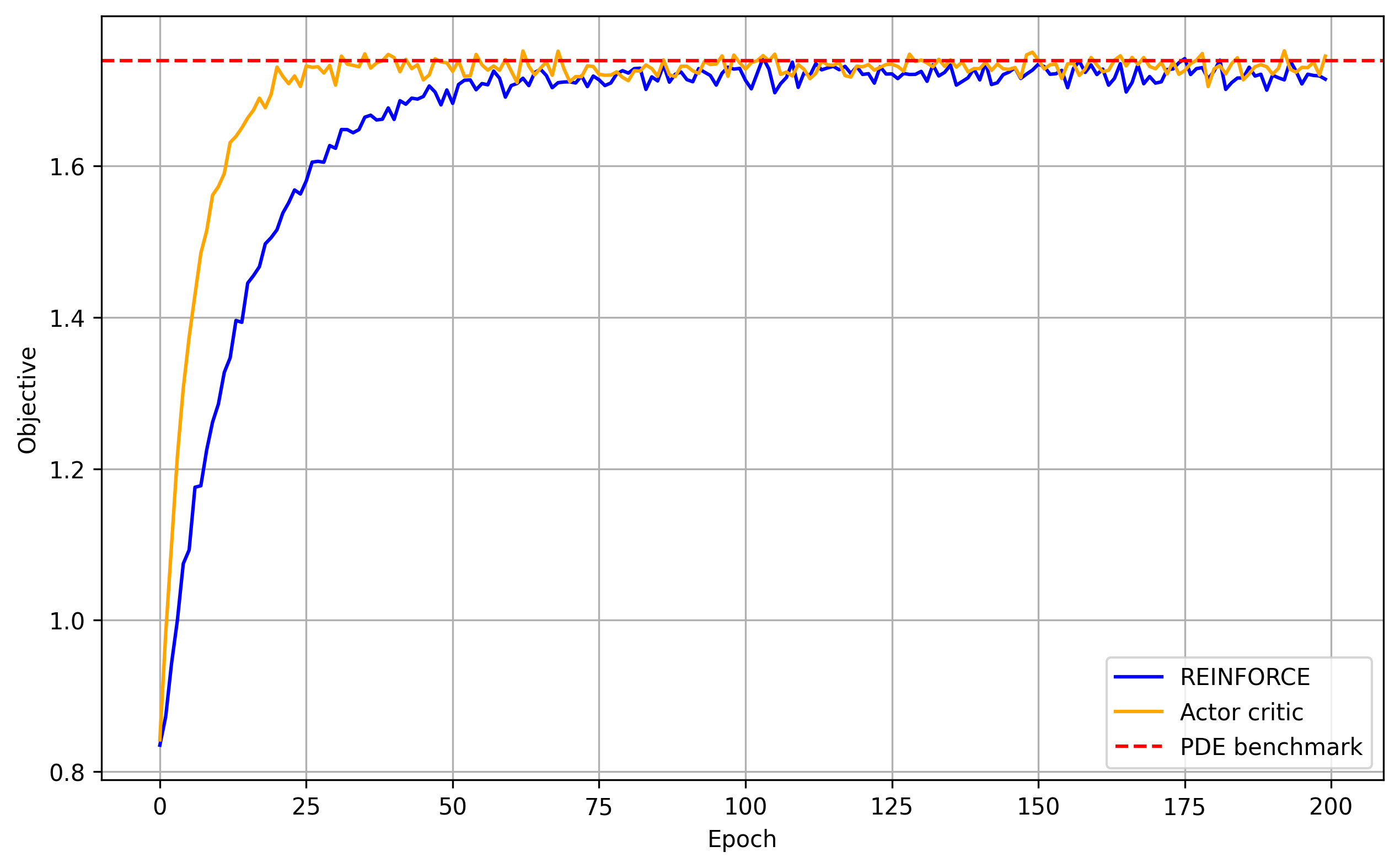}
					\caption{Training performance for $(x_0,y_0)=(1,2.8)$.}
					\label{FIG:learning_metrics_X0_1}
				\end{subfigure}
				\begin{subfigure}[t]{0.49\textwidth}
					\centering
					\includegraphics[scale=0.38]{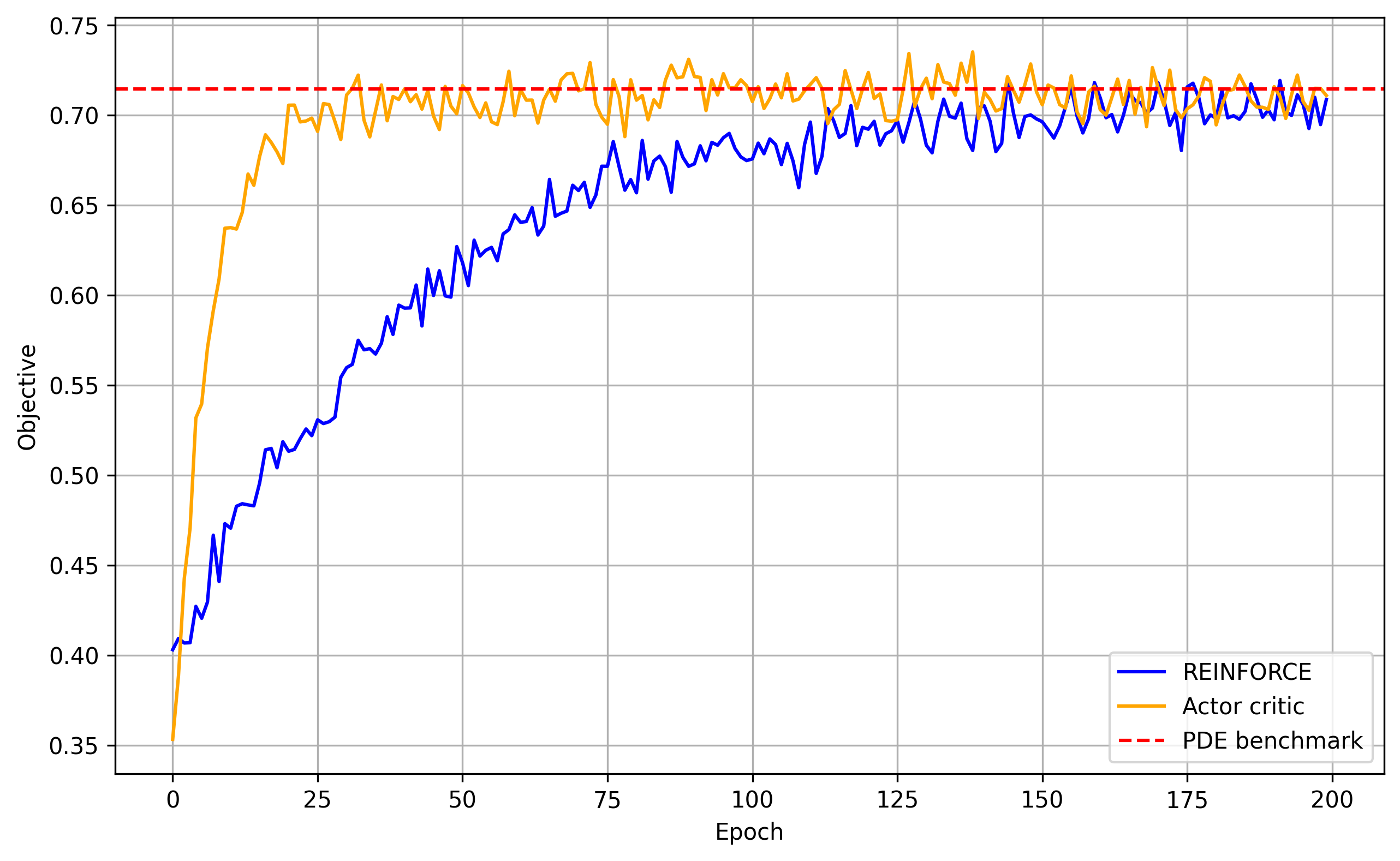}
					\caption{Training performance for $(x_0,y_0)=(0,2.8)$.}
					\label{FIG:learning_metrics_X0_0}
				\end{subfigure}
				\caption{Convergence of the learned objective toward the PDE benchmark value.}
			\end{figure}

\begin{figure}[H]
				\centering
				\begin{subfigure}[t]{0.49\textwidth}
					\centering
					\includegraphics[scale=0.4]{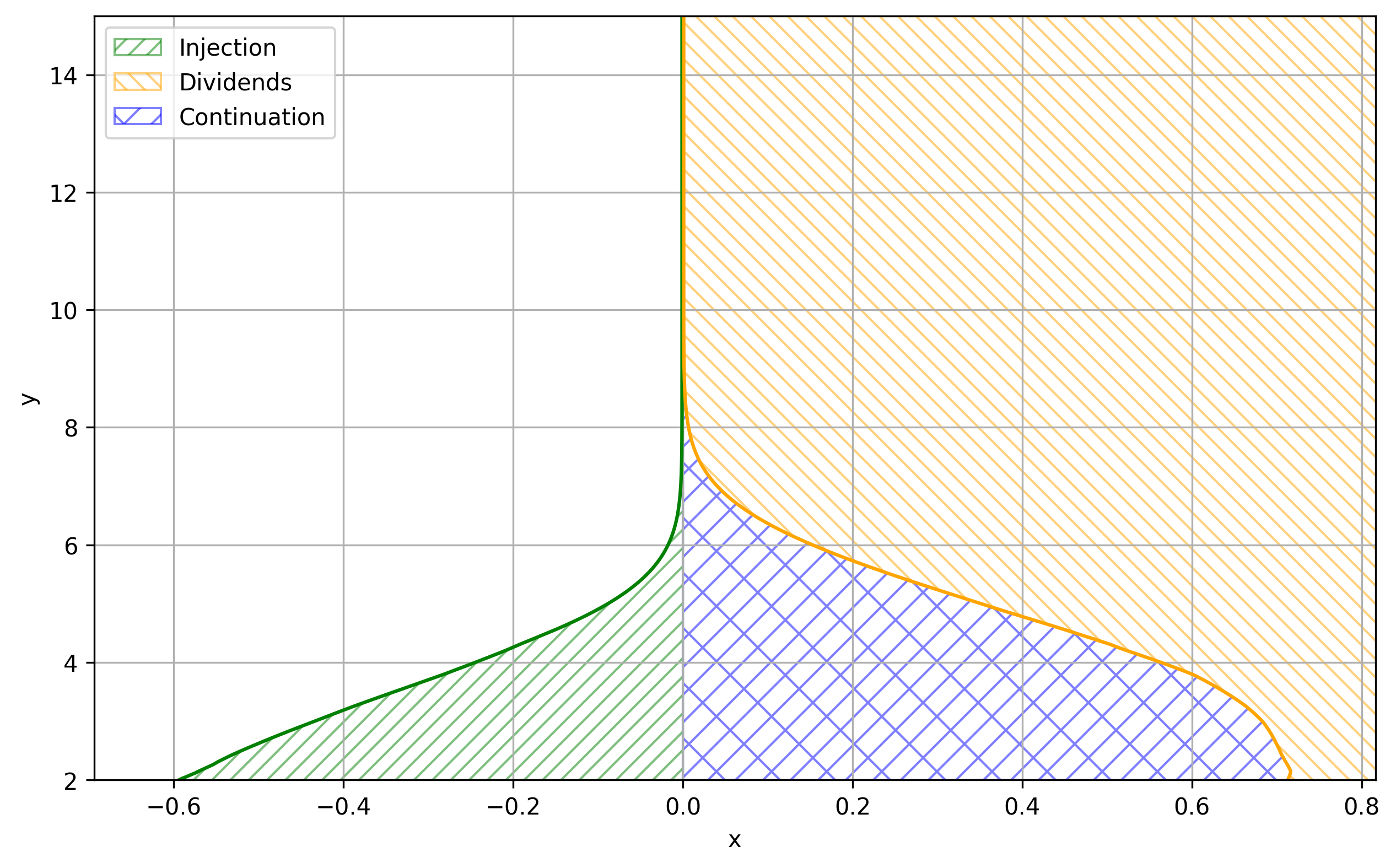}
					\caption{Learned optimal policy $(x_0, y_0) = (1, 2.8)$.}
					\label{FIG:policy_X0_1}
				\end{subfigure}
				\begin{subfigure}[t]{0.49\textwidth}
					\centering
					\includegraphics[scale=0.4]{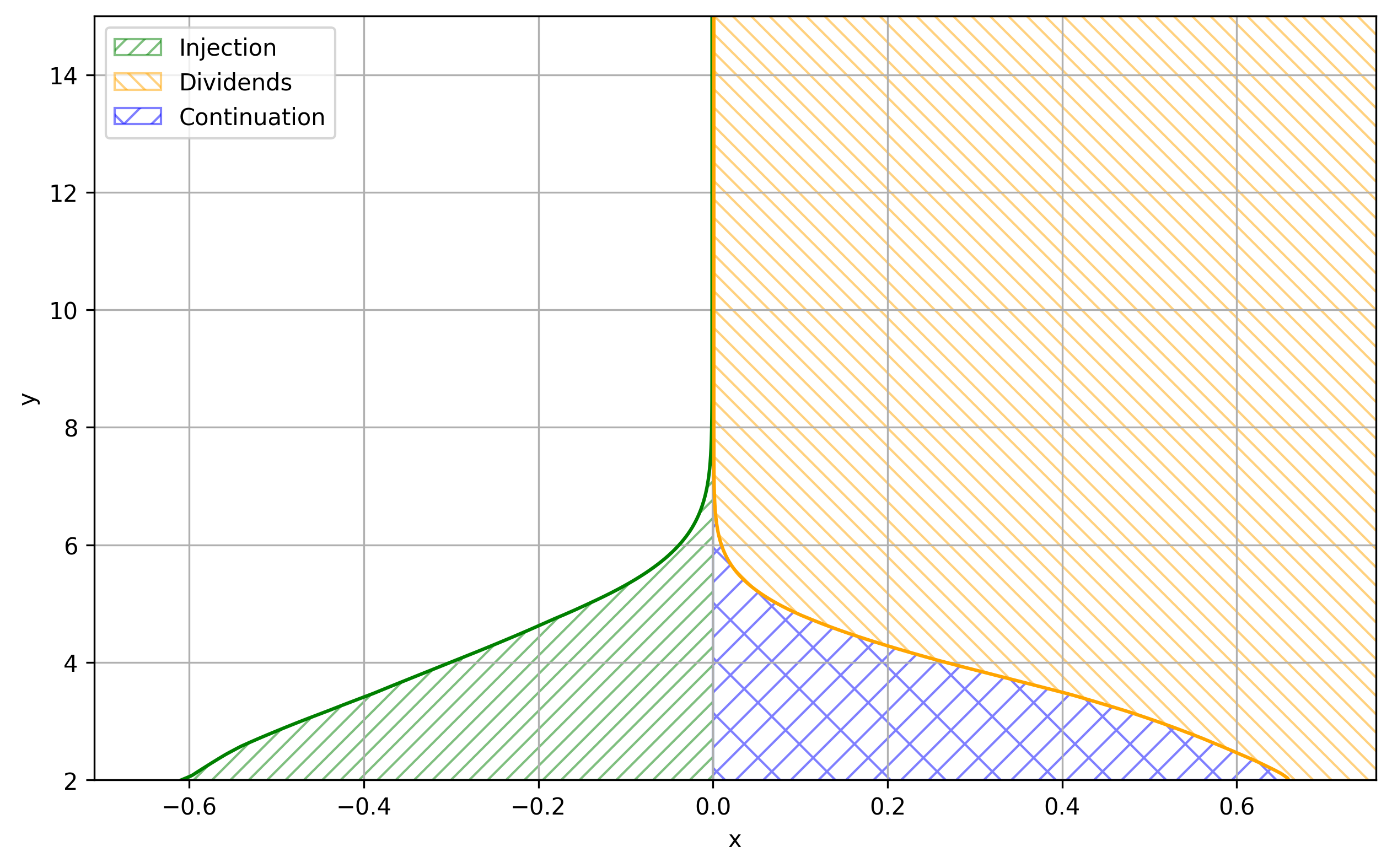}
					\caption{Learned optimal policy $(x_0, y_0) = (0, 2.8)$.}
					\label{FIG:policy_X0_0}
				\end{subfigure}
				\caption{Learned control regions obtained by the policy gradient algorithm for two initial states.}
			\end{figure}

			Figures~\ref{FIG:learning_metrics_X0_1}--\ref{FIG:learning_metrics_X0_0} display the evolution of the empirical objective $J$ during training, together with the PDE benchmark value.  
			In both initial configurations, the actor–critic method exhibits the fastest and most stable convergence, reaching the PDE benchmark within relatively few epochs.  
			The REINFORCE estimator also converges toward the correct value, although with slightly higher variance, which is expected for Monte Carlo policy gradients.  
			Note that the policy gradient method optimises over a restricted class of parametrised stochastic policies, which structurally yields a lower bound on the true value. The residual fluctuations around the PDE benchmark are due to the exploration variance combined with the Monte Carlo approximation of the objective.
			The actor--critic method exhibits faster stabilisation and lower epoch-to-epoch variance owing to the learned critic baseline, while REINFORCE, though unbiased, displays higher variance due to its reliance on full trajectory returns.
			Overall, both algorithms succeed in learning policies whose performance matches the PDE solution, thereby validating the discrete-time formulation and the policy gradient estimators developed in this section.
			In Figures~\ref{FIG:policy_X0_1} and~\ref{FIG:policy_X0_0}, we display the control regions learned by the reinforcement learning agent.  
			The colour map indicates the action selected in each state: the blue region corresponds to inaction, the green region to capital injection, and the yellow region to dividend distribution.
			The learned strategies display the same qualitative structure as the optimal policy obtained through the variational inequality formulation in the PDE section. This close agreement provides strong evidence for the validity of the reinforcement learning approach. 
			Some discrepancies between the two training runs can be observed in the precise location of the control boundaries. This behaviour is expected: since policy-gradient methods optimise over a restricted class of parametrised stochastic policies, they converge to near-optimal strategies rather than an exact optimum. For such quasi-optimal policies, the control boundary is not uniquely defined. Our Monte Carlo experiments confirm that the estimated value is only slightly sensitive to variations in the dividend boundary, provided that the global structure of the optimal strategy is preserved.
			\begin{table}[H]
				\centering
				\begin{tabular}{cc|c||ccc||ccc}
					\toprule
					$x$ & $y$ & PDE & MC (Opt.) & IC95\% (MC Opt.) & Rel. err. & MC (RL) & IC95\% (RL) & Rel. err. \\
					\midrule
					\midrule
					$0$ & $2$ & $0.8588$ & $0.8414$ & [$0.8023$, $0.8805$] & $-2.07\%$ & $0.8677$ & [$0.8263$, $0.9090$] & $1.03\%$ \\
					$0$ & $3$ & $0.6811$ & $0.6642$ & [$0.6269$, $0.7014$] & $-2.54\%$ & $0.6840$ & [$0.6455$, $0.7225$] & $0.44\%$ \\
					$0$ & $4$ & $0.5298$ & $0.5181$ & [$0.4833$, $0.5528$] & $-2.27\%$ & $0.5360$ & [$0.5016$, $0.5705$] & $1.17\%$ \\
					\midrule
					\midrule
					$0.5$ & $2$ & $1.3874$ & $1.3412$ & [$1.2987$, $1.3838$] & $-3.44\%$ & $1.3890$ & [$1.3467$, $1.4313$] & $0.12\%$ \\
					$0.5$ & $3$ & $1.2031$ & $1.1581$ & [$1.1166$, $1.1995$] & $-3.89\%$ & $1.2368$ & [$1.1948$, $1.2788$] & $2.80\%$ \\
					$0.5$ & $4$ & $1.0360$ & $0.9882$ & [$0.9514$, $1.0249$] & $-4.84\%$ & $1.0324$ & [$0.9920$, $1.0728$] & $-0.35\%$ \\
					\midrule
					\midrule
					$1.0$ & $2$ & $1.8881$ & $1.8673$ & [$1.8257$, $1.9089$] & $-1.11\%$ & $1.8872$ & [$1.8451$, $1.9293$] & $-0.05\%$ \\
					$1.0$ & $3$ & $1.7033$ & $1.6886$ & [$1.6477$, $1.7294$] & $-0.87\%$ & $1.7143$ & [$1.6727$, $1.7558$] & $0.64\%$ \\
					$1.0$ & $4$ & $1.5360$ & $1.4894$ & [$1.4527$, $1.5261$] & $-3.13\%$ & $1.5312$ & [$1.4911$, $1.5714$] & $-0.31\%$ \\
					\bottomrule
				\end{tabular}
				\caption{Comparison of the PDE and RL estimates of the value function.}
				\label{TAB:value_PDE_RL_comparison}
			\end{table}
			Finally, Table~\ref{TAB:value_PDE_RL_comparison} reports three sets of values for representative state pairs. The first column (PDE) shows the benchmark value computed from the numerical solution of the HJB variational inequality. 
			The second block provides a Monte Carlo estimate of the value obtained when applying the theoretically optimal policy to the discretised environment. 
			The third block reports the corresponding estimate obtained using the policy learned by reinforcement learning. Both Monte Carlo values are computed from $4{,}096$ simulated trajectories, and the reported confidence intervals are the standard asymptotic $95\%$ confidence intervals.
			The relative errors reported in the table are computed by comparing respectively the Monte Carlo estimate of the value function applying theoretical optimal policy and RL Monte Carlo estimate to the PDE benchmark value.
			The policy learned by RL consistently outperforms the theoretically optimal continuous-time policy when both are evaluated on the discretised environment. This result is empirical, yet expected by construction. The RL policy is trained directly on the discretised environment and therefore implicitly adapts to the time discretisation, while the continuous-time optimal policy targets the original continuous-time problem.

			Beyond this qualitative agreement, the reinforcement learning framework offers two significant advantages. First, it scales naturally to higher-dimensional settings in which PDE-based methods become impractical or computationally prohibitive. Second, it provides a flexible modelling environment: changes to the claim distribution, richer dependence structures between claims and intensity, or more complex interactions in the dynamics can be incorporated with minimal modifications to the learning procedure. In this regard, reinforcement learning constitutes a powerful and adaptable tool for approximating optimal strategies in stochastic control problems with complex or high-dimensional dynamics.
			
	\subsection{Scaling to higher-dimensional portfolios}
	\label{SEC:multiline_extension}
		
		The previous numerical experiments demonstrate that the reinforcement learning methodology successfully recovers the PDE benchmark in the two-dimensional setting.
		However, the practical advantage of policy-gradient methods over grid-based solvers becomes most apparent in higher-dimensional problems, where the curse of dimensionality renders finite-difference schemes computationally intractable.
		To illustrate this, we extend the model to an insurer operating $d$ distinct lines of business with correlated claim dynamics.
			
		\subsubsection{Multi-line surplus model}
			
			Consider an insurer managing $d \geq 2$ insurance lines. Each line $i \in \{1,\dots,d\}$ collects premiums at a constant rate $c^i > 0$ and is exposed to claims arriving according to a counting process $N^i = (N^i_t)_{t \geq 0}$.
			The aggregate surplus process, in the absence of any control, is given by
			\begin{equation}
				\label{EQN:multiline_surplus}
				R_t = x + \left(\sum_{i=1}^{d} c^i\right) t - \sum_{i=1}^{d} \sum_{k=1}^{N^i_t} Y^i_k,
			\end{equation}
			where $x \in \mathbb{R}^+$ is the initial capital and $(Y^i_k)_{k \geq 1}$ are i.i.d.\ positive claim sizes associated with line $i$, with density $f^i$ and independent of the counting processes.
			
			To capture contagion effects across insurance lines, we model the vector of claim intensities $\lambda_t = (\lambda^1_t, \dots, \lambda^d_t)$ as a mutually exciting Hawkes process.
			Each intensity $\lambda^i$ evolves according to
			\begin{equation}
				\label{EQN:multiline_hawkes}
				d\lambda^i_t = a^i(b^i - \lambda^i_t)\,dt + \sum_{j=1}^{d} \eta^{ij}\,dN^j_t, \qquad i = 1,\dots,d,
			\end{equation}
			where, for each line $i$:
			\begin{itemize}
				\item $a^i > 0$ is the mean-reversion rate toward the baseline intensity $b^i > 0$,
				\item $\eta^{ii} > 0$ captures the self-excitation effect within line $i$,
				\item $\eta^{ij} \geq 0$, for $j \neq i$, captures the cross-excitation from line $j$ to line $i$.
			\end{itemize}
			The matrix $\eta = (\eta^{ij})_{1 \leq i,j \leq d}$ encodes the full dependency structure of the claim arrivals.
			Off-diagonal entries $\eta^{ij} > 0$ reflect contagion channels through which a claim in one business line increases the likelihood of subsequent claims in another line.
			Such cross-excitation naturally arises in multi-line portfolios exposed to common risk factors: for instance, a natural disaster may simultaneously trigger property, business-interruption, and liability claims.

			As in the two-dimensional setting, the insurer distributes dividends and injects capital.
			Under a control strategy $\alpha = (Z_t, K_t)_{t \geq 0}$, the controlled surplus process reads
			\begin{equation*}
				X_t = x + \left(\sum_{i=1}^{d} c^i\right) t - \sum_{i=1}^{d} \sum_{k=1}^{N^i_t} Y^i_k - Z_t + K_t,
			\end{equation*}
			and the ruin time is defined as before by $T^\alpha = \inf\{t \geq 0 : X_{t^+} < 0\}$.
			The shareholders' objective is
			\begin{equation}
				\label{EQN:multiline_control_problem}
				v(x, y) = \sup_{\alpha \in \mathcal{A}(x,y)} \mathbb{E}\left[\int_0^{T^\alpha} e^{-\rho s}\,dZ_s - \delta \int_0^{T^\alpha} e^{-\rho s}\,dK_s\right], \qquad (x,y) \in \mathbb{R} \times \prod_{i=1}^d [b^i,+\infty),
			\end{equation}
			where $y = (y^1,\dots,y^d)$ denotes the initial intensity vector.
			
			The state space of this problem is $(d+1)$-dimensional: the scalar surplus $x$ and the $d$ intensity components.
			Already for $d = 3$, the state space has dimension four, making finite-difference schemes impractical due to the exponential growth in the number of grid points.
			
		\subsubsection{Adaptation of the learning procedure}
			
			The reinforcement learning setup extends naturally from the single-line case.
			The observation space becomes $S_{t_i} = (\lambda^1_{t_i}, \dots, \lambda^d_{t_i})$, and the policy $\pi_\theta$ maps the current intensity vector to the parameters of the control boundaries.
			Specifically, the neural network takes $\lambda_{t_i} \in \mathbb{R}^d$ as input and outputs $4$ real numbers parametrising the distributions from which the dividend barrier $x^\star$ and the injection threshold $\kappa^\star$ are sampled.
			
			We emphasise that the only modifications required compared to the two-dimensional implementation are the event simulation, which now involves $d$ interacting counting processes, and the input dimension of the neural network, which increases from $1$ to $d$. The output dimension, the learning algorithms (Algorithms~\ref{ALG:off_line_policy_gradient} and~\ref{ALG:actor_critic_gradient_desc}), and the policy gradient estimators all remain unchanged, illustrating the scalability of the approach.
			
		\subsubsection{Numerical experiments}
		\label{SEC:multiline_numerics}
			
			\paragraph*{Parameter configuration} ~\\
			We consider $d = 3$ insurance lines, corresponding to a state space of dimension $4$.
			The model parameters are chosen to reflect a diversified portfolio with heterogeneous risk profiles and moderate cross-excitation.
			The configuration is reported in Tables~\ref{TAB:multiline_model_params} and~\ref{TAB:multiline_excitation_matrix}.
			
			\begin{table}[H]
				\centering
				\begin{tabular}{c|ccc}
					\toprule
					Line $i$ & $c^i$ & $a^i$ & $b^i$ \\
					\midrule
					$1$ & $0.40$ & $2.0$ & $0.8$ \\
					$2$ & $0.35$ & $1.5$ & $0.7$ \\
					$3$ & $0.30$ & $2.5$ & $0.5$ \\
					\bottomrule
				\end{tabular}
				\qquad
				\begin{tabular}{c|c}
					\toprule
					Line $i$ & $\beta^i$ \\
					\midrule
					$1$ & $3.0$ \\
					$2$ & $2.5$ \\
					$3$ & $4.0$ \\
					\bottomrule
				\end{tabular}
				\caption{Model parameters for the multi-line extension. The discount rate $\rho = 0.1$ and injection cost $\delta = 1.8$ are common to all lines.}
				\label{TAB:multiline_model_params}
			\end{table}
			
			\begin{table}[H]
				\centering
				$\eta = \begin{pmatrix} 0.30 & 0.10 & 0.05 \\ 0.10 & 0.35 & 0.08 \\ 0.05 & 0.08 & 0.25 \end{pmatrix}$
				\caption{Cross-excitation matrix $\eta$.}
				\label{TAB:multiline_excitation_matrix}
			\end{table}
			
			\paragraph*{Results} ~\\
			In the absence of a PDE reference solution, we compare the learned strategies against an optimised constant-barrier heuristic.
			This heuristic applies a state-independent dividend barrier~$\bar x$ and a state-independent injection threshold~$\bar\kappa$.
			The pair $(\bar x, \bar\kappa)$ is selected by a two-stage grid search over the parameter space.
			This procedure yields the best constant-barrier strategy attainable without any dependence on the intensity process, and therefore provides a meaningful lower bound against which the RL policies can be assessed.
			
			Both the REINFORCE and actor--critic algorithms are trained using the same hyperparameters as in the two-dimensional experiments: a time step $h = 1/50$, a horizon $T = 50$, and $M = 2048$ parallel trajectories per epoch.
			The neural network architecture is identical to the single-line case, with the only modification being the input dimension, which increases from $1$ to~$d$. Training is carried out over $300$ epochs.
			\begin{figure}[H]
				\centering
				\includegraphics[scale=0.5]{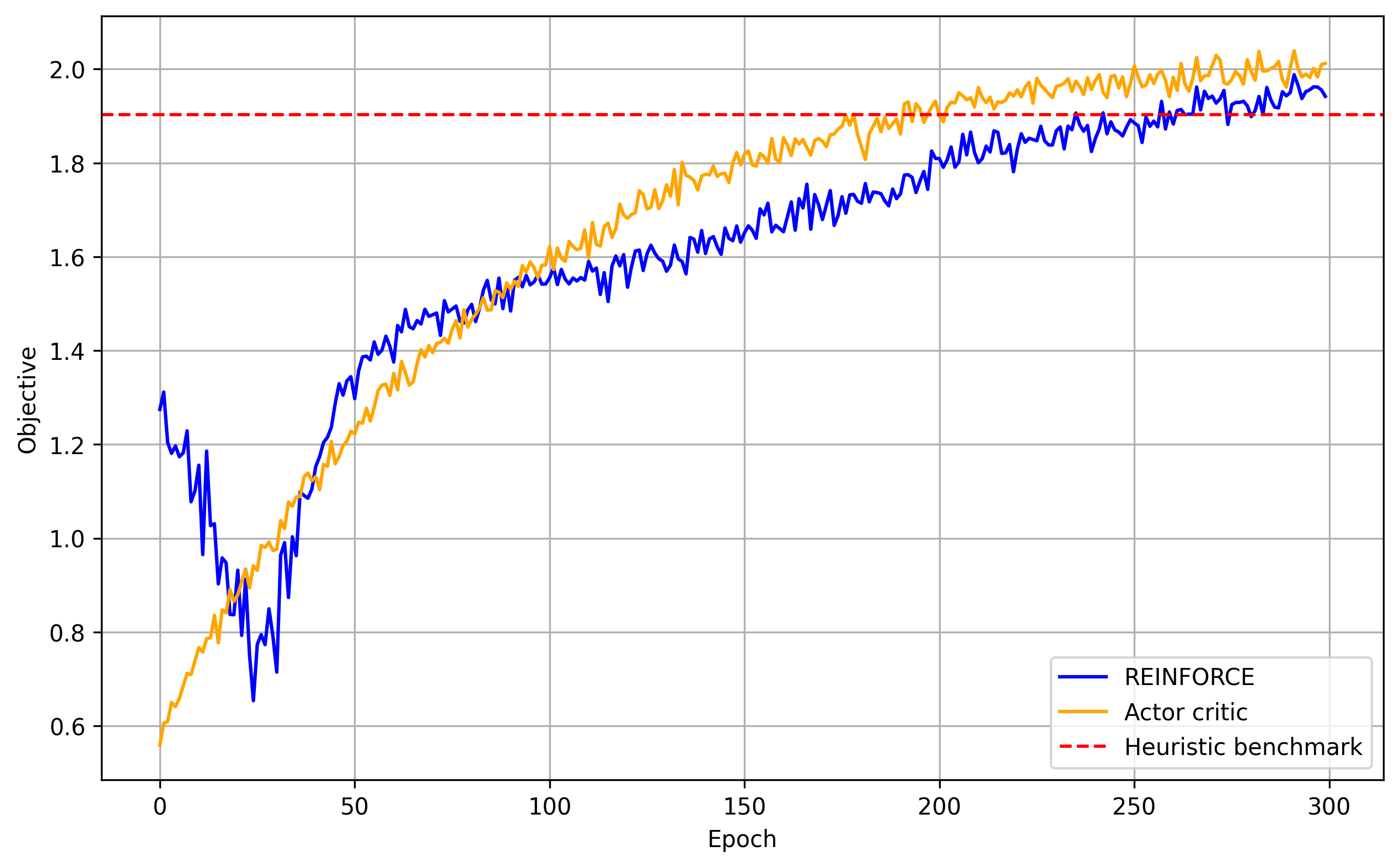}
				\caption{Comparison of the learned objective and the heuristic reference value. \\ $x_0 = 1$ and $y^i = b^i$,  $\forall i$.}
				\label{FIG:learning_metric_multidim}
			\end{figure}
			Figure~\ref{FIG:learning_metric_multidim} displays the evolution of the empirical objective during training, together with the constant-barrier heuristic as a reference. 
			Both algorithms cross the heuristic threshold around epoch $200$ and continue improving beyond it. The actor--critic method converges more smoothly from the start, while REINFORCE exhibits higher initial variance before stabilising at a comparable level. 
			Both algorithms reach similar final performance, slightly above the heuristic, validating the reinforcement learning approach in this four-dimensional setting.
			
			In this work, the reinforcement learning methodology is used as a numerical solver with fully specified model dynamics and simulation-based trajectory generation. 
			An alternative perspective, developed by Jia and Zhou~\cite{JiaZhou22_AC, JiaZhou22_TD, JiaZhou23_QL} in the context of diffusion-driven control problems, is to learn optimal policies directly from observed trajectories without knowledge of the model coefficients. 
			The recent work of Liang, Luo and Yu~\cite{LiangLuoYu2025} has extended this framework to a class of singular stochastic control problems, characterising the optimal control as a pair of regions via policy improvement on region iterations. 
			Extending this model-free approach to our setting would be particularly appealing, as it would allow learning optimal dividend and injection strategies directly from observed claim histories, without requiring knowledge of the Hawkes parameters. 
			The main challenges specific to our model, namely the jump dynamics and the need to infer the latent intensity from claim arrival data alone, make this a natural and promising direction for future research.
			
	\printbibliography[heading=bibintoc, title={Bibliography}]

\section*{Statements and declarations}
	
	\subsection*{Funding}
	The authors declare that no funds, grants, or other support were received during the preparation of this manuscript.
	
	\subsection*{Competing interests}
	The authors have no relevant financial or non-financial interests to disclose.
	
	\subsection*{Author contributions}
	All authors listed contributed equally to this work.

\end{document}